\documentclass[11pt]{article}
\usepackage[T1]{fontenc}
\usepackage{lmodern,amsmath,amssymb,amsthm,mathtools,mathrsfs}
\usepackage[margin=1in]{geometry}
\usepackage{enumitem,microtype,xcolor,booktabs}
\usepackage{tikz}
\usetikzlibrary{arrows.meta,calc}
\definecolor{xiBlue}{RGB}{211,235,243}
\definecolor{xiRose}{RGB}{247,214,223}
\definecolor{xiGold}{RGB}{255,231,174}
\definecolor{xiInk}{RGB}{20,35,48}
\usepackage[colorlinks=true,linkcolor=blue,citecolor=blue,urlcolor=blue]{hyperref}
\hypersetup{pdftitle={The monochromatic three-arm exponent for critical planar percolation},
 pdfauthor={Gefei Cai},pdfsubject={CLE, Liouville quantum gravity, and conformal welding}}
\numberwithin{equation}{section}
\newtheorem{theorem}{Theorem}[section]
\newtheorem{proposition}[theorem]{Proposition}
\newtheorem{lemma}[theorem]{Lemma}

\newtheorem{definition}[theorem]{Definition}
\theoremstyle{remark}
\newtheorem{remark}[theorem]{Remark}

\renewcommand{\P}{\mathbb P}

\newcommand{\CLE}{\mathrm{CLE}}
\newcommand{\SLE}{\mathrm{SLE}}

\newcommand{\diam}{\operatorname{diam}}

\allowdisplaybreaks
\title{Monochromatic arm exponents for critical planar percolation}
\author{Gefei Cai}
\date{\today}
\begin{document}
\maketitle

\begin{abstract}
We determine the monochromatic three-arm exponent for critical planar percolation through the exact spectral threshold of an explicit positive integral operator family with elementary kernels, obtaining $\xi_3=0.76928375\cdots$. Unlike the monochromatic two-arm (backbone) event, the three-arm event can be obstructed by two distinct white clusters, so a Markovian exploration must retain the past connections. We introduce a two-point boundary observable that records these connections. Using the coupling between conformal loop ensemble (CLE) and Liouville quantum gravity (LQG), this observable satisfies a renewal equation in terms of quantum boundary lengths, and yields the explicit operator whose spectral information determines $\xi_3$. Our approach also gives an explicit spectral characterization of monochromatic $n$-arm exponents for all $n\ge4$.
\end{abstract}

\setcounter{tocdepth}{1} \tableofcontents

\section{Introduction}\label{sec:introduction}

Bernoulli percolation is one of the simplest and most fundamental models of statistical mechanics. In two dimensions, it has been central to the study of critical phenomena through conformal invariance. For critical site percolation on the triangular lattice, Smirnov~\cite{Smirnov} rigorously proved the conformal invariance of crossing probabilities and Cardy's formula. Together with the introduction of Schramm--Loewner evolution (SLE)~\cite{SchrammSLE}, this establishes the convergence of the critical percolation interface to $\SLE_6$. The full collection of interfaces was subsequently shown to have a scaling limit~\cite{CNFull}, described by the nested conformal loop ensemble (CLE) $\CLE_6$~\cite{Sheffield09}.

The arm exponents are a family of basic scaling exponents in critical percolation that describe the power-law decay rate of the probabilities of disjoint paths crossing annuli with specific colors. An arm event is called \emph{polychromatic} if its prescribed colors include both black and white, and \emph{monochromatic} if all arms have the same color. The exact evaluation of these exponents was an early application of SLE. The one-arm exponent is $5/48$~\cite{LSWOneArm}, and the polychromatic $(n\ge2)$-arm exponents are $\alpha_n=\frac{n^2-1}{12}$
for every color sequence containing both colors~\cite{SWArms}. These results confirmed predictions from statistical physics~\cite{ADA}.

On the other hand, for the monochromatic $n$-arm exponents with $n\ge2$, the exact evaluation remained open for many years. In~\cite{BN}, it is shown that $\xi_n$ exists and satisfies $\alpha_n<\xi_n<\alpha_{n+1}$ for $n\ge2$. It is only recently~\cite{NQSZ} that the monochromatic two-arm exponent $\xi_2=0.3566668\cdots$ (also known as the backbone exponent) has been determined explicitly as a root of an elementary function. Since the backbone exponent is determined by the obstruction caused by a single white cluster, the derivation there reduces to a conformal-radius computation for SLE bubbles, which is carried out using the coupling between SLE and Liouville quantum gravity (LQG).

However, for three or more monochromatic arms, distinct white clusters may jointly obstruct the crossings even when none of them does so alone. In terms of the exploration process in~\cite{NQSZ}, the previously discovered loops become boundary arcs, and a newly discovered loop can separate the target together with one of these arcs. Therefore, the stopping rule depends on past connections, which makes essential difficulties in extending their derivation to general $\xi_n$; see~\cite[Remark~2.3]{NQSZ}.

In this paper, we determine the monochromatic three-arm exponent $\xi_3$ through the unique spectral threshold of an explicit family of positive integral operators with elementary kernels, obtaining $\xi_3=0.76928375\cdots$. 
The key is a new two-point boundary observable that records the connections required for the monochromatic three-arm event. After coupling CLE with an independent LQG field, we use the domain Markov property to derive a renewal equation for this observable in terms of quantum boundary lengths. Taking the Mellin transform then yields an integro-differential equation with specific boundary conditions. Finally, a general spectral threshold theorem concludes the spectral characterization of $\xi_3$. Our approach is robust and also gives an explicit spectral characterization of monochromatic $n$-arm exponents for all $n\ge4$; see Section~\ref{sec:higher}.

\subsection{The exponent $\xi_3$ and its analytic characterization}\label{subsec:intro-critical}

Let $\delta\mathbb T$ be the triangular lattice with mesh size $\delta>0$, with each vertex independently colored black or white with probability $1/2$ each. For $0<r<R$, write $p_3(r,R)$ for the limit of the probability of three vertex-disjoint black crossings of $\{z\in\mathbb C:r<|z|<R\}$ as $\delta\downarrow0$. Let $\tau=(2\pi)^{-1}\log(R/r)>0$ be the conformal modulus, and set $p_3(\tau)=p_3(e^{-2\pi\tau},1)$. The results of~\cite{BN} give
\begin{equation}\label{eq:known-exponent}
 p_3(\tau)=\exp\{-2\pi\xi_3\tau+o(\tau)\},
 \qquad \frac23<\xi_3<\frac54.
\end{equation}

We express $\xi_3$ through the following family of integral operators on $\mathbb R$. For $t\in\mathbb R$, put
\begin{equation}\label{eq:intro-kernel}
 \begin{gathered}
 x(t)=\frac{e^t}{1+e^t},\qquad a(t)=1-x(t),\qquad
 L_0=-\frac{d^2}{dt^2}+\frac14,\\
 H(t)=\frac{9\sqrt3}{16\pi}
       \frac{2\cosh(t/2)+3}{(2\cosh(t/2)+1)^3}.
 \end{gathered}
\end{equation}
For $s\in\mathbb C$ with $|\operatorname{Re}s|<2$ and a bounded measurable function $f:\mathbb R\to\mathbb C$, define
\begin{equation}\label{eq:intro-T}
 (T_s f)(t)=\int_{\mathbb R}H(t-u)
 \bigl[x(t)^{1+s}x(u)^{1-s}
       +a(t)^{1+s}a(u)^{1-s}\bigr]f(u)\,du.
\end{equation}
The powers use the real logarithms of $x$ and $a$.

\begin{theorem}\label{thm:main}
There is a unique positive parameter $s_*\in(0,2)$ for which
\begin{equation}\label{eq:main-resonance}
                        L_0r=T_{s_*}r
\end{equation}
has a strictly positive bounded even solution $r\in C^2(\mathbb R)$. Moreover, $3/2<s_*<2$, and
\begin{equation}\label{eq:main-value}
             \displaystyle\xi_3=\frac{s_*^2-1/4}{3}.
\end{equation}
The solution is unique up to multiplication by a positive constant and has a positive common limit at both $-\infty$ and $+\infty$.
\end{theorem}

All coefficients and the boundary condition in \eqref{eq:main-resonance} are explicit.  Numerical evaluation gives
\[
                 s_*\approx1.59932838,\qquad \xi_3\approx0.76928375.
\]
See Appendix~\ref{subsec:critical-numerics}. This agrees with the estimate $0.77\pm0.02$ obtained by the transfer matrix~\cite{JZJ}.
In a recent independent numerical study~\cite{WDL26}, the authors obtained the estimate $\xi_3=0.7692(3)$ from the fractal dimension of triconnected clusters in critical site percolation on the triangular lattice. Their numerical result is in excellent agreement with our exact formula for $\xi_3$.

The parameter $s_*$ is a simple zero of the Fredholm determinant defined as follows. Put
\[
  g(t)=e^{-|t|/2},\qquad
  q_0(t)=\frac{3\sqrt3}{8\pi(\cosh(t/2)+1/2)},\qquad
 A(t)=\frac{27}{32\pi^2}t\coth\frac{3t}{4},
\]
where $A(0)=9/(8\pi^2)$ by continuity. Write $(f*h)(t)=\int_{\mathbb R}f(t-u)h(u)\,du$ for convolution, and set
\[
 R_s f=g*(T_sf)-q_0*f,\qquad
 Sf=f-A*f+A\int_{\mathbb R}f.
\]
The latter formula is applied to exponentially decaying even functions. Let $L^2_{\mathrm{even}}(\mathbb R)$ be the space of square-integrable even functions, let $I$ be its identity operator, and let $1$ denote the constant function.

\begin{proposition}\label{prop:intro-determinant}
The block operator below is the identity plus a trace-class operator on $\mathbb C\oplus L^2_{\mathrm{even}}(\mathbb R)$, and its Fredholm determinant
\begin{equation}\label{eq:intro-determinant}
 \mathscr D(s)=\det
 \begin{pmatrix}
  \displaystyle\int_{\mathbb R}R_s1 &
       \displaystyle f\longmapsto\int_{\mathbb R}R_sf\\[2pt]
  -SR_s1&I-SR_s
 \end{pmatrix},\qquad |\operatorname{Re}s|<2,
\end{equation}
is holomorphic in this strip. The parameter $s_*$ of Theorem~\ref{thm:main} is its smallest positive zero:
\[
 s_* =\min\{s\in(0,2):\mathscr D(s)=0\}.
\]
This zero is simple.
\end{proposition}

The proof of Proposition~\ref{prop:intro-determinant} is given in Appendix~\ref{sec:determinant}.

\subsection{Proof strategy}\label{subsec:intro-strategy}

We first explain why the approach used to compute the backbone exponent in~\cite{NQSZ} does not directly extend to three arms. For the monochromatic two-arm (backbone) event, the obstruction comes from a single white cluster, whose outer boundary surrounds the origin. In~\cite{NQSZ}, the authors used a radial exploration of CLE via a Poisson point process of SLE bubbles, stopping at the first bubble whose outer boundary surrounds the origin; in particular, such a stopping time is independent of the previously discovered CLE loops. Then the backbone exponent is reduced to the moment of the conformal radius of SLE bubbles, which can be computed using the coupling between CLE and LQG. For three arms, however, two distinct white clusters may jointly separate the origin from the boundary, even when neither cluster separates them individually. When doing the same exploration of CLE, previously discovered CLE loops become boundary arcs after conformal mapping, and a newly discovered loop may separate the marked point together with one of these boundary arcs. The stopping time in the monochromatic three-arm case therefore depends on the previously discovered loops, and one would have to keep track of the boundary arcs inherited from the full exploration history. This difficulty was pointed out in~\cite[Remark~2.3]{NQSZ}.

Our approach circumvents this difficulty by introducing a two-point boundary observable $G$ that records the corresponding connections, including separating configurations formed by a white cluster together with an exterior boundary arc; see Section~\ref{sec:geometry}. The key geometric observation is that, after cutting a boundary-touching CLE loop avoiding the inner set, the same observable appears in the remaining region with the new boundary marks given by the CLE loop. More precisely, we identify the observable in the remaining region with the indicator that the monochromatic three-arm event holds and the discarded boundary arc is maximal; see Lemma~\ref{geo:restart}. Averaging over the original boundary marks then expresses the averaged two-point boundary observable as a weighted sum of the same observable in the remaining components, together with a contribution from white clusters meeting the inner set. This is the formula of Proposition~\ref{prop:geometric-cutting}.

To evaluate this sum, we need a tractable description of the retained component together with its inherited boundary arc and the newly exposed arc along the CLE interface. We couple CLE with an independent LQG field and equip both arcs with their quantum boundary lengths.
The coupling of SLE and CLE with LQG was originally developed in~\cite{Sheffield16,DMS21,MSW20,MSW21} and has since been widely used to derive exact formulas in statistical mechanics; see e.g.~\cite{AHSInt,ACSW,SXZ,LSYZ,CFSX}. In Section~\ref{sec:lqg}, we use this coupling to identify the quantum surface associated with the retained component as a CLE-decorated quantum disk with one interior mark and two boundary marks. It also gives an explicit kernel for the quantum boundary lengths of the newly exposed and the discarded boundary arcs. Integrating the sum in Proposition~\ref{prop:geometric-cutting} against the LQG field therefore gives a closed equation in terms of quantum boundary lengths; see Proposition~\ref{prop:decorated-cut}.

We then transfer this equation to the annulus setup, and obtain an equation for the annulus partition function in Section~\ref{sec:welding}. This follows from the construction of the Brownian annulus~\cite{ARSAnnulus} and its application to percolation in~\cite{SXZ}. After Mellin transformation, this yields an integro-differential equation in Proposition~\ref{weld:transform} with an explicit kernel, which will lead to the integral operator in Theorem~\ref{thm:main}.

The remaining task is to analyze the integro-differential equation of Proposition~\ref{weld:transform} and derive an exact spectral formula for $\xi_3$. In Section~\ref{sec:endpoint}, we combine arm event estimates with moment bounds for Gaussian multiplicative chaos (GMC)~\cite{Remy,RobertVargas} to identify the boundary condition satisfied by the transformed observable. The general spectral threshold theorem, Theorem~\ref{thm:critical-spectral}, then identifies the convergence endpoint of the Mellin transform with the spectral threshold, proving Theorem~\ref{thm:main}. Appendix~\ref{sec:determinant} gives an alternative expression via Fredholm determinant and proves Proposition~\ref{prop:intro-determinant}.

The above approach extends to the monochromatic three-arm exponents of $\CLE_\kappa$ for all $4<\kappa<8$. The computation in Section~\ref{sec:lqg} is valid for general $4<\kappa<8$, and we provide the candidate exponent in Appendix~\ref{sec:general}. The additional technical input is the corresponding spectral threshold theorem and analogs of the estimates determining the boundary condition in Section~\ref{sec:endpoint} for general $\kappa$. We leave the details to future work.

\begin{remark}[SLE-martingale approach]\label{rem:direct-sle}
In~\cite[Appendix~B]{LSWOneArm}, based on the SLE exploration, the backbone exponent $\xi_2$ is reduced to an eigenvalue problem for a second-order PDE. Recently, its explicit solution was obtained with the assistance of AI in~\cite{MW26}, thereby recovering the result of~\cite{NQSZ}. For $\xi_3$, one might also use the geometric observation in Proposition~\ref{prop:geometric-cutting} to develop an SLE-martingale approach that tracks both the conformal radius and the harmonic measure of the inherited boundary arc in the component retained after cutting a complete boundary-touching loop. It would be interesting to determine whether this approach yields an explicit characterization of $\xi_3$ and, if so, how that relates to the one in Theorem~\ref{thm:main}.
\end{remark}

\subsection{Monochromatic $n$-arm exponents for $n\ge4$}\label{sec:higher}
The approach above provides a robust method for determining all of the monochromatic $n$-arm exponents $\xi_n$ for critical percolation. The geometric observation in Section~\ref{sec:geometry} extends to $n\ge4$ by introducing additional boundary marks and recording the possible disjoint black connections from these marks to the inner boundary. Cutting complete boundary-touching CLE loops then yields a closed system of identities for the marked events, and conformal welding expresses these identities in terms of quantum boundary lengths. See the supplementary note~\cite{Cai26}.

Here we briefly present the spectral characterization for every $\xi_n$ with $n\ge4$. Write $r(A)$ for the spectral radius of an operator $A$, and let
\[
 \Delta_{n-1}=\left\{(x_1,\ldots,x_{n-1})\in(0,\infty)^{n-1}:\sum_{i=1}^{n-1}x_i=1\right\},
 \qquad dx=dx_1\cdots dx_{n-2},
\]
where $x_{n-1}=1-\sum_{i=1}^{n-2}x_i$. Here $x_i$ is the quantum length of the $i$th successive outer boundary arc divided by the outer perimeter. We write $L^1(\Delta_{n-1},dx)$ for the space of functions integrable with respect to this coordinate Lebesgue measure.

\begin{theorem}[Higher monochromatic arm exponents]\label{thm:higher-arms}
For each $n\ge4$ and $0<\beta<1$, let $\mathcal L_{n,\beta}$ be the positive operator on $L^1(\Delta_{n-1},dx)$ defined in~\cite[(1.1)]{Cai26}. There is a unique $\beta_n\in(0,1)$ such that
\begin{equation}\label{ha:exponent}
 r(\mathcal L_{n,\beta_n})=1,\qquad
 \xi_n=\frac{(n+\beta_n)^2-1}{12}.
\end{equation}
The spectral radius is less than one for $0<\beta<\beta_n$ and greater than one for $\beta_n<\beta<1$. At $\beta=\beta_n$, the eigenvalue one is algebraically simple.
\end{theorem}

In~\cite{Cai26}, we explicitly define the operator $\mathcal L_{n,\beta}$ and sketch how the proof for $\xi_3$ extends to $\xi_n$ with $n\ge4$. The Laplace transform of the joint boundary-length measure defining $\mathcal L_{n,\beta}$ is expressed in terms of one elementary kernel, two endpoint functions and a finite recursion.

The main new geometric step is the construction of a connection event requiring, for each choice of a marked point, two black paths from that point and one from every other mark to distinct points of the inner boundary, disjoint except at their common starting point; see \cite[Section~2]{Cai26}. The additional analytic work consists of estimates for simultaneous collisions of marks and control of the essential spectral radius, including contributions from cutting sequences that leave some original arc lengths unchanged. This is sketched in~\cite[Sections~4 and~5]{Cai26}.

Finite-dimensional approximations of $\mathcal L_{n,\beta}$ give the following numerical evaluations for $\xi_n$:
\begin{center}
\begin{tabular}{@{}lcccc@{}}
\toprule
$n$ & $4$ & $5$ & $6$& $7$\\
Numerical values for $\xi_n$ & $\approx1.35$ & $\approx2.10$ & $\approx3.0$& $\approx4.1$\\
\bottomrule
\end{tabular}
\end{center}
For $n=4,5,6$, these values are consistent with the transfer-matrix estimates in~\cite{JZJ}.
We discuss the difficulties in improving numerical precision for $n\ge4$ in~\cite[Remark~1.1]{Cai26}.

\paragraph{Acknowledgment.}
We are grateful to Xin Sun for encouraging us to work on this problem, and to Zhuoyan Xie and Zijie Zhuang for helpful discussions. We also thank Qi Wang, Youjin Deng, and Ming Li for sharing their recent numerical results~\cite{WDL26}. This work was supported by National Key R\&D Program of China (No.2023YFA1010700) and the National Natural Science Foundation of China (Grant No.12526204).

\paragraph{AI statement.}
The author developed the proof strategy, including the geometric exploration, the LQG derivation, and the reduction to a renewal equation and a spectral threshold problem. Within this framework, AI tools were used in carrying out calculations, technical estimates, and proofs in functional analysis, especially for $n\ge4$. The author also used AI tools for language editing and proofreading. The author takes full responsibility for the content of this paper.

\section{The annular event and a two-point boundary observable}
\label{sec:geometry}

Throughout this section we use the nested $\CLE_6$ describing the full scaling limit of critical Bernoulli percolation~\cite{Smirnov,CNFull,GPS}. The black and white sides of each interface are inherited from percolation.

Let $D$ be a simply connected domain and let $K\Subset D$ be the closure of a connected open set such that $D\setminus K$ is doubly connected. We impose black boundary conditions on both boundary components, and boundary points and arcs of $D$ are understood as prime ends. White connections after removing $K$ are defined by the path components in \eqref{eq:geo:induce} below.

\subsection{White clusters and necklaces}\label{subsec:necklaces}

Consider critical site percolation on the triangular lattice of mesh size $\delta$. The crossings are vertex-disjoint, including at boundary vertices, and each stops at its first visit to the inner boundary. Let $H_3^\delta(D,K)$ be the event of three vertex-disjoint black crossings from the outer to the inner boundary.

A \emph{quasi-white circuit} is a simple circuit separating the two boundary components whose vertices are white except for a specified finite set of black vertices. Menger's theorem and planar duality give
\begin{equation}\label{eq:geo:defect}
 (H_3^\delta(D,K))^c
 =\{\text{there is a separating circuit with at most two black vertices}\}.
\end{equation}
If both white portions of a quasi-white circuit with two black vertices belong to the same white cluster in $D\setminus K$, join them there. Cutting along a simple subpath of this connection decomposes the separating homology class into circuits with at most one black vertex, one of which still separates. Thus a minimal obstruction uses one white cluster or two distinct white clusters.

We use the \emph{external frontier} of each white cluster, obtained at finite mesh by closing passages of width one in its black external vertex boundary, as in~\cite[Section~2.2]{NQSZ}. At each passage, keep the side specified by the quasi-white circuit. In the limit, this frontier lies in the complete interface trace. A separating circuit with at most two black vertices becomes a \emph{necklace} of one or two complete interfaces, each with winding number one in the complementary regions on its white side and with the prescribed side at every contact.

Let $H_3(D,K)$ denote the absence of such a necklace. At the inner boundary, enlarge $K$ slightly before passing to the limit, as specified below. Members of a necklace may touch. For round annuli, the necklace argument of \cite[Section~4]{BN} identifies this measurable event in the full percolation scaling limit and gives
\begin{equation}\label{eq:geo:H3limit}
 \P[H_3^\delta(D,K)]\longrightarrow \P[H_3(D,K)]
 \qquad (\delta\downarrow0).
\end{equation}
For a rough inner set, such as the filled metric ball used below, first uniformize $D$ by a disk and $D\setminus K$ by $\{r<|w|<1\}$. Replace $r$ by $r+\varepsilon$, retain the original outer boundary and its prime ends, and let $\varepsilon\downarrow0$. For lattice comparisons, first fix these boundary approximations and a positive cutoff on interface diameters. A change in a separating circuit containing at most two black vertices across a fixed annulus forces six alternating arms: two macroscopic black passages through a contact corresponding to a single black vertex in the lattice approximation must be separated by white sectors. The same event controls penetration into a macroscopic fjord closed when forming the external frontier. The bulk six-arm exponent is greater than two, so a union bound over $O(\varepsilon^{-2})$ boxes of side $\varepsilon$ makes these errors vanish. At the outer boundary, use the boundary three-arm exponent greater than one and $O(\varepsilon^{-1})$ half-boxes; see \cite[Theorems~3 and~4]{SWArms}.

Apply these estimates to the finitely many interfaces before removing the boundary approximations, as in the frontier comparison of \cite[Lemma~2.3]{NQSZ}. Then decrease the diameter cutoff. A vanishing interface in a macroscopic obstruction involving two white clusters either leaves a one-cluster obstruction or forces the same six-arm event at the two black vertices, as in \cite[Section~4]{BN}.

To define white connections after removing $K$, we complete each white cluster in $D$ as in~\cite[Section~1.2]{MYBM}. Orient its outer interface $\eta$ to have winding number one in the complementary regions on its white side. Let $C_\eta$ be the union of those regions and $\Gamma_\eta$ the next generation of loops in them. For points $u,v$ in a set $V$ and a collection of loops $\Sigma$, let $P(u,v;V;\Sigma)$ be the set of paths in $V$ from $u$ to $v$ that do not cross any loop in $\Sigma$. Take the completion of the points of $C_\eta$ lying neither on nor inside a loop of $\Gamma_\eta$ with respect to the distance
\[
 d_{\mathrm{path},\eta}(u,v)=\inf\{\diam_{\mathrm E}(\gamma):
   \gamma\in P(u,v;\overline{C_\eta};\{\eta\}\cup\Gamma_\eta)\},
\]
where $\diam_{\mathrm E}$ denotes Euclidean diameter. Write $X_\eta$ for the completion of this white cluster and $\Pi_\eta\colon X_\eta\to\overline D$ for the continuous extension of the planar embedding.

The \emph{white clusters in $D\setminus K$} are the path components of
\begin{equation}\label{eq:geo:induce}
                       X_\eta\setminus\Pi_\eta^{-1}(K).
\end{equation}
Write $\eta_C$ for the original interface loop of a white cluster $C$. For a white cluster $C$ in $D$, let $X_C=X_{\eta_C}$ be its completion, and set $\Pi_C=\Pi_{\eta_C}$ and $d_{\mathrm{path},C}=d_{\mathrm{path},\eta_C}$. For a white cluster $C$ in $D\setminus K$, let $X_C$ be its completion, given by its closure in $X_{\eta_C}$, and set $\Pi_C=\Pi_{\eta_C}|_{X_C}$. The prescribed side at a contact is the side from which the white path approaches the interface, and we retain the corresponding endpoint in $X_C$.
The outer interface of $X_\eta$ is $\eta$, also when $\eta$ touches $\partial D$; see~\cite[Sections~1.2 and~2.4]{MYBM}.

The space $X_\eta$ admits a compact geodesic metric inducing its $d_{\mathrm{path},\eta}$ topology \cite[Section~2.4]{MYBM}. It is therefore separable and locally path connected, so the components in \eqref{eq:geo:induce} are open and countable. Enumerate them by the first points they contain in an independent sequence sampled from the normalized gasket measure of full support \cite[Section~2.3]{MYBM}. This sequence is almost surely dense. Membership in one component is tested by finite chains of geodesic balls with rational radii centered at these points, whose closures avoid $\Pi_\eta^{-1}(K)$. The necklace event is universally measurable in the space of loop collections of \cite[Section~2.2]{CNFull}, with the inherited orientations.

We use the following contact properties of $\CLE_6$: almost surely every interior point has at most two interface visits, counting repeated visits separately, and every boundary point has at most one. The interior assertion is the absence of triple points in \cite[Theorem~2(1)]{CNFull}; the boundary assertion follows from the half-plane three-arm estimate \cite[Theorem~3]{SWArms}.

\subsection{The two-point boundary observable}\label{subsec:two-points}

Fix distinct $p,q\in\partial D$. Attach a deterministic white arc outside each boundary arc from $p$ to $q$, leaving the configuration inside $D$ unchanged. The two exterior arcs meet only at $p,q$. On the lattice, keep the boundary black. Several paths may start at either $p$ or $q$, and may follow the black boundary before entering the interior. Neither marked point may occur elsewhere along a path.

In the following, an \emph{original white cluster} is a white cluster in $D\setminus K$ before the exterior arcs are added. Let $\lambda_p$, $\lambda_q$, and $\lambda_{pq}$ be the maximum numbers of vertex-disjoint black paths (except at common starting points) from $p$, from $q$, and from $\{p,q\}$, respectively, to the inner boundary, with distinct endpoints there. The corresponding lattice connection conditions are
\begin{equation}\label{eq:geo:fan}
\lambda_p\ge2,\qquad \lambda_q\ge2,\qquad \lambda_{pq}\ge3.
\end{equation}
Write $P_0,P_1$ for the exterior white arcs.

\begin{definition}[Two-point boundary observable]\label{def:two-point-observable}
A \emph{separating configuration} is of one of the following four types and separates the marked points $p,q$ from $K$ (see Figure~\ref{fig:four-witnesses}):
\begin{enumerate}[label=(\alph*)]
 \item one original white cluster;
 \item two original white clusters;
 \item an original white cluster together with $P_0$;
 \item an original white cluster together with $P_1$.
\end{enumerate}
Define $G(D,K;p,q)$ as the indicator that no separating configuration occurs. We also call a separating configuration involving one original white cluster and one exterior arc a \emph{mixed separating configuration}.
\end{definition}

\begin{figure}[!htb]
\centering
\begingroup
\definecolor{xiBlue}{RGB}{211,235,243}
\definecolor{xiRose}{RGB}{247,214,223}
\definecolor{xiGold}{RGB}{255,231,174}
\begin{tikzpicture}[x=.8cm,y=.8cm,line width=.65pt,line cap=round,line join=round,every node/.style={font=\small,inner sep=1pt},>={Latex[length=1.5mm,width=1.1mm]}]
\path[use as bounding box] (-1.65,-2.62) rectangle (18.00,2.05);
\node at (1.75,-2.5) {(a) One cluster: two possibilities};
\begin{scope}[shift={(0,0)}]
\draw (0,0) circle (1.45);
\path[fill=xiBlue,draw,even odd rule]
 (0,1.08) .. controls (.67,1.14) and (1.16,.61) .. (1.12,.03) .. controls (1.13,-.66) and (.57,-1.12) .. (-.05,-1.08) .. controls (-.72,-1.09) and (-1.17,-.57) .. (-1.10,.02) .. controls (-1.10,.65) and (-.64,1.10) .. (0,1.08) -- cycle
 (0,.64) .. controls (-.40,.65) and (-.70,.36) .. (-.66,-.02) .. controls (-.66,-.42) and (-.31,-.69) .. (.04,-.65) .. controls (.45,-.65) and (.70,-.35) .. (.65,.04) .. controls (.63,.43) and (.34,.65) .. (0,.64) -- cycle;
\filldraw[fill=black!7] (0,0) circle (.25); \node at (0,0) {$K$}; \node at (0,-.88) {$C$};
\end{scope}
\begin{scope}[shift={(3.50,0)}]
\draw (0,0) circle (1.45);
\path[fill=xiBlue,draw,even odd rule] (0.92185,0.74894)
 .. controls (1.05565,0.54348) and (1.17720,0.20263) .. (1.13000,-0.13000)
 .. controls (1.07000,-0.73000) and (0.53000,-1.14000) .. (-0.10000,-1.08000)
 .. controls (-0.77000,-1.07000) and (-1.16000,-0.53000) .. (-1.09000,0.08000)
 .. controls (-1.06000,0.80000) and (-0.51000,1.28000) .. (-0.03821,1.21329)
 .. controls (0.34907,1.20832) and (0.71600,1.06505) .. (0.92185,0.74894)
 .. controls (0.74882,0.75711) and (0.57578,0.73116) .. (0.40275,0.67108)
 .. controls (0.24076,0.61951) and (0.07791,0.54997) .. (-0.09742,0.44422)
 .. controls (-0.40000,0.62000) and (-0.74000,0.23000) .. (-0.64000,-0.11000)
 .. controls (-0.62000,-0.48000) and (-0.26000,-0.71000) .. (0.12000,-0.65000)
 .. controls (0.52000,-0.62000) and (0.72000,-0.26000) .. (0.64119,0.07892)
 .. controls (0.47986,0.39681) and (0.17559,0.48056) .. (-0.09742,0.44422)
 .. controls (0.07761,0.87069) and (0.24074,0.93241) .. (0.40275,0.67108)
 .. controls (0.57578,0.98176) and (0.74882,1.00771) .. (0.92185,0.74894) -- cycle;
\filldraw[fill=black!7] (0,0) circle (.25); \node at (0,0) {$K$};
\fill (-0.09742,0.44422) circle (.025) (0.40275,0.67108) circle (.025) (0.92185,0.74894) circle (.025);

\node at (.03,-.89) {$C$};
\end{scope}
\begin{scope}[shift={(7.60,0)}]
\draw (0,0) circle (1.45);
\path[fill=xiBlue,draw] (-1.18440,0.03023)
 .. controls (-1.29720,0.61100) and (-0.64860,1.09040) .. (-0.03760,1.06220)
 .. controls (0.62980,1.11860) and (1.30660,0.61100) .. (1.18440,0.05688)
 .. controls (1.04653,0.48745) and (0.90867,0.54403) .. (0.77080,0.22662)
 .. controls (0.63920,0.55943) and (0.50760,0.57043) .. (0.37600,0.25963)
 .. controls (0.33840,0.53580) and (0.11280,0.63920) .. (-0.11280,0.56400)
 .. controls (-0.35720,0.50760) and (-0.37600,0.12220) .. (-0.37600,-0.15651)
 .. controls (-0.50760,0.17268) and (-0.63920,0.18281) .. (-0.77080,-0.12611)
 .. controls (-0.90867,0.20557) and (-1.04653,0.25769) .. (-1.18440,0.03023) -- cycle;
\path[fill=xiRose,draw] (-1.18440,0.03023)
 .. controls (-1.31600,-0.54520) and (-0.64860,-1.12800) .. (0.00000,-1.06220)
 .. controls (0.67680,-1.08100) and (1.32540,-0.45120) .. (1.18440,0.05688)
 .. controls (1.04653,0.13651) and (0.90867,0.19309) .. (0.77080,0.22662)
 .. controls (0.63920,0.25863) and (0.50760,0.26963) .. (0.37600,0.25963)
 .. controls (0.60160,-0.16920) and (0.46060,-0.65800) .. (0.01880,-0.61100)
 .. controls (-0.29140,-0.63920) and (-0.47000,-0.40420) .. (-0.37600,-0.15651)
 .. controls (-0.50760,-0.16572) and (-0.63920,-0.15559) .. (-0.77080,-0.12611)
 .. controls (-0.90867,-0.09523) and (-1.04653,-0.04311) .. (-1.18440,0.03023) -- cycle;
\fill (1.18440,0.05688) circle (.025) (0.77080,0.22662) circle (.025) (0.37600,0.25963) circle (.025) (-0.37600,-0.15651) circle (.025) (-0.77080,-0.12611) circle (.025) (-1.18440,0.03023) circle (.025);

\filldraw[fill=black!7] (0,0) circle (.25); \node at (0,0) {$K$}; \node at (0,.78) {$C_1$}; \node at (0,-.78) {$C_2$};
\node at (0,-2.5) {(b) Two clusters};
\end{scope}
\begin{scope}[shift={(11.7,0)}]
\draw (0,0) circle (1.45);
\draw[line width=2.3pt,black!38] (0:1.56) arc[start angle=0,end angle=180,radius=1.56];
\draw[densely dashed] (180:1.56) arc[start angle=180,end angle=360,radius=1.56];
\node at (0,1.9) {$P_0$}; \node at (0,-1.9) {$P_1$};
\path[fill=xiBlue,draw] (-0.49593,1.36255)
 .. controls (-0.55000,1.10000) and (-0.76000,1.02000) .. (-0.93204,1.11076)
 .. controls (-0.93000,0.85000) and (-1.08000,0.91000) .. (-1.22967,0.76838)
 .. controls (-1.40000,-0.30000) and (-0.85000,-1.14000) .. (0.00000,-1.14000)
 .. controls (0.85000,-1.14000) and (1.40000,-0.30000) .. (1.22967,0.76838)
 .. controls (1.08000,0.91000) and (0.93000,0.85000) .. (0.93204,1.11076)
 .. controls (0.76000,1.02000) and (0.55000,1.10000) .. (0.49593,1.36255)
 .. controls (0.53000,0.85000) and (0.83000,0.20000) .. (0.66000,-0.29000)
 .. controls (0.52000,-0.87000) and (-0.52000,-0.87000) .. (-0.66000,-0.29000)
 .. controls (-0.83000,0.20000) and (-0.53000,0.85000) .. (-0.49593,1.36255) -- cycle;
\fill (-0.49593,1.36255) circle (.030) (-0.93204,1.11076) circle (.030) (-1.22967,0.76838) circle (.030) (1.22967,0.76838) circle (.030) (0.93204,1.11076) circle (.030) (0.49593,1.36255) circle (.030);
\filldraw[fill=black!7] (0,0) circle (.25); \node at (0,0) {$K$};
\node at (0,-0.92) {$C$};
\fill (-1.45,0) circle (.034) (1.45,0) circle (.034); \node[left] at (-1.67,0) {$p$}; \node[right] at (1.67,0) {$q$};
\node at (0,-2.5) {(c) With $P_0$};
\end{scope}
\begin{scope}[shift={(15.85,0)}]
\draw (0,0) circle (1.45);
\draw[densely dashed] (0:1.56) arc[start angle=0,end angle=180,radius=1.56];
\draw[line width=2.3pt,black!38] (180:1.56) arc[start angle=180,end angle=360,radius=1.56];
\node at (0,1.9) {$P_0$}; \node at (0,-1.9) {$P_1$};
\begin{scope}[yscale=-1]
\path[fill=xiBlue,draw] (-0.49593,1.36255)
 .. controls (-0.55000,1.10000) and (-0.76000,1.02000) .. (-0.93204,1.11076)
 .. controls (-0.93000,0.85000) and (-1.08000,0.91000) .. (-1.22967,0.76838)
 .. controls (-1.40000,-0.30000) and (-0.85000,-1.14000) .. (0.00000,-1.14000)
 .. controls (0.85000,-1.14000) and (1.40000,-0.30000) .. (1.22967,0.76838)
 .. controls (1.08000,0.91000) and (0.93000,0.85000) .. (0.93204,1.11076)
 .. controls (0.76000,1.02000) and (0.55000,1.10000) .. (0.49593,1.36255)
 .. controls (0.53000,0.85000) and (0.83000,0.20000) .. (0.66000,-0.29000)
 .. controls (0.52000,-0.87000) and (-0.52000,-0.87000) .. (-0.66000,-0.29000)
 .. controls (-0.83000,0.20000) and (-0.53000,0.85000) .. (-0.49593,1.36255) -- cycle;
\fill (-0.49593,1.36255) circle (.030) (-0.93204,1.11076) circle (.030) (-1.22967,0.76838) circle (.030) (1.22967,0.76838) circle (.030) (0.93204,1.11076) circle (.030) (0.49593,1.36255) circle (.030);
\end{scope}
\filldraw[fill=black!7] (0,0) circle (.25); \node at (0,0) {$K$};
\node at (0,0.92) {$C$};
\fill (-1.45,0) circle (.034) (1.45,0) circle (.034); \node[left] at (-1.67,0) {$p$}; \node[right] at (1.67,0) {$q$};
\node at (0,-2.5) {(d) With $P_1$};
\end{scope}
\end{tikzpicture}
\endgroup
\caption{The four separating configurations. Blue and pink indicate distinct original white clusters. In (c)--(d), closing a crosscut in the white cluster along the thick exterior arc separates $K$ from the marked points $p,q$; the dashed arc is the other exterior arc. At contact points, the paths approach along their specified white sides.}
\label{fig:four-witnesses}
\end{figure}

Away from neighborhoods of the marked points carrying a macroscopic white cluster, Menger's theorem and planar duality show that failure of $\lambda_{pq}\ge3$ produces one of the four listed separating configurations. For the condition $\lambda_p\ge2$, delete the marked point $q$, which is not used as a starting point, joining the exterior arcs. A set of at most one black vertex other than $p$ that meets every such path then gives an obstruction formed by one white cluster, since the exterior arcs have no additional common contact point. The argument for $\lambda_q\ge2$ is identical. Conversely, every listed separating configuration violates at least one condition in~\eqref{eq:geo:fan}.

For fixed $p,q$, the probability of the excluded configurations tends to zero: a macroscopic white cluster at a marked point requires a boundary one-arm event in a shrinking half-annulus. If a microscopic cluster belonged to an omitted macroscopic separating configuration, another component of that separating configuration would cross the same half-annulus, with the same vanishing bound. Thus exclusion of the four listed separating configurations is the continuum version of \eqref{eq:geo:fan}.

\subsection{Boundary arcs and integration of the marks}
\label{subsec:shadow-allocation}

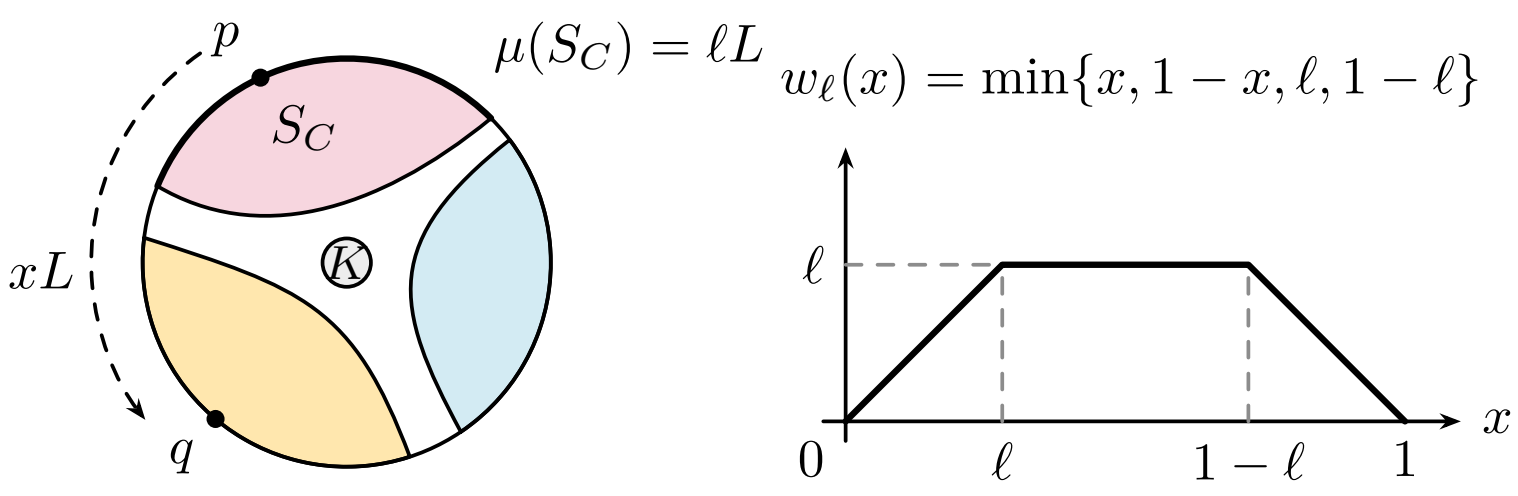
\begin{figure}[!htb]
\centering
\begin{tikzpicture}[x=1cm,y=1cm,line width=.7pt,
  every node/.style={font=\small},line cap=round,line join=round]
\begin{scope}
 \path[fill=xiRose,draw] (45:1.35)
  arc[start angle=45,end angle=158,radius=1.35]
  .. controls (-.60,.12) and (.18,.32) .. (45:1.35);
 \path[fill=xiGold,draw] (173:1.35)
  arc[start angle=173,end angle=288,radius=1.35]
  .. controls (.04,-.28) and (-.29,-.17) .. (173:1.35);
 \path[fill=xiBlue,draw] (304:1.35)
  arc[start angle=304,end angle=397,radius=1.35]
  .. controls (.23,.14) and (.30,-.28) .. (304:1.35);
 \draw (0,0) circle (1.35);
 \filldraw[fill=black!7] (0,0) circle (.16);
 \node[font=\footnotesize] at (0,0) {$K$};
 \node at (-.28,.89) {$S_C$};
 \draw[very thick] (45:1.35)
  arc[start angle=45,end angle=158,radius=1.35];
 \node[above right] at (53:1.38) {$\mu(S_C)=\ell L$};
 \fill (115:1.35) circle (1.7pt);
 \fill (230:1.35) circle (1.7pt);
 \node[above left] at (115:1.35) {$p$};
 \node[below left] at (230:1.35) {$q$};
 \draw[-{Stealth[length=4pt]},dashed] (125:1.69)
   arc[start angle=125,end angle=218,radius=1.69];
 \node[left] at (-1.66,-.05) {$xL$};
 \node[align=center] at (0,-2.12) {Different maximal arcs\\give $G=1$ on $H_3$.};
\end{scope}
\begin{scope}[shift={(3.3,-1.05)},x=3.7cm,y=3.7cm]
 \draw[-{Stealth[length=4pt]}] (-.04,0)--(1.10,0) node[right] {$x$};
 \draw[-{Stealth[length=4pt]}] (0,-.035)--(0,.49);
 \draw[dashed,black!45] (.28,0)--(.28,.28);
 \draw[dashed,black!45] (.72,0)--(.72,.28);
 \draw[dashed,black!45] (0,.28)--(.28,.28);
 \draw[line width=1.2pt] (0,0)--(.28,.28)--(.72,.28)--(1,0);
 \node[below left] at (0,0) {$0$};
 \node[below] at (.28,0) {$\ell$};
 \node[below] at (.72,0) {$1-\ell$};
 \node[below] at (1,0) {$1$};
 \node[left] at (0,.28) {$\ell$};
 \node[align=center] at (.51,.61)
   {$w_\ell(x)=\min\{x,1-x,\ell,1-\ell\}$};
 \node[align=center] at (.5,-.285)
   {The contribution of one arc\\after averaging the original mark.};
\end{scope}
\end{tikzpicture}
\caption{Averaging boundary marks gives the piecewise linear weight in Lemma~\ref{geo:allocation}. The left panel shows three maximal discarded boundary arcs, and omitted arcs fill the gaps up to a null set. Average $p$ against $\mu$ and place $q$ a counterclockwise $\mu$-length $xL$ from $p$. The right panel gives the contribution of an arc of normalized length $\ell<1/2$. The measure $-w_\ell''$ has unit atoms at $\ell,1-\ell$, and clusters meeting $K$ contribute them to the inhomogeneous term. The thick arc is the discarded boundary arc $S_C$.}
\label{fig:shadow-allocation}
\end{figure}

For $x,\ell\in(0,1)$, set $w_\ell(x)=\min\{x,1-x,\ell,1-\ell\}$. Averaging the two-point observable $G$ over the marks gives a sum of these piecewise linear weights, one for each maximal discarded boundary arc (Figure~\ref{fig:shadow-allocation}).

Let $C$ be a white cluster in $D\setminus K$ whose planar image $\Pi_C(X_C)$ meets the outer boundary but does not surround the inner set. Retain the component of $D\setminus\Pi_C(X_C)$ containing the connected interior of $K$, and contacts with $\partial K$ use the inner-boundary approximation described above. The part of $\partial D$ accessible from this component is the \emph{inherited boundary arc} $J_C$. Put
\[
                      S_C=\partial D\setminus J_C.
\]
The \emph{discarded boundary arc} $S_C$ is the arc of the original boundary cut off by $C$. For a white cluster in $D\setminus K$, it is the union of the arcs cut off by its white crosscuts. When $C$ is a white cluster in $D$ with $\Pi_C(X_C)\cap K=\varnothing$, cutting its interface loop $\eta_C$ discards $S_C$ from the component containing $K$.

Write $H_3=H_3(D,K)$ when the domain and inner set are fixed. On $H_3$, any two discarded boundary arcs have disjoint interiors or one contains the other. Indeed, two arcs with overlapping interiors, neither containing the other, give a two-component necklace if they belong to different white clusters. If they belong to the same cluster, joining the corresponding white portions gives a separating configuration consisting of one white cluster. Two arcs covering the whole boundary give the same forbidden separating configuration. It follows that the maximal arcs have disjoint interiors. Let $\mathcal M$ be the set of white clusters whose discarded boundary arcs are maximal.

Let $\mu$ be a finite nonatomic boundary measure with full support, and write $L=\mu(\partial D)$ and $\ell_C=\mu(S_C)/L$. Every discarded boundary arc of positive $\mu$-measure lies in a maximal one. Indeed, fix such an arc $S_0$ in a disk uniformization. All arcs containing $S_0$ have at least its mass. Since $K$ is separated from the original boundary, an interface loop of small diameter cuts off only a small arc of that boundary away from $K$. The mass of an arc tends to zero uniformly with its diameter, so the interface loops corresponding to arcs containing $S_0$ have diameters bounded away from zero. CLE local finiteness leaves finitely many such loops.

Each such loop also has only finitely many white clusters in $D\setminus K$ meeting the original boundary. Indeed, if $\eta$ is the original loop, in the compact geodesic space $X_\eta$ put $F_o=\Pi_\eta^{-1}(\partial D)$ and $F_i=\Pi_\eta^{-1}(K)$. If both are nonempty, their distance is $\Delta>0$. Cover $F_o$ by finitely many balls of radius $\Delta/3$ centered in $F_o$. Each ball is path connected and avoids $F_i$, so it meets only one component of $X_\eta\setminus F_i$. If either set is empty the assertion is immediate. Thus only finitely many discarded boundary arcs contain $S_0$, so one of them is maximal.

Each nondegenerate maximal arc belongs to a unique white cluster in $D\setminus K$. Two distinct original interface loops would both visit its endpoints on the original boundary, contrary to the boundary contact property above. For a fixed loop $\eta$, each contact with the original boundary has a unique preimage under $\Pi_\eta$: the argument in the proof of Lemma~\ref{geo:cut-connectivity}, applied before deleting $K$ in a half-disk neighborhood of the boundary, would otherwise force a second boundary visit. This preimage lies outside $\Pi_\eta^{-1}(K)$ and in a unique white cluster in $D\setminus K$, so two such clusters cannot have the same discarded arc.

\begin{lemma}\label{geo:allocation}
On $H_3$, assume that the union of the maximal discarded arcs has full $\mu$-measure in $\partial D$. Identify $\partial D$ with $\mathbb R/L\mathbb Z$ using counterclockwise $\mu$-length. For $0<x<1$, the point $p+xL$ is a counterclockwise $\mu$-length $xL$ from $p$. Define
\[
 \overline G(x)=\frac1L\int_{\partial D}
       G(D,K;p,p+xL)\,\mu(dp).
\]
Then for any nonnegative integrable function $f$ on $(0,1)$,
\begin{align}
 \overline G(x)
   &=\mathbf1_{H_3}\sum_{C\in\mathcal M}
        w_{\ell_C}(x),\label{eq:geo:pointallocation}\\
 \int_0^1 f(x)\overline G(x)\,dx
   &=\mathbf1_{H_3}\sum_{C\in\mathcal M}\Psi_f(\ell_C),\qquad
 \Psi_f(\ell)=\int_0^1f(x)w_\ell(x)\,dx.
 \label{eq:geo:allocation}
\end{align}
\end{lemma}

\begin{proof}
Fix $x\in(0,1)$ and work on $H_3$. The uncovered boundary points and the arc endpoints form a $\mu$-null set. Exclude this set and its inverse image under the $\mu$-preserving translation $p\mapsto p+xL$. For the remaining marks $p,q=p+xL$, absence of mixed separating configurations is equivalent to $p$ and $q$ lying in distinct maximal discarded boundary arcs.

Suppose both marks lie in the interior of $S_C$, and put $\eta=\eta_C$. Choose a compact subarc $I\Subset S_C$ containing them. Crosscuts in the completion of $C$ cut off arcs covering $I$. Take a finite chain with consecutive overlaps and join the corresponding contact points on their prescribed sides by paths in $C$. Their union $A$ is compact in $X_\eta\setminus\Pi_\eta^{-1}(K)$, so its projection stays a positive distance from $K$. Each point of $I$ is cut off by a selected crosscut from the component of $D\setminus\Pi_\eta(A)$ containing $K$. Together with the appropriate exterior arc, these projected paths give a mixed separating configuration. Conversely, a mixed separating configuration cuts off both marks in an arc contained in one maximal arc. Separating configurations using only original white clusters account exactly for $H_3^c$.

For a circle of length one and an arc of length $\ell$, the set of rotations for which both $p$ and $p+x$ belong to that arc has length $(\ell-x)_++(\ell-1+x)_+$, where $y_+=\max\{y,0\}$. Since $\sum_C\ell_C=1$, the probability that the two points belong to different maximal arcs is
\[
 \sum_C\bigl(\ell_C-(\ell_C-x)_+-(\ell_C-1+x)_+\bigr)
 =\sum_C\min\{x,1-x,\ell_C,1-\ell_C\}.
\]
This gives \eqref{eq:geo:pointallocation}, and integrating against $f$ and applying Tonelli gives \eqref{eq:geo:allocation}.
\end{proof}

For the quantum disk with an independent CLE used below, condition on the Liouville field, which determines the quantum boundary length measure $\mu$ and the inner set $K$. Fix $p$ in a disk uniformization and a boundary neighborhood whose closure avoids $K$. Remove the forced black boundary row when testing a black arm from scale $\varepsilon$ to the fixed width of that neighborhood. If there is no such arm, duality gives a white crosscut with endpoints on opposite sides of $p$; on $H_3$, the white cluster in $D\setminus K$ containing this crosscut cuts off an arc containing $p$. An uncovered point therefore requires a black arm entering the domain at every scale. The boundary one-arm estimate \cite[Theorem~3]{SWArms} makes this event null for each fixed $p$. Applying Tonelli to the conditionally fixed measure $\mu$ gives
\begin{equation}\label{eq:geo:coverage}
 \mathbf1_{H_3}\,\mu\Bigl(\partial D\setminus
                  \bigcup_{C\in\mathcal M}S_C\Bigr)=0
 \quad\hbox{almost surely}.
\end{equation}
The countable set of arc endpoints has zero $\mu$-mass by nonatomicity. Thus Lemma~\ref{geo:allocation} applies to the quantum disk with its independent CLE.

\subsection{Cutting a boundary-touching loop}\label{subsec:geometric-restart}

We rewrite the terms in~\eqref{eq:geo:allocation} with $\Pi_C(X_C)\cap K=\varnothing$ in terms of the boundary two-point observable in the retained components.
\begin{proposition}[Averaged cutting identity]\label{prop:geometric-cutting}\label{cor:geometric-cutting}
Let $\mu$ be a finite nonatomic boundary measure of full support satisfying~\eqref{eq:geo:coverage}. For each white cluster $C$ in $D$ whose complete interface loop $\eta_C$ touches the boundary and whose planar image $\Pi_C(X_C)$ avoids $K$, let $U_C$ be the component of $D\setminus\eta_C$ containing $K$. Consider those $C$ for which $U_C$ is exterior to $\eta_C$ and has a nondegenerate inherited boundary arc $J_C$, with endpoints $p_C,q_C$. Then with $L$, $\overline G$, $\mathcal M$ and $\Psi_f$ as defined above, every nonnegative integrable function $f$ on $(0,1)$ satisfies
\begin{align}\label{eq:geo:decomposition}
 \int_0^1 f(x)\overline G(x)\,dx
 ={}&\mathbf1_{H_3}\sum_{\substack{C\in\mathcal M:\\ \Pi_C(X_C)\cap K\ne\varnothing}}
                  \Psi_f(\ell_C)\\
 &+\sum_{\substack{C:\\ \Pi_C(X_C)\cap K=\varnothing}}
          \Psi_f\bigl(\ell_C\bigr)
                         G(U_C,K;p_C,q_C).\nonumber
\end{align}
The second sum is over the white clusters $C$ with boundary-touching interface loops $\eta_C$ specified above.
\end{proposition}

A white cluster $C$ in $D\setminus K$ whose planar image avoids $K$ coincides with a white cluster in $D$ before removing $K$. Indeed, let $\eta=\eta_C$ and put $F=\Pi_\eta^{-1}(K)$. If a path component $W$ of $X_\eta\setminus F$ were proper and its closure avoided $F$, it would be open and closed in the connected space $X_\eta$. Thus every proper component has closure meeting $F$.

For a cluster $C$ in the second sum of~\eqref{eq:geo:decomposition}, $U_C$ is also the component of $D\setminus\Pi_C(X_C)$ used to define $J_C$. Its prime-end boundary consists of $J_C$ and a \emph{newly exposed boundary arc} $I_C$ along the external frontier of $C$. Attach exterior white arcs $P_I,P_J$ along $I_C,J_C$, respectively.

We first describe the white connections in $U_C$ needed to evaluate $G(U_C,K;p_C,q_C)$.
\begin{lemma}\label{geo:cut-connectivity}
For a boundary-touching interface loop $\eta_C$ and an exterior component $U_C$ as above, with $\Pi_C(X_C)\cap K=\varnothing$, restricting white connections to $U_C$ commutes with removing $K$. The planar images of the completions of two distinct white clusters other than $C$, one retained and one discarded, are disjoint along the part of $\eta_C$ bounding $U_C$. In a separating configuration involving finitely many contacts, portions of $P_I$ may be replaced by projections under $\Pi_C$ of paths in $X_C$, preserving the prescribed sides at contacts and the separating homology class.
\end{lemma}

\begin{proof}
Write $\eta=\eta_C$. Nesting and the contact properties of Section~\ref{subsec:necklaces} imply that $\eta$ is outermost. If $E$ is another white cluster in $D$ and $z\in\Pi_E(X_E)\cap\eta$ lies in $D$, then $E$ lies on an exterior side of $\eta$, and $\eta$ and the outer interface $\eta_E$ each visit $z$ once. Noncrossing prevents $E$ from passing between complementary regions through such contacts, and passage through a double point of $\eta$ would give a third visit.

Each such contact $z$ has a unique preimage under $\Pi_E$. Indeed, two preimages at positive $d_{\mathrm{path},E}$-distance would require at least four strands of $\eta_E$ and the loops of the next generation across annuli centered at $z$, with fixed small outer radius and inner radius tending to zero, by the argument in the proof of \cite[Lemma~1.10]{AMYMetrics}. The single visit of $\eta_E$ supplies only two such strands. The remaining strands belong to loops of the next generation whose diameters are bounded away from zero. Local finiteness therefore forces one such loop to contain $z$, contradicting the absence of triple points. The same argument applies to any white cluster in $D$ at a contact with the original boundary: in a half-disk neighborhood, two preimages would force a second boundary visit.

By compactness and local path connectedness of $X_E$, all preimages of points sufficiently near $z$ lie in one path-connected neighborhood avoiding $\Pi_E^{-1}(K)$. Thus the completions of distinct white clusters in $D\setminus K$ obtained from $E$ have disjoint planar images along the part of $\eta$ bounding $U_C$. A contact there between two distinct white clusters in $D$, both different from $C$, would give three interface visits, or a second boundary visit. This proves the first two assertions.

Finally, for each frontier portion to be replaced, choose the preimages of its endpoints according to the prescribed sides, join them by a path in $X_C$, and project under $\Pi_C$. The projected paths lie on the discarded side and avoid $K$. Each replacement and the frontier portion it replaces form a closed path in the filled hull of $\eta$, namely its trace together with the bounded components of its complement. Since $K$ lies outside this hull, the path has winding number zero about $K$. The replacements therefore preserve the prescribed sides and the separating homology class.
\end{proof}

\begin{figure}[!htb]
\centering
\begingroup
\definecolor{xiBlue}{RGB}{211,235,243}
\definecolor{xiRose}{RGB}{247,214,223}
\definecolor{xiGold}{RGB}{255,231,174}
\begin{tikzpicture}[x=1cm,y=1cm,line width=.65pt,line cap=round,line join=round,every node/.style={font=\small,inner sep=1pt},>={Latex[length=1.5mm,width=1.1mm]}]
\path[use as bounding box] (-2.65,-3.55) rectangle (9.72,2.79);
\node at (0,2.64) {(a) A larger discarded arc}; \node at (7.08,2.64) {(b) Replacement along $I_C$};
\begin{scope}[shift={(0,0)}]
\path[fill=xiRose,draw] (-1.17583,1.67926) arc[start angle=125,end angle=235,radius=2.05]
 .. controls (-0.94139,-0.55975) and (-0.94139,0.55975) .. (-1.17583,1.67926) -- cycle;
\path[fill=xiBlue,draw] (-0.96242,1.81004)
 .. controls (-0.90355,1.84134) and (-0.65160,1.54001) .. (-0.57117,1.56929)
 .. controls (-0.49074,1.59856) and (-0.49143,1.99134) .. (-0.42622,2.00520)
 .. controls (-0.36101,2.01906) and (-0.20397,1.68989) .. (-0.11859,1.69586)
 .. controls (-0.03320,1.70183) and (0.07650,2.04966) .. (0.14300,2.04501)
 .. controls (0.20951,2.04036) and (0.02395,1.48925) .. (-0.05000,1.25000)
 .. controls (-0.12395,1.01075) and (-0.30500,0.70500) .. (-0.35000,0.45000)
 .. controls (-0.39500,0.19500) and (-0.39500,-0.19500) .. (-0.35000,-0.45000)
 .. controls (-0.30500,-0.70500) and (-0.12395,-1.01075) .. (-0.05000,-1.25000)
 .. controls (0.02395,-1.48925) and (0.20951,-2.04036) .. (0.14300,-2.04501)
 .. controls (0.07650,-2.04966) and (-0.03320,-1.70183) .. (-0.11859,-1.69586)
 .. controls (-0.20397,-1.68989) and (-0.36101,-2.01906) .. (-0.42622,-2.00520)
 .. controls (-0.49143,-1.99134) and (-0.49074,-1.59856) .. (-0.57117,-1.56929)
 .. controls (-0.65160,-1.54001) and (-0.90355,-1.84134) .. (-0.96242,-1.81004)
 .. controls (-1.02128,-1.77874) and (-1.00136,-1.53651) .. (-0.98000,-1.34000)
 .. controls (-0.95864,-1.14349) and (-0.84400,-0.77600) .. (-0.82000,-0.50000)
 .. controls (-0.79600,-0.22400) and (-0.79600,0.22400) .. (-0.82000,0.50000)
 .. controls (-0.84400,0.77600) and (-0.95864,1.14349) .. (-0.98000,1.34000)
 .. controls (-1.00136,1.53651) and (-1.02128,1.77874) .. (-0.96242,1.81004) -- cycle;
\draw (0,0) circle (2.05);
\fill (-1.17583,1.67926) circle (.036) (-1.17583,-1.67926) circle (.036);
\node[above left] at (-1.17583,1.67926) {$p_C$}; \node[below left] at (-1.17583,-1.67926) {$q_C$};
\node[left] at (-2.09,.05) {$S_C$}; \node at (1.39,1.16) {$J_C$};
\filldraw[fill=black!7] (0,0) circle (.28); \node at (0,0) {$K$};
\fill (-0.96242,1.81004) circle (.030) (-0.42622,2.00520) circle (.030) (0.14300,2.04501) circle (.030) (0.14300,-2.04501) circle (.030) (-0.42622,-2.00520) circle (.030) (-0.96242,-1.81004) circle (.030);
\node at (-1.65,.02) {$C$}; \node at (-1.22,.25) {$I_C$}; \node at (-.59,0) {$E$}; \node at (.85,-.90) {$U_C$};
\draw[densely dashed] (-1.17583,1.67926) -- (-1.24466,1.77756) arc[start angle=125,end angle=-125,radius=2.17] -- (-1.17583,-1.67926);
\draw[densely dashed] (-0.96242,1.81004) -- (-1.01875,1.91600);
\draw[densely dashed] (-0.42622,2.00520) -- (-0.45117,2.12258);
\draw[densely dashed] (0.14300,2.04501) -- (0.15137,2.16471);
\draw[densely dashed] (0.14300,-2.04501) -- (0.15137,-2.16471);
\draw[densely dashed] (-0.42622,-2.00520) -- (-0.45117,-2.12258);
\draw[densely dashed] (-0.96242,-1.81004) -- (-1.01875,-1.91600);
\node[right] at (2.20,.17) {$P_J$}; \node[align=center] at (0,-2.90) {on $H_3$: $S_C\subsetneq S_E$\\$\longleftrightarrow$ mixed separating configuration\\using $P_J$};
\end{scope}
\begin{scope}[shift={(7.08,0)}]
\path[fill=xiRose,draw] (-1.17583,1.67926) arc[start angle=125,end angle=235,radius=2.05]
 .. controls (-0.94139,-0.55975) and (-0.94139,0.55975) .. (-1.17583,1.67926) -- cycle;
\path[fill=xiBlue,draw] (-1.09743,1.25000)
 .. controls (-1.10782,1.31667) and (-0.45029,1.43467) .. (-0.15000,1.42000)
 .. controls (0.15029,1.40533) and (0.70507,1.34653) .. (0.95000,1.15000)
 .. controls (1.19493,0.95347) and (1.50973,0.40707) .. (1.52000,0.08000)
 .. controls (1.53027,-0.24707) and (1.25760,-0.85853) .. (1.02000,-1.08000)
 .. controls (0.78240,-1.30147) and (0.21056,-1.40507) .. (-0.10000,-1.43000)
 .. controls (-0.41056,-1.45493) and (-1.10782,-1.31667) .. (-1.09743,-1.25000)
 .. controls (-1.08704,-1.18333) and (-0.88586,-1.15253) .. (-0.88000,-1.11000)
 .. controls (-0.87414,-1.06747) and (-1.06545,-1.02667) .. (-1.05746,-0.96000)
 .. controls (-1.04948,-0.89333) and (-0.88468,-0.84693) .. (-0.88000,-0.80000)
 .. controls (-0.87532,-0.75307) and (-1.03086,-0.70667) .. (-1.02554,-0.64000)
 .. controls (-1.02022,-0.57333) and (-0.60081,-0.78853) .. (-0.38000,-0.79000)
 .. controls (-0.15919,-0.79147) and (0.29520,-0.76587) .. (0.48000,-0.65000)
 .. controls (0.66480,-0.53413) and (0.88880,-0.19800) .. (0.88000,0.00000)
 .. controls (0.87120,0.19800) and (0.60773,0.57973) .. (0.42000,0.70000)
 .. controls (0.23227,0.82027) and (-0.18799,0.82880) .. (-0.40000,0.82000)
 .. controls (-0.61201,0.81120) and (-1.02022,0.57333) .. (-1.02554,0.64000)
 .. controls (-1.03086,0.70667) and (-0.87532,0.75307) .. (-0.88000,0.80000)
 .. controls (-0.88468,0.84693) and (-1.04948,0.89333) .. (-1.05746,0.96000)
 .. controls (-1.06545,1.02667) and (-0.87414,1.06747) .. (-0.88000,1.11000)
 .. controls (-0.88586,1.15253) and (-1.08704,1.18333) .. (-1.09743,1.25000) -- cycle;
\draw (0,0) circle (2.05);
\fill (-1.17583,1.67926) circle (.036) (-1.17583,-1.67926) circle (.036);
\node[above left] at (-1.17583,1.67926) {$p_C$}; \node[below left] at (-1.17583,-1.67926) {$q_C$};
\node[left] at (-2.09,.05) {$S_C$}; \node at (1.39,1.16) {$J_C$};
\filldraw[fill=black!7] (0,0) circle (.28); \node at (0,0) {$K$};
\fill (-1.09743,1.25000) circle (.030) (-1.09743,-1.25000) circle (.030) (-1.05746,-0.96000) circle (.030) (-1.02554,-0.64000) circle (.030) (-1.02554,0.64000) circle (.030) (-1.05746,0.96000) circle (.030);
\node at (-1.84,.02) {$C$}; \node at (-.76,-.10) {$I_C$}; \node at (1.21,.05) {$E$}; \node at (.65,-1.67) {$U_C$};
\draw[densely dashed] (-1.17583,-1.67926) -- (-1.29328,-1.65467);
\draw[densely dashed] plot[domain=-1.67926169:1.67926169,samples=81,smooth,variable=\t] ({-1.00000000-0.06235344*\t*\t-.12/sqrt(1+4*0.06235344*0.06235344*\t*\t)},{\t-.24*0.06235344*\t/sqrt(1+4*0.06235344*0.06235344*\t*\t)});
\draw[densely dashed] (-1.29328,1.65467) -- (-1.17583,1.67926);
\draw[line width=1.0pt] (-1.05746,-0.96000) .. controls (-1.78,-.83) and (-1.78,.83) .. (-1.05746,0.96000);
\draw[->] (-1.16,.02) -- (-1.49,.02); \node at (-1.44,1.10) {$P_I$}; \node[align=center] at (0,-2.90) {mixed separating configuration using $P_I$\\$\longleftrightarrow$ separating configuration in $D$\\using $E,C$};
\end{scope}
\end{tikzpicture}
\endgroup
\caption{The two mixed separating configurations used in the proof of Lemma~\ref{geo:restart}. In $U_C$, $J_C$ is the inherited boundary arc and $I_C$ the newly exposed boundary arc. The dashed curves are the added exterior white arcs $P_I,P_J$, and the blue pieces belong to one white cluster $E$. In (a), on $H_3$, a crosscut in $E$ separates $I_C$ from $K$. Restoring the pink region on the side of $C$ gives $S_C\subsetneq S_E$; closing the crosscut along $P_J$ gives a mixed separating configuration in $U_C$. In (b), replacing the used portion of $P_I$ by the projection of a path in $X_C$ gives a separating configuration in $D$ using $E,C$, with the same separating homology class and prescribed sides at contacts.}
\label{fig:restart-witnesses}
\end{figure}

\begin{lemma}[Pointwise cutting identity]\label{geo:restart}
For each boundary-touching interface loop $\eta_C$ with the exterior component $U_C$ specified above and with $\Pi_C(X_C)\cap K=\varnothing$,
\begin{equation}\label{eq:geo:restart}
 \mathbf1_{H_3(D,K)}\mathbf1_{\{C\in\mathcal M\}}
       =G(U_C,K;p_C,q_C).
\end{equation}
The observable on the right uses the original white connections restricted to $U_C\setminus K$.
\end{lemma}

\begin{proof}
Apply Lemma~\ref{geo:cut-connectivity} (see Figure~\ref{fig:restart-witnesses} for the two mixed cases). If $H_3(D,K)$ fails through a separating configuration using one white cluster, that separating configuration remains in $U_C$. A separating configuration using two clusters other than $C$ also remains there. A separating configuration using $C$ becomes a mixed separating configuration with the exterior white arc along $I_C$. Each case gives $G(U_C,K;p_C,q_C)=0$. If $H_3(D,K)$ holds but $C\notin\mathcal M$, a cluster $E$ with $S_C\subsetneq S_E$ has a crosscut in $U_C$ separating $I_C$ from $K$. Closing this crosscut along the exterior arc $P_J$ gives a mixed separating configuration. We have therefore proved
\[
 \mathbf1_{H_3(D,K)}\mathbf1_{\{C\in\mathcal M\}}=0
 \quad\Longrightarrow\quad G(U_C,K;p_C,q_C)=0.
\]

Conversely, suppose $G(U_C,K;p_C,q_C)=0$. A separating configuration in $U_C$ using one or two original white clusters remains a separating configuration in $D$ after restoring $D\setminus U_C$. For a mixed separating configuration using $P_I$, replace its used portions by projections of paths in $X_C$, as in Lemma~\ref{geo:cut-connectivity}. The prescribed sides at the contacts give a separating configuration in $D$ using two white clusters. For a mixed separating configuration using $P_J$, let $E$ be its original white cluster. A crosscut in $E$ has endpoints on $J_C$ and separates $I_C$ from $K$. By Jordan separation, the inherited boundary arc accessible from $K$ after removing $E$ is strictly contained in $J_C$. Hence $S_C\subsetneq S_E$, unless the same crosscut already gives a separating configuration in $D$ using one or two white clusters. In each case, $\mathbf1_{H_3(D,K)}\mathbf1_{\{C\in\mathcal M\}}=0$. This proves the reverse implication and hence~\eqref{eq:geo:restart}.
\end{proof}

\begin{proof}[Proof of Proposition~\ref{prop:geometric-cutting}]
Split the sum in Lemma~\ref{geo:allocation} according to whether $\Pi_C(X_C)$ meets $K$. For the terms avoiding $K$, identify $C$ with a white cluster in $D$ as above and apply Lemma~\ref{geo:restart}. All summands are nonnegative, so this regrouping is justified by Tonelli's theorem.
\end{proof}

The first term in \eqref{eq:geo:decomposition} gives the nonnegative inhomogeneous term in the renewal equation of Proposition~\ref{weld:renewal}. Each cluster contributing to this term meets both boundary components of the annulus and, on $H_3$, supplies a white crossing alongside the three black crossings. After removing fixed neighborhoods of the two boundary components, these contributions are supported on the $BBBW$ event, meaning three disjoint black crossings and one white crossing.

\section{Conformal welding and the boundary-length kernel}
\label{sec:lqg}

For every $4<\kappa<8$, we compute the boundary-length kernel for cutting a complete boundary-touching $\mathrm{CLE}_\kappa$ loop. The computation uses conformal welding of quantum surfaces and the coupling of CLE with LQG originating from~\cite{Sheffield16,DMS21,MSW20,MSW21} and later developed in e.g.~\cite{AHSInt,AHSY,ACSW,LSYZ}. In particular, we integrate the radial welding identity of~\cite[Theorem~5.19]{LSYZ} over all quantum surfaces other than the retained disk, and obtain the explicit joint measure of the boundary lengths and the retained disk. This is given in Proposition~\ref{prop:decorated-cut}.

Throughout this section,
\[
 4<\kappa<8,\qquad \gamma=\frac4{\sqrt\kappa}\in(\sqrt2,2),\qquad
 \alpha=\frac\kappa4,\qquad e=\alpha-1,\qquad
 Q_{\rm L}=\frac2\gamma+\frac\gamma2.
\]
Here $Q_{\rm L}$ is the constant in the LQG coordinate change. The quantum boundary and area measures associated with the field $h$ are $\nu_h$ and $\mu_h$. Under a conformal map $f:D\to\widetilde D$, the field changes to $h\circ f^{-1}+Q_{\rm L}\log|(f^{-1})'|$. Write $\mathrm{QD}_{i,j}$ for the ordinary quantum disk measure with $i\ge0$ area marks and $j\ge0$ ordered boundary marks, where $i,j$ are integers. Area and boundary marks are integrated against $\mu_h$ and $\nu_h$, respectively. Write $\mathrm{QD}_{1,0}(L)$ and $\mathrm{QD}_{1,2}(a,b)$ for the disintegrations of $\mathrm{QD}_{1,0}$ and $\mathrm{QD}_{1,2}$ by perimeter $L>0$ and boundary arc lengths $a,b>0$, with respect to $dL$ and $da\,db$. We write $\mathrm{QD}_{i,j}\otimes\mathrm{CLE}_\kappa$ when an independent CLE is sampled on this quantum surface. For a measure $M$, write $|M|$ for its total mass. By~\cite{DMS21}, we have $|\mathrm{QD}_{1,0}(L)|\propto L^{-\alpha}$, and
\[
 |\mathrm{QD}_{1,2}(a,b)|=(a+b)|\mathrm{QD}_{1,0}(a+b)|\propto(a+b)^{1-\alpha}.
\]

We use $\mathcal M_2^{\mathrm{disk}}(W)$ for the quantum disk measure of weight $W>0$ with two marked boundary points, with the conventions of~\cite[Definitions~3.4--3.6]{LSYZ}. Independent forests attached along its boundary give the corresponding forested quantum surface, and the original quantum surface is called its spine. Ordinary boundary lengths are measured on this spine, whereas the boundary of the attached forest carries generalized quantum length, as in~\cite[Definitions~5.1--5.3]{LSYZ}. For the exterior component containing the area mark, write $a$ for the length of the inherited boundary arc, $b$ for the length of the newly exposed boundary arc, and $s$ for the length of the discarded boundary arc. The original disk has perimeter $a+s$, and the retained disk has perimeter $a+b$.

\subsection{Cutting a complete loop}\label{subsec:exterior-bead}

Put $W=\gamma^2-2$ and $k=(\alpha-1)^{-2}$. Thus $0<W<\gamma^2/2$ and $\gamma^2-W=2$. A pinched thin quantum annulus of weight $W$ is a cyclic Poisson chain of ordinary weight-$2$ disks. Its inner and outer boundaries are forested independently, as in~\cite[Definition~3.12 and equation~(5.2)]{LSYZ}. We index the boundary-touching CLE loops by $\eta_C$; when $\kappa=6$, $C$ is the corresponding white cluster from Section~\ref{sec:geometry}.

\begin{lemma}\label{lem:prewelding-product}
Under $\mathrm{QD}_{1,0}\otimes\mathrm{CLE}_\kappa$, count each complete boundary-touching loop $\eta_C$ once and retain the component containing the area mark when it is exterior to the loop and has a nondegenerate inherited boundary arc. Let $a$ be the quantum length of the inherited boundary arc, $b$ the length of the newly exposed boundary arc, and $s$ the length of the discarded boundary arc. The resulting measure is the pushforward under welding of
\begin{equation}\label{eq:prewelding-product}
 ds\,da\,db\;\Xi_{s,b}(dE)\bigl(\mathrm{QD}_{1,2}(a,b)\otimes\mathrm{CLE}_\kappa\bigr),
\end{equation}
where $\Xi_{s,b}$ is a nonnegative $\sigma$-finite measure kernel independent of $a$. The data $E$ consists of the other annulus beads, the inner forested disk, the forests along the boundary to be welded, and the welding identifications. The two boundary marks of the retained disk are the endpoints of its inherited boundary arc.
\end{lemma}

\begin{proof}
The loops touching the original boundary are outermost. Conditional on each generation, independent CLEs are sampled in the complementary domains. For such a loop $\eta_C$, let $\mathcal I_C$ be the union of the components it surrounds, and let $A_C=\mu_h(\mathcal I_C)$. First sample a disk from $\mathrm{QD}_{1,0}$ with area mark $w$ and select the outermost CLE loop surrounding $w$, restricted to the event that it touches the boundary. By~\cite[Remark~1.9]{LSYZ}, this is the clockwise true loop of the boundary conformal loop ensemble $\mathrm{BCLE}_\kappa(0)$ surrounding $w$. Let $\operatorname{Weld}$ denote the pushforward to the resulting loop-decorated quantum surface. After foresting the boundary, \cite[Theorem~5.19, equation~(5.18)]{LSYZ} gives
\begin{equation}\label{eq:lsyz-loop-welding}
 \mathrm{QD}_{1,0}^{f}\otimes\mu^{\circlearrowright}
 =c_\kappa\int_0^\infty \operatorname{Weld}\!\left(\widetilde{\mathrm{QA}}^{f}(W;\ell'),\,\ell'\mathrm{QD}_{1,0}^{f}(\ell')\right)d\ell'.
\end{equation}
Here $\mu^{\circlearrowright}$ is the subprobability law of the boundary-touching loop selected above, $\ell'$ is the generalized boundary length used in welding, and $c_\kappa>0$ is the constant in~\cite[equation~(5.16)]{LSYZ}. The superscript $f$ denotes independent boundary foresting. The forested pinched thin annulus measure $\widetilde{\mathrm{QA}}^{f}(W;\ell')$ and disk measure $\mathrm{QD}_{1,0}^{f}(\ell')$ are disintegrated by the generalized lengths of the boundaries to be welded.

The ordinary beads of the annulus are the exterior components with an inherited boundary arc. In the construction of~\cite[Definition~5.3 and equation~(5.2)]{LSYZ}, each bead has one arc on the annulus's outer boundary and one on its inner boundary. The outer forest becomes the auxiliary forest of the reconstructed disk and is forgotten. The inner boundary is welded to the marked forested disk, and its image is $\eta_C$. The beads meet only at the annular cut points, which lie on this loop. The forests along the inner boundary attach at points and cannot connect the open interiors of two ordinary beads after the loop is removed. Their additional complementary components inherit no nondegenerate arc of the original boundary. An ordinary bead with arc lengths $(a,b)$ therefore becomes the retained component $U_C$, with inherited boundary arc $J_C$ of length $a$ and newly exposed boundary arc $I_C$ of length $b$. The sum of the outer arc lengths of all other beads is $s$.

We now place a second area mark $z$ in such an exterior bead. By~\cite[Definition~3.12]{LSYZ}, the annulus is obtained by sampling $T$ with measure $k\,dT/T$ and a Poisson process on a circle of length $T$ with intensity $du\,\mathcal M_2^{\mathrm{disk}}(2)$. Writing $\mu_U$ for the quantum area measure of a bead $U$, Mecke's formula selects one bead with measure $\mathcal M_2^{\mathrm{disk}}(2)(dU)\mu_U(dz)=\mathrm{QD}_{1,2}(dU,dz)$ and leaves an independent Poisson process of the same intensity. Rotating the circle to the selected point integrates its location over a set of length $T$, so $k\,dT/T$ becomes $k\,dT$. Disintegrating the selected bead by its two arc lengths gives $da\,db\,\mathrm{QD}_{1,2}(a,b)$. See Figure~\ref{fig:welding-product}.

\begin{figure}[!htb]
\centering
\begin{tikzpicture}[x=1cm,y=1cm,line width=.65pt,draw=black,
  line cap=round,line join=round,
  every node/.style={font=\fontsize{9}{10.5}\selectfont,inner sep=1.5pt},
  >={Latex[length=1.6mm,width=1.2mm]}]
\path[use as bounding box] (0,0) rectangle (15.55,6.10);
\node at (2.30,5.95) {(a) Before welding};
\node at (8.00,5.95) {(b) Complete loop};
\node at (13.50,5.95) {(c) Retained disk};
\begin{scope}[shift={(2.30,3.20)}]
 \coordinate (al) at (-1.30,.30);
 \coordinate (ar) at (1.30,.30);
 \coordinate (abl) at (-.70,-1.20);
 \coordinate (abr) at (.70,-1.20);
 \path[fill=xiGold,draw] (al)
   .. controls (-1.30,2.30) and (1.30,2.30) .. (ar)
   .. controls (.70,1.00) and (-.70,1.00) .. cycle;
 \path[fill=xiBlue,draw] (ar)
   .. controls (2.00,-.40) and (1.70,-1.40) .. (abr)
   .. controls (.65,-.70) and (.75,-.45) .. cycle;
 \path[fill=xiBlue,draw] (abr)
   .. controls (.50,-1.95) and (-.50,-1.95) .. (abl)
   .. controls (-.40,-1.00) and (.40,-1.00) .. cycle;
 \path[fill=xiBlue,draw] (abl)
   .. controls (-1.70,-1.40) and (-2.00,-.40) .. (al)
   .. controls (-.75,-.45) and (-.65,-.70) .. cycle;
 \filldraw[fill=black!7] (0,1.92) circle (.12);
 \filldraw[fill=black!7] (0,2.13) circle (.09);
 \filldraw[fill=black!7] (0,.725) circle (.10);
 \filldraw[fill=black!7] (0,.555) circle (.07);
 \filldraw[fill=black!7] (1.7504,-.6781) circle (.12);
 \filldraw[fill=black!7] (1.9480,-.7493) circle (.09);
 \filldraw[fill=xiRose] (0,-.10) circle (.44);
 \filldraw[fill=black!7] (-.56,-.10) circle (.12);
 \filldraw[fill=black!7] (-.75,-.10) circle (.07);
 \filldraw[fill=black!7] (.56,-.10) circle (.12);
 \filldraw[fill=black!7] (.75,-.10) circle (.07);
 \filldraw[fill=black!7] (0,-.65) circle (.11);
 \filldraw[fill=black!7] (0,-.83) circle (.07);
 \fill (0,-.10) circle (1.4pt);
 \node[below right] at (0,-.10) {$w$};
 \fill (-.08,1.34) circle (1.4pt);
 \node[right] at (-.08,1.34) {$z$};
 \node at (-.84,1.62) {$a$};
 \node at (.73,.86) {$b$};
\end{scope}
\node[align=center] at (2.30,.62)
  {ordinary annulus beads\\and an inner forested disk};

\draw[->] (4.65,3.55) -- (5.94,3.55);
\node[above] at (5.30,3.58) {weld};
\path[fill=xiBlue,draw] (8,3.50) circle (1.70);
\coordinate (wp) at ({8+1.70*cos(220)},{3.50+1.70*sin(220)});
\coordinate (wq) at ({8+1.70*cos(320)},{3.50+1.70*sin(320)});
\coordinate (wl) at ({8+1.70*cos(245)},{3.50+1.70*sin(245)});
\coordinate (wr) at ({8+1.70*cos(295)},{3.50+1.70*sin(295)});
\path[fill=xiGold,draw] (wp)
  arc[start angle=220,end angle=-40,radius=1.70]
  .. controls (8.80,3.52) and (7.20,3.52) .. cycle;
\path[fill=xiRose,draw] (wp)
  .. controls (7.20,3.52) and (8.80,3.52) .. (wq)
  .. controls (9.06,2.55) and (8.93,2.20) .. (wr)
  .. controls (8.60,2.26) and (8.35,2.55) .. (8,2.27)
  .. controls (7.65,2.55) and (7.40,2.26) .. (wl)
  .. controls (7.07,2.20) and (6.94,2.55) .. cycle;
\draw (8,3.50) circle (1.70);
\fill (wp) circle (1.4pt) (wq) circle (1.4pt);
\fill (8.13,4.08) circle (1.4pt);
\node[above right] at (8.13,4.08) {$z$};
\node at (7.42,4.24) {$U_C$};
\fill (7.87,2.76) circle (1.4pt);
\node[right] at (7.87,2.76) {$w$};
\node[above] at (8,5.25) {$J_C:\ a$};
\node at (8,3.48) {$I_C:\ b$};
\node[left] at (wp) {$p_C$};
\node[right] at (wq) {$q_C$};
\node at (9.27,2.00) {$\eta_C$};
\draw (9.14,2.10) -- (8.99884,2.32707);
\node[below] at (8,1.70) {$\partial D\setminus J_C:\ s$};
\node[align=center] at (8,.62)
  {forget the outer forest\\and remove $w$ with weight $A_C^{-1}$};

\draw[->] (10.16,3.55) -- (11.43,3.55);
\node[above] at (10.80,3.58) {retain $U_C$};
\coordinate (cp) at (12.03,2.64);
\coordinate (cq) at (14.97,2.64);
\path[draw,fill=xiGold] (cp)
  .. controls (11.90,4.93) and (15.10,4.93) .. (cq)
  .. controls (14.16,1.55) and (12.84,1.55) .. cycle;
\fill (cp) circle (1.4pt) (cq) circle (1.4pt);
\fill (13.50,3.33) circle (1.4pt);
\node[above right] at (13.50,3.33) {$z$};
\node[above] at (13.50,4.46) {$J_C:\ a$};
\node[below] at (13.50,1.72) {$I_C:\ b$};
\node[left] at (cp) {$p_C$};
\node[right] at (cq) {$q_C$};
\node at (13.50,.62)
  {$\mathrm{QD}_{1,2}(a,b)\otimes\mathrm{CLE}_\kappa$};
\end{tikzpicture}
\caption{Radial welding and the retained component in Lemma~\ref{lem:prewelding-product}.
The colored annulus beads have weight $2$, and the gray lobes are attached forests. The inner forested disk carries $w$, and the selected yellow bead carries $z$. After welding, weighting by $A_C^{-1}$ and forgetting $w$ counts $\eta_C$ once, where $A_C$ is the total quantum area surrounded by $\eta_C$. The retained disk has inherited boundary arc $J_C$ of length $a$, newly exposed boundary arc $I_C$ of length $b$, and endpoint prime ends $p_C,q_C$; the discarded boundary arc has length $s$.}
\label{fig:welding-product}
\end{figure}

Conditional on $a$, the selected bead's outer forest has total mass one. Its inner forest and the remaining factors depend on the bead only through $b$. To disintegrate by $s$, write $R(a+b)^{-1-\alpha}\,da\,db$ for the joint density of the two boundary arc lengths under $\mathcal M_2^{\mathrm{disk}}(2)$, where $R>0$, and let $\Gamma$ denote Euler's gamma function. Conditional on $T$, the sum $S_T$ of the other beads' outer lengths satisfies
\[
 \mathbb E[e^{-pS_T}]=\exp\!\left\{-T\frac{R[-\Gamma(1-\alpha)]}{\alpha}p^{\alpha-1}\right\},\qquad p>0.
\]
The variable $S_T$ has a density on $(0,\infty)$. Disintegrating by $S_T=s$ and retaining the other quantum surfaces and welding identifications gives a measure kernel $\Xi^{\mathrm{in}}_{s,b}$ for the measure with both area marks.

The union $\mathcal I_C$ is the image of the entire inner forested disk. Hence $A_C$ depends only on $E$. Since $\eta_C$ surrounds a nonempty open set and the reconstructed disk has finite area, $0<A_C<\infty$ almost everywhere. The identity
\[
 \int_{\mathcal I_C}A_C^{-1}\mu_h(dw)=1
\]
removes the mark $w$. Weight $\Xi^{\mathrm{in}}_{s,b}$ by $A_C^{-1}$ and forget $w$ to obtain $\Xi_{s,b}$.

The domain Markov property gives independent CLEs in the exterior components, conditional on the complete loop and its complementary domains. Area marking preserves this conditional independence, since the marks depend only on the quantum surfaces and the selected loop. Sampling the CLE in the selected bead gives~\eqref{eq:prewelding-product}. Tonelli's theorem and monotone convergence extend the identity from finite restrictions to all nonnegative measurable test functions.
\end{proof}

\subsection{Computation of the boundary-length kernel}\label{subsec:cut-kernel}

We next compute $|\Xi_{s,b}|$. Set
\begin{equation}\label{eq:cut-C}
 C_\alpha=-\frac\alpha\pi\sin\frac{2\pi}\alpha>0
\end{equation}
and define the positive boundary-length kernel
\begin{equation}\label{eq:cut-B-positive}
 B_\alpha(s,b)=-\frac{2C_\alpha\sin(\pi\alpha)}\pi
 \int_0^\infty
 \frac{u^\alpha(1+u)}
 {(u^{2\alpha}-2\cos(\pi\alpha)u^\alpha+1)(b+su)^3}\,du.
\end{equation}

\begin{proposition}[Boundary-length kernel]
\label{prop:decorated-cut}
Under $\mathrm{QD}_{1,0}\otimes\mathrm{CLE}_\kappa$, count each complete boundary-touching loop $\eta_C$ once and retain the component containing the area mark when it is exterior to the loop and has a nondegenerate inherited boundary arc. The joint measure of $(s,a,b)$ and the retained component with its field, marked points and CLE is
\begin{equation}\label{eq:decorated-cut}
 B_\alpha(s,b)\,ds\,da\,db\;
                   \bigl(\mathrm{QD}_{1,2}(a,b)\otimes\mathrm{CLE}_\kappa\bigr).
\end{equation}
The boundary marks are the endpoints of the inherited boundary arc. The averaged cutting identity in Proposition~\ref{prop:geometric-cutting} imposes maximality of the discarded boundary arc.
\end{proposition}

\begin{proof}
Write $\widehat B_\alpha(s,b)=|\Xi_{s,b}|$.
Let $\mathrm{GQD}$ be the unmarked generalized quantum disk measure (see~\cite[Definition~2.27]{ACSWLoop}), and let $\mathrm{GQD}_{1,0}$ be its area-marked version, using the area of the entire generalized disk. By~\cite[Theorem~2.22 and Proposition~2.34]{ACSWLoop}, after forgetting the distinguished ordinary disk, $\mathrm{QD}_{1,0}^{f}$ is a constant multiple of $\mathrm{GQD}_{1,0}$. Here one first forgets the generalized boundary mark in Proposition~2.34 and divides by the generalized perimeter. Weighting by $A_C^{-1}$ and forgetting $w$ therefore gives a constant multiple of $\mathrm{GQD}$. Absorb this constant and the constant in~\cite[(5.18)]{LSYZ} into $\lambda_\kappa>0$.

\emph{Integration over the remaining quantum surfaces.} The annulus construction uses the measure $k\,dT/T$, where $k=(\alpha-1)^{-2}$. The two boundary arc lengths under $\mathcal M_2^{\mathrm{disk}}(2)$ have joint density $R(a+b)^{-1-\alpha}\,da\,db$. For $p,r>0$, the Poisson process of these disks has Laplace exponent
\[
 \Psi(p,r)=\iint_{(0,\infty)^2}(1-e^{-pa-rb})R(a+b)^{-1-\alpha}\,da\,db
 =\frac{R[-\Gamma(1-\alpha)]}{\alpha}\frac{p^\alpha-r^\alpha}{p-r}.
\]
Consequently
\begin{equation}\label{eq:spine-transform}
 \frac{k}{\Psi(p,r)}=E_\alpha\frac{p-r}{p^\alpha-r^\alpha},\qquad E_\alpha=\frac{\alpha}{(\alpha-1)R\Gamma(2-\alpha)}.
\end{equation}
The singularity at $p=r$ is removable.

As in Lemma~\ref{lem:prewelding-product}, selecting the exterior bead by its area measure changes $k\,dT/T$ to $k\,dT$. Let $(S_T,V_T)$ be the sums of the ordinary outer and inner arc lengths of the remaining beads. Then
\[
 \mathbb E[e^{-pS_T-rV_T}]=e^{-T\Psi(p,r)}.
\]
The selected bead contributes $b$ to the ordinary inner length. The independent inner forests together have generalized length $Y_{b+V_T}$, where $(Y_v)_{v\ge0}$ is the stable subordinator of index $1/\alpha$, normalized by $\mathbb E[e^{-tY_v}]=e^{-v t^{1/\alpha}}$ for $t>0$ as in~\cite[Definition~5.1]{LSYZ}. The outer forests integrate to one. Disintegrating $\mathrm{GQD}$ by generalized perimeter $\ell>0$ gives $|\mathrm{GQD}(\ell)|=c_g\ell^{-2-1/\alpha}$ for a constant $c_g>0$~\cite[Lemma~2.28]{ACSWLoop}. The factor $\ell$ in~\eqref{eq:lsyz-loop-welding} therefore gives
\[
 \ell|\mathrm{GQD}(\ell)|=\frac{c_g}{\Gamma(1+1/\alpha)}\int_0^\infty t^{1/\alpha}e^{-t\ell}\,dt.
\]
Tonelli's theorem, integration over $T$, and the substitution $t=r^\alpha$ now yield
\[
 \int_0^\infty e^{-ps}\widehat B_\alpha(s,b)\,ds
 =\lambda_\kappa K_{\rm ext}\int_0^\infty r^\alpha\frac{p-r}{p^\alpha-r^\alpha}e^{-br}\,dr,
 \qquad K_{\rm ext}=\frac{c_g\alpha E_\alpha}{\Gamma(1+1/\alpha)}.
\]
Thus $\widehat B_\alpha=\lambda_\kappa K_{\rm ext}R_\alpha$, where
\begin{equation}\label{eq:kernel-inverse}
 R_\alpha(s,b)=\mathcal L^{-1}_{p\to s}\left[\int_0^\infty r^\alpha\frac{p-r}{p^\alpha-r^\alpha}e^{-br}\,dr\right].
\end{equation}
Here $\mathcal L^{-1}_{p\to s}$ denotes inverse Laplace transform. With the principal branch of $p^\alpha$, the Stieltjes density obtained from the boundary values across the negative real axis is
\[
 \begin{aligned}
 j_\alpha(x,r)&=-\frac1\pi\operatorname{Im}\frac{-x-r}{x^\alpha e^{i\pi\alpha}-r^\alpha}\\
 &=-\frac{\sin(\pi\alpha)}\pi\frac{(x+r)x^\alpha}{x^{2\alpha}-2\cos(\pi\alpha)x^\alpha r^\alpha+r^{2\alpha}},\qquad x,r>0.
 \end{aligned}
\]
Cauchy's formula gives the Stieltjes representation $(p-r)/(p^\alpha-r^\alpha)=\int_0^\infty j_\alpha(x,r)(p+x)^{-1}\,dx$: the ratio is bounded at zero, is $O(|p|^{1-\alpha})$ at infinity, and its singularity on the positive real axis is removable. Writing $(p+x)^{-1}=\int_0^\infty e^{-s(p+x)}\,ds$, setting $x=ur$ and integrating $r^2e^{-(b+su)r}$ gives
\[
 R_\alpha(s,b)=-\frac{2\sin(\pi\alpha)}\pi\int_0^\infty\frac{u^\alpha(1+u)}{(u^{2\alpha}-2\cos(\pi\alpha)u^\alpha+1)(b+su)^3}\,du.
\]
This function is positive for $s,b>0$. The inverse Laplace measure has no atom at $s=0$, since its transform tends to zero as $p\to\infty$.

\emph{Normalization.} Let $A_oL^2$ be the mean area of an ordinary quantum disk of perimeter $L$, and let $A_g\ell^{2/\alpha}$ be the mean area of a generalized quantum disk of generalized perimeter $\ell$. By~\cite[Theorem~2.23 and Proposition~2.29]{ACSWLoop},
\[
 A_o=\frac1{4(\alpha-1)\sin(\pi/\alpha)},\qquad A_g=\frac{\alpha}{4\sin(\pi/\alpha)}\frac{\Gamma(1-1/\alpha)}{\Gamma(1/\alpha)}.
\]
Removing the area weight in~\cite[Theorem~2.23]{ACSWLoop} gives the unmarked ordinary disk area an inverse-gamma distribution of shape $\alpha$. Sampling the area mark from the whole inner generalized disk changes the inner welding factor to
\[
 \ell|\mathrm{GQD}_{1,0}(\ell)|=c_gA_g\ell^{-1+1/\alpha}=\frac{c_gA_g}{\Gamma(1-1/\alpha)}\int_0^\infty t^{-1/\alpha}e^{-t\ell}\,dt.
\]
Integrate $Le^{-pL}$ against the measure of the reconstructed disk, where $L$ is its perimeter. In the cyclic annulus, multiplication by its outer perimeter $L$ differentiates the Laplace transform, and
\[
 \int_0^\infty\frac{k\,dT}{T}[-\partial_p e^{-T\Psi(p,r)}]=k\,\partial_p\log\Psi(p,r).
\]
On the welding side of~\eqref{eq:lsyz-loop-welding}, the integral of this test function is $\lambda_\kappa K_{\rm in}p^{\alpha-2}$, where
\[
 K_{\rm in}=\frac{kc_gA_g\alpha I_\alpha}{\Gamma(1-1/\alpha)},\qquad I_\alpha=\int_0^\infty u^{\alpha-2}\left(\frac\alpha{1-u^\alpha}-\frac1{1-u}\right)du=-\pi\left(\cot\frac\pi\alpha+\cot(\pi\alpha)\right)>0.
\]
The singularities at $u=1$ cancel. Setting $v=u^\alpha$ in the first term and applying the principal-value identity $\operatorname{pv}\int_0^\infty v^{q-1}(1-v)^{-1}\,dv=\pi\cot(\pi q)$ for $0<q<1$, with $q=1-1/\alpha$ and $q=\alpha-1$, gives the displayed value of $I_\alpha$. In the reconstructed disk, the outermost CLE loop surrounding the area mark touches the boundary with probability
\[
 p_{\rm touch}=1+\frac{\sin(\pi(\alpha+2/\alpha))}{\sin(\pi\alpha)}
\]
by~\cite[Theorem~1.1]{ASYZ}. This probability is independent of the marked point by conformal invariance. Since the ordinary disk perimeter measure is $RL^{-2-\alpha}\,dL$, the same integral equals $p_{\rm touch}RA_o\Gamma(2-\alpha)p^{\alpha-2}$. It follows that
\[
 \lambda_\kappa K_{\rm ext}=\frac{p_{\rm touch}RA_o\Gamma(2-\alpha)}{K_{\rm in}}K_{\rm ext}=\frac{\alpha p_{\rm touch}}{I_\alpha}=-\frac\alpha\pi\sin\frac{2\pi}\alpha=C_\alpha.
\]
These integrals are finite and positive for $1<\alpha<2$. Thus $\widehat B_\alpha=B_\alpha$, and Lemma~\ref{lem:prewelding-product} gives~\eqref{eq:decorated-cut}.
\end{proof}

The boundary-length kernel $B_\alpha$ also has the elementary form
\begin{equation}\label{eq:cut-B-elementary}
 B_\alpha(s,b)=C_\alpha\partial_s\partial_b
                   \frac{s^{\alpha-1}-b^{\alpha-1}}{s^\alpha-b^\alpha}.
\end{equation}
It is symmetric, homogeneous of degree $-3$, and has diagonal value $B_\alpha(s,s)=C_\alpha(\alpha^2-1)/(6\alpha s^3)$. Its Mellin transform is
\begin{equation}\label{eq:cut-Mellin}
 \int_0^\infty s^zB_\alpha(s,b)\,ds
 =b^{z-2}\frac{-\sin(2\pi/\alpha)\sin(\pi/\alpha)z(1-z)}
 {\sin(\pi z/\alpha)\sin(\pi(1-z)/\alpha)},
 \qquad 1-\alpha<\operatorname{Re}z<\alpha.
\end{equation}
To obtain this formula, integrate $(b+us)^{-3}s^z$ in $s$, giving $\tfrac12b^{z-2}u^{-z-1}\Gamma(1+z)\Gamma(2-z)$, then set $v=u^\alpha$ and factor the quadratic denominator. For fixed $b>0$, the kernel has order $s^{\alpha-2}$ at zero and $s^{-\alpha-1}$ at infinity, so the strip of absolute convergence is precisely $1-\alpha<\operatorname{Re}z<\alpha$. The formula is understood by continuity at removable singularities.

\section{The renewal equation on the Brownian annulus}\label{sec:welding}

We derive an integro-differential equation from the Mellin transform of the Brownian annulus measure weighted by the averaged two-point observable $\overline G$. Throughout this section and through Section~\ref{sec:identification}, $\kappa=6$, $\gamma=\sqrt{8/3}$ and $B=B_{3/2}$, where $B_\alpha$ is defined in~\eqref{eq:cut-B-positive}.

Let $\mathrm{BA}$ denote the Brownian annulus measure of~\cite[Definition~1.1 and equation~(1.1)]{ARSAnnulus} (see also~\cite{BettinelliMiermontII,LeGallMetzDonnadieu} for other equivalent constructions). Write $L,m>0$ for the outer and inner quantum boundary lengths, and let $x\in(0,1)$ be the normalized quantum boundary length of the counterclockwise arc from the first outer mark to the second. We choose a homogeneous density $z(L,x;m)$ with respect to $dL\,dm$ for $\mathrm{BA}$ weighted by the averaged two-point observable $\overline G(x)$ of Section~\ref{sec:geometry}, as follows. With $r=m/L$, the boundary-length distribution of $\mathrm{BA}$ is
\[
             \frac{dL}{L}\,\frac{dr}{2\sqrt r(1+r)}.
\]
Adding a constant to the field rescales both lengths and leaves $\overline G(x)$ unchanged. On the quotient by this scaling, let $\mathcal Z(r,x)$ be the conditional mean of $\overline G(x)$ given $r$. The averaging identity in Lemma~\ref{geo:allocation} makes $\mathcal Z(r,\cdot)$ nonnegative, symmetric, concave and bounded by one. Construct it first at rational $x$, choose a common full-measure set of $r$, and extend continuously in $x\in(0,1)$. Set it to zero on the exceptional set of $r$, and define
\[
 z(L,x;m)=\frac{L^{-2}\mathcal Z(m/L,x)}
                         {2\sqrt{m/L}(1+m/L)}.
\]
This density satisfies
\begin{equation}\label{weld:densitybound}
 \begin{gathered}
 z(cL,x;cm)=c^{-2}z(L,x;m),\qquad c>0,\\
 0\le z(L,x;m)\le\frac1{2\sqrt{Lm}(L+m)}.
 \end{gathered}
\end{equation}
For real $\sigma$, define the Mellin transform
\begin{equation}\label{eq:physical-F}
 F(\sigma)=\int_0^\infty m^{-\sigma}z(1,1/2;m)\,dm\in[0,\infty].
\end{equation}
Recall $x(t)=e^t/(1+e^t)$, $a(t)=1-x(t)$ and the operators $L_0,T_s$ from~\eqref{eq:intro-kernel}--\eqref{eq:intro-T}. For real $s\in(-2,2)$ and nonnegative measurable $f$, the integral in~\eqref{eq:intro-T} also defines $T_sf$, with values in $[0,\infty]$.

\begin{proposition}[Integro-differential equation]\label{weld:transform}
For $\sigma\in(-2,2)$ such that $F(\sigma)<\infty$, define
\begin{equation}\label{weld:Mellin}
 \varphi_{-\sigma}(x)=\int_0^\infty m^{-\sigma}z(1,x;m)\,dm,
 \qquad
 u_\sigma(t)=\frac{\varphi_{-\sigma}(x(t))}{\sqrt{x(t)a(t)}}.
\end{equation}
Then $u_\sigma$ is continuous, nonnegative and even, with $F(\sigma)=u_\sigma(0)/2$. There is a nonnegative even locally finite measure $j_\sigma$ such that, in distributions on $\mathbb R$,
\begin{equation}\label{weld:operatoridentity}
                    (L_0-T_{-\sigma})u_\sigma=j_\sigma.
\end{equation}
The measure $j_\sigma$ is the contribution of white clusters meeting the inner boundary, as specified in~\eqref{eq:j-positive-definition}.
\end{proposition}

\subsection{Preservation of the filled metric ball}\label{subsec:metric-cutoff}

Following~\cite[Section~5]{SXZ}, we use the Brownian annulus to transfer the measure decomposition of Proposition~\ref{prop:decorated-cut} to the annulus setup. The construction of the Brownian annulus removes a filled metric ball and reweights by a function of its perimeter. The same weight applies before and after cutting a loop disjoint from the ball; see Figure~\ref{fig:safe-ball}.

Sample $(D,h,z)$ from $\mathrm{QD}_{1,0}$, let $d_h$ be its LQG metric, and consider
\[
 E_D=\{d_h(z,\partial D)>1\}.
\]
On this event let $K_D$ be the closed filled metric ball of radius one about $z$, obtained by adjoining the complementary components not meeting $\partial D$, and write $m_D$ for its quantum boundary length. For a subdomain $U$ containing $z$, define $E_U,K_U,m_U$ in the same way using the restricted field and the internal LQG metric on $U$. The CLE is sampled independently of the field.
\begin{figure}[!htb]
\centering
\begingroup
\definecolor{xiBlue}{RGB}{224,242,247} \definecolor{xiRose}{RGB}{251,229,233} \definecolor{xiGold}{RGB}{255,244,208} \definecolor{xiInk}{RGB}{20,35,48}
\begin{tikzpicture}[x=1cm,y=1cm,line width=.7pt,draw=black,
  every node/.style={font=\fontsize{9}{10.5}\selectfont,inner sep=1.5pt},
  >={Latex[length=1.6mm,width=1.2mm]}]
\path[use as bounding box] (0,0) rectangle (15.55,7.15);
\node at (3.40,6.85) {Before the cut: $(D,h,z)$};
\node at (11.80,6.85) {After the cut: $(U,h|_U,z)$};
\foreach \shift in {0,8.40}{
  \begin{scope}[shift={(\shift,0)}]
    \coordinate (bp) at ({3.40+2.50*cos(134)},{3.80+2.50*sin(134)});
    \coordinate (bq) at ({3.40+2.50*cos(226)},{3.80+2.50*sin(226)});
    \path[draw=black,fill=white] (bp) arc[start angle=134,end angle=-134,radius=2.50]
      .. controls (2.35,2.40) and (2.87,2.96) .. (2.62,3.80)
      .. controls (2.87,4.64) and (2.35,5.20) .. (bp) -- cycle;
    \path[draw=black,fill=xiBlue,draw=none] (4.03,3.80) ellipse (1.12 and 1.15);
    \draw[densely dashed] (4.03,3.80) ellipse (1.12 and 1.15);
    \path[draw=black,fill=xiGold] (3.41,3.39)
      .. controls (3.22,3.81) and (3.43,4.29) .. (3.81,4.47)
      .. controls (4.28,4.67) and (4.62,4.36) .. (4.66,3.92)
      .. controls (4.91,3.64) and (4.52,3.12) .. (4.10,3.14)
      .. controls (3.80,3.03) and (3.59,3.18) .. (3.41,3.39) -- cycle;
    \fill (4.02,3.81) circle (.05);
    \node[above right] at (4.02,3.81) {$z$};
    \node at (3.85,3.44) {$K$};
    \node[above right] at (4.52,4.30) {$m$};
    \node[right] at (5.13,3.50) {$\rho$};
    \node at (3.70,5.85) {inherited arc $a$};
    \fill (bp) circle (.05);
    \fill (bq) circle (.05);
    \node[above left] at (bp) {$p_C$};
    \node[below left] at (bq) {$q_C$};
  \end{scope}
}
\coordinate (sp) at ({3.40+2.50*cos(134)},{3.80+2.50*sin(134)});
\coordinate (sq) at ({3.40+2.50*cos(226)},{3.80+2.50*sin(226)});
\path[draw=black,fill=xiRose] (sp) arc[start angle=134,end angle=226,radius=2.50]
  .. controls (2.35,2.40) and (2.87,2.96) .. (2.62,3.80)
  .. controls (2.87,4.64) and (2.35,5.20) .. (sp) -- cycle;
\fill (sp) circle (.05);
\fill (sq) circle (.05);
\node[align=center] at (1.73,3.77) {discarded\\side of $A$};
\node[right] at (2.66,4.38) {$b$};
\node[left,align=center] at (11.06,3.78) {newly exposed\\arc $b$};
\node[left] at (.89,3.78) {$s$};
\draw[->] (6.35,3.78) -- (8.41,3.78);
\node[above,align=center] at (7.38,3.87) {cut along the\\frontier of $A$};
\node[below,align=center] at (7.38,3.61) {restrict $h$\\to $U$};

\node at (3.40,1.03) {$L=a+s$};
\node at (11.80,1.03) {$L'=a+b$};
\node at (7.77,.48) {$1<\rho<d_{h|_U}(z,\partial U),\qquad
                 K_D=K_U=K,\qquad m_D=m_U=m$};
\end{tikzpicture}
\endgroup
\caption{Cutting a complete boundary-touching loop outside the filled ball. The planar image $A=\Pi_C(X_C)$ meets the original boundary and avoids the filled unit ball $K_D$. Cutting along its frontier removes the pink region and leaves the gold filled ball unchanged. The inherited boundary arc has length $a$, the newly exposed boundary arc has length $b$, and the discarded boundary arc has length $s$. The blue region shows the common metric neighborhood beyond radius one. The drawing shows the part of $\partial A$ bordering the retained component.}
\label{fig:safe-ball}
\end{figure}

\begin{lemma}\label{weld:safeball}
On $E_D$, let $C$ be a white cluster whose interface $\eta_C$ is a complete boundary-touching CLE loop, with $\Pi_C(X_C)\cap K_D=\varnothing$. Let $U$ be the component of $D\setminus\Pi_C(X_C)$ containing $z$, equipped with the restricted field.  Then $E_U$ holds and
\[
 K_U=K_D,\qquad m_U=m_D.
\]
Conversely, in the welding of Proposition~\ref{prop:decorated-cut}, $E_U$ implies $E_D$ and the condition $\Pi_C(X_C)\cap K_D=\varnothing$.
\end{lemma}

\begin{proof}
Write $A=\Pi_C(X_C)$, and by Section~\ref{subsec:necklaces}, it is compact, connected and meets $\partial D$. By locality of the LQG metric~\cite[Section~1.2 and Remark~1.5]{GMMetric}, the LQG metric of $h|_U$ is the internal LQG metric $d_h(\cdot,\cdot;U)$. A first-exit argument therefore gives
\[
 d_{h|_U}(z,\partial U)=d_h(z,A\cup\partial D)=:\delta,
\]
and the unfilled balls in $D$ and $U$ agree at every radius below $\delta$. Indeed, any path reaching $A\cup\partial D$ first approaches a boundary point of the metric completion of $U$. Conversely, the inclusion of $U$ into $D$ is $1$-Lipschitz, and every boundary point of this completion projects to $A\cup\partial D$: a limit in $U$ also converges in the internal metric, by taking a small metric neighborhood contained in $U$.

If $E_D$ holds and $A\cap K_D=\varnothing$, compactness gives $\delta>1$, hence $E_U$. If $E_U$ holds, the same equality gives $E_D$ and shows that $A$ avoids the common unfilled unit ball. Since $A$ is connected and meets the original boundary, it lies in the exterior component of the ball's complement. A path from a point of $U$ through this exterior component to $\partial D$ either stays in $U$ or first reaches $A$; in either case it reaches $\partial U$ without crossing the ball. Thus precisely the same complementary components are filled in the two domains. This proves $K_U=K_D$ and $A\cap K_D=\varnothing$ in both directions.

Finally choose $1<\rho<\delta$. The common balls through radius $\rho$ give the same field on a neighborhood of the common inner boundary. The associated quantum boundary measures, and hence $m_U$ and $m_D$, agree.
\end{proof}

By~\cite[Proposition~6.10]{ARSAnnulus}, there are a positive constant $c_2$ and a positive function $\mathfrak f:(0,\infty)\to(0,\infty)$ such that removing $K_D$ under the measure
\begin{equation}\label{weld:BAsurgery}
 \frac{1_{E_D}}{c_2m_D\mathfrak f(m_D)}\,
 \mathrm{QD}_{1,0}
\end{equation}
gives the Brownian annulus measure $\mathrm{BA}$, after fixing the common normalization of $\mathrm{QD}_{1,0}$. Since the weight depends only on $m_D$, Lemma~\ref{weld:safeball} allows us to use the same weight on both sides of the measure decomposition in Proposition~\ref{prop:decorated-cut}.

Condition on the field and use the inner-boundary approximation of Section~\ref{subsec:necklaces}. Removing $K_D$ and imposing black boundary conditions on the newly created inner boundary erases precisely the white connections through $K_D$: every compact connecting path avoiding $K_D$ lies in one of the approximating annuli. The white clusters in $D\setminus K_D$ are therefore the path components in~\eqref{eq:geo:induce}. Conformal invariance at each approximation, followed by its limit, gives the usual percolation law on the annulus, conditionally independent of the field given its conformal modulus. Let $\mathcal C_\tau=[0,\tau]\times(\mathbb R/\mathbb Z)$ be the cylinder of circumference one and modulus $\tau>0$. Write $\mathrm{LF}_\tau$ for the measure obtained from a free boundary Gaussian free field on this cylinder by integrating its additive constant against Lebesgue measure; this measure is independent of the choice of additive normalization of the Gaussian free field. Let $\eta$ denote the Dedekind eta function,
\[
 \eta(z)=e^{\pi iz/12}\prod_{n=1}^\infty(1-e^{2\pi inz}),\qquad \operatorname{Im}z>0.
\]
The Brownian annulus identity~\cite[Theorem~1.3]{ARSAnnulus} is
\begin{equation}\label{weld:BAfield}
 \mathrm{BA}
 =\frac1{\sqrt2}\eta(2i\tau)\,
 \mathrm{LF}_{\tau}(dh)\,1_{\{\tau>0\}}\,d\tau.
\end{equation}

\subsection{The renewal equation}
\label{subsec:annulus-renewal}

We apply Proposition~\ref{prop:geometric-cutting} to the Brownian annulus. Averaging the original marks gives the piecewise linear weights $w_\ell$, $\ell\in(0,1)$, of Lemma~\ref{geo:allocation}, and integrating the first boundary mark of the retained component contributes its perimeter.

For a nonnegative measurable symmetric function $f$ on $(0,1)$, put
\[
 z_f(L;m)=\int_0^1 f(x)z(L,x;m)\,dx,
 \qquad
 \Psi_f(\ell)=\int_0^1
 f(x)\min\{x,1-x,\ell,1-\ell\}\,dx.
\]
Let $t_f(L;m)$ be the density, with respect to $dL\,dm$, obtained by integrating the first term of~\eqref{eq:geo:decomposition} against $\mathrm{BA}$. It includes precisely the contributions from white clusters that meet the inner boundary. Choose it on the quotient by scaling as for $z_f$, with the same homogeneity.

\begin{proposition}[Renewal equation]\label{weld:renewal}
As nonnegative measures in $(L,m)$,
\begin{equation}\label{weld:renewaleq}
\begin{split}
 z_f(L;m)={}&t_f(L;m)\\
 &+\int_0^L\!ds\int_0^\infty\!db\,
 (L-s+b)B(s,b)\Psi_f(s/L)
 z\left(L-s+b,\frac{L-s}{L-s+b};m\right).
\end{split}
\end{equation}
For the densities defined above, the identity holds almost everywhere in $(L,m)$, with values in $[0,\infty]$. Homogeneity also gives the identity for almost every $m$ at $L=1$.
\end{proposition}

\begin{proof}
Apply the averaged cutting identity of Proposition~\ref{prop:geometric-cutting} with $K=K_D$ and $\mu=\nu_h$. The clusters meeting the inner boundary give $t_f$. Each remaining term has weight $\Psi_f(s_C/L)$, where $s_C=\nu_h(S_C)$, and evaluates $G$ in the component retained after cutting $\eta_C$.

The retained component is marked at the endpoints $p_C,q_C$ of its inherited boundary arc. The factor $\Psi_f(s_C/L)$ comes from averaging the two original marks.

Integrate the loop sum using the measure decomposition in Proposition~\ref{prop:decorated-cut}. By Lemma~\ref{weld:safeball}, the original and retained components have the same inner perimeter $m$, so the factor in \eqref{weld:BAsurgery} is unchanged.  Write $a$ and $b$ for the lengths of the inherited and newly exposed boundary arcs. At fixed perimeter $L'=a+b$ of the retained component, integrating the first ordered boundary mark gives $L'$, and the second is specified by the arc length $a$.  Thus the $(a,b)$-disintegration of the marked retained component has density $L'z(L',a/L';m)$.  Finally set $a=L-s$ and $L'=L-s+b$.  This gives \eqref{weld:renewaleq} for integrable $f$ by Tonelli's theorem. For an arbitrary nonnegative measurable symmetric $f$, apply this identity to $f\wedge n$ and let $n\to\infty$ by monotone convergence.
\end{proof}

\subsection{The Mellin transform}\label{subsec:mellin-sign}

We prove Proposition~\ref{weld:transform} from~\eqref{weld:renewaleq} by taking the Mellin transform and differentiating twice with respect to $x$ in the sense of distributions.

\begin{proof}[Proof of Proposition~\ref{weld:transform}]
Symmetry and concavity give $\mathcal Z(r,x)\le\mathcal Z(r,1/2)$. Hence finiteness of $F(\sigma)$ and dominated convergence imply that $\varphi_{-\sigma}$ is continuous, nonnegative and symmetric about $1/2$. This proves the stated properties of $u_\sigma$, and its value at zero is $2F(\sigma)$. Homogeneity gives
\[
 \int_0^\infty m^{-\sigma}z(L,x;m)\,dm
       =L^{-1-\sigma}\varphi_{-\sigma}(x).
\]
Let $\vartheta_\sigma$ be the Mellin transform of the contribution from clusters meeting the inner boundary, defined by
\[
 \int_0^1 f(x)\vartheta_\sigma(x)\,dx
 =\int_0^\infty m^{-\sigma}t_f(1;m)\,dm.
\]
As a nonnegative integral of the symmetric concave weights $w_\ell$, $\vartheta_\sigma$ is symmetric and concave, with $0\le\vartheta_\sigma\le\varphi_{-\sigma}$. Hence $-\vartheta_\sigma''$ is a nonnegative locally finite measure on $(0,1)$.

Apply the Mellin transform to~\eqref{weld:renewaleq}, set $L=1$ and differentiate twice in distributions. The two slope changes of the piecewise linear weight $w_\ell$ give
\[
 \begin{split}
 -\varphi_{-\sigma}''(x)={}&-\vartheta_\sigma''(x)\\
 &+\int_0^\infty (1-x+b)^{-\sigma}B(x,b)
                \varphi_{-\sigma}\!\left(\frac{1-x}{1-x+b}\right)db\\
 &+\int_0^\infty (x+b)^{-\sigma}B(1-x,b)
                \varphi_{-\sigma}\!\left(\frac{x}{x+b}\right)db .
 \end{split}
\]
Truncate the arc sums, differentiate the weights $w_\ell$ and apply monotone convergence. Concavity and finiteness of $\varphi_{-\sigma}$ make all three measures locally finite. The factor $(1-x+b)^{-\sigma}$ comes from multiplying the perimeter factor $L'$ in~\eqref{weld:renewaleq} by the factor $L'^{-1-\sigma}$ from the Mellin transform.

Use $B(s,b)=(sb)^{-3/2}H(\log(s/b))$, where $H$ is given in~\eqref{eq:intro-kernel}. The substitutions $b=(1-x)(1-y)/y$ and $b=x(1-y)/y$, together with symmetry of $\varphi_{-\sigma}$ in the first term, turn these two integrals into $T_{-\sigma}u_\sigma$. For a smooth function $\varphi$, the differential part follows from
\[
 L_0\left[\frac{\varphi(x(t))}{\sqrt{x(t)a(t)}}\right]
 =-[x(t)a(t)]^{3/2}\varphi''(x(t)).
\]
Define $j_\sigma$ by requiring that, for every continuous compactly supported test function $\chi$,
\[
 \int_{\mathbb R}\chi(t)\,j_\sigma(dt)
 =\int_0^1\sqrt{x(1-x)}\,
 \chi\!\left(\log\frac{x}{1-x}\right)(-\vartheta_\sigma'')(dx).
\]
It is nonnegative, even and locally finite. This gives~\eqref{weld:operatoridentity}. In particular, an arc of normalized quantum boundary length $\ell$ contributes atoms of mass $\sqrt{\ell(1-\ell)}$ at $t=\pm\log(\ell/(1-\ell))$, counted twice when $\ell=1/2$.
\end{proof}

\section{Estimates for the Mellin transform}
\label{sec:endpoint}

We prove boundedness of $u_\sigma$, analyticity of $j_\sigma$, and comparison of the expectation in~\eqref{eq:weighted-selector} below with the unmarked three-arm probability. The proofs use GMC moment estimates and the boundary one-arm and polychromatic four-arm exponents. Throughout, $\kappa=6$ and $\gamma^2=8/3$.

Consider the cylinder $\mathcal C_\tau=[0,\tau]\times(\mathbb R/\mathbb Z)$, whose outer and inner boundaries are $\partial_o\mathcal C_\tau=\{0\}\times(\mathbb R/\mathbb Z)$ and $\partial_i\mathcal C_\tau=\{\tau\}\times(\mathbb R/\mathbb Z)$, respectively. We transport the geometric notation of Section~\ref{sec:geometry} conformally to $\mathcal C_\tau$. By conformal invariance, $p_3(\tau)$ is the probability of three disjoint black crossings of this cylinder. Let $h^{\mathrm{lat}}$ be the lateral part of the free boundary Gaussian free field, that is, its component with zero average on each circle $\{t\}\times(\mathbb R/\mathbb Z)$, and let $M_o,M_i$ be its outer and inner quantum boundary length measures. Set $Z_o=M_o(\partial_o\mathcal C_\tau)$ and $Z_i=M_i(\partial_i\mathcal C_\tau)$. For $0<x<1$, sample the first mark according to $M_o/Z_o$ and place the second mark so that the counterclockwise arc between them has quantum boundary length $xZ_o$. We denote the average of the two-point observable $G$ over the first mark by $\overline G_\tau(x)$ and, for a real parameter $\sigma$, put
\begin{equation}\label{eq:weighted-selector}
 W_\sigma(\tau,x)=\mathbb E\big[(Z_o/Z_i)^\sigma
                 \overline G_\tau(x)\big].
\end{equation}
The expectation is over the field and independent percolation, with $\tau$ fixed.

Let $\mathcal M_i$ be the set of white clusters $C\in\mathcal M$ whose images $\Pi_C(X_C)$ meet $\partial_i\mathcal C_\tau$. Set $\mathcal M_i=\varnothing$ off $H_3$, and for $C\in\mathcal M_i$ put
\[
 \ell_C=\frac{M_o(S_C)}{Z_o}\in(0,1),\qquad
 t_C=\log\frac{\ell_C}{1-\ell_C}.
\]
For every real $\sigma$ and every nonnegative measurable function $\chi$ on $\mathbb R$, define
\begin{equation}\label{eq:j-positive-definition}
\begin{split}
 \int_{\mathbb R}\chi(t)\,j_\sigma(dt)
 ={}&\frac{\sqrt2}{\gamma}\int_0^\infty
       \eta(2i\tau)e^{2\pi\sigma^2\tau/3}\\
 &\quad\times\mathbb E\!\left[(Z_o/Z_i)^\sigma\mathbf1_{H_3}
       \sum_{C\in\mathcal M_i}\sqrt{\ell_C(1-\ell_C)}
       \big(\chi(t_C)+\chi(-t_C)\big)\right]d\tau.
\end{split}
\end{equation}
This defines a nonnegative even measure with values in $[0,\infty]$. Integrating the additive constant and circle-average process as in~\eqref{eq:zero-mode-endpoint}, and then using~\eqref{weld:BAfield}, shows that it agrees with the measure in Proposition~\ref{weld:transform} for $\sigma\in(-2,2)$ with $F(\sigma)<\infty$. Recall the Green kernel $g(t)=e^{-|t|/2}$ of $L_0$, and set
\begin{equation}\label{eq:xi-sigma-conversion-endpoint}
 X(\sigma)=\frac{\sigma^2-1/4}{3}.
\end{equation}

\begin{proposition}\label{est:endpoint}
For every compact interval $J\subset(3/2,2)$ there are constants $C,c,K>0$ and $\tau_0>c$ such that, uniformly for $\sigma\in J$ and $0<x<1$,
\begin{align}
 W_\sigma(\tau,x)&\le C\sqrt{x(1-x)}p_3(\tau-c),
                        &&\tau\ge\tau_0,\label{eq:long-endpoint}\\
 W_\sigma(\tau,x)&\le C\tau^{-K}\sqrt{x(1-x)},
                        &&0<\tau\le\tau_0.\nonumber
\end{align}
Consequently, if $\sigma\in J$ and $X(\sigma)<\xi_3$, the function $\varphi_{-\sigma}(x)/\sqrt{x(1-x)}$ is bounded on $(0,1)$, or equivalently $u_\sigma$ is bounded on $\mathbb R$.
\end{proposition}

\begin{proposition}\label{est:inhomogeneous}
For every compact interval $J\subset(3/2,2)$ there is $a>0$ such that $\sigma\mapsto j_\sigma$ extends holomorphically to a complex neighborhood of $J$, with values in the Banach space of complex Borel measures $j$ on $\mathbb R$ with norm $\int e^{a|t|}|j|(dt)<\infty$, where $|j|$ denotes total variation. For real $\sigma\in J$, the measure is even, nonnegative and nonzero. Consequently $g*j_\sigma$ is a strictly positive continuous function that decays exponentially at both ends.
\end{proposition}

\begin{proposition}\label{est:comparison}
For every compact interval $J\subset(3/2,2)$, there are constants $c_J,C_J>0$ and $c,c'>0$ such that, for all sufficiently large $\tau$ and $\sigma\in J$,
\begin{equation}\label{eq:marked-unmarked-comparison}
 c_Jp_3(\tau-c)
       \le W_\sigma(\tau,1/2)
       \le C_Jp_3(\tau-c').
\end{equation}
The right endpoint of the interval of convergence of the Mellin transform $F$ is therefore
\begin{equation}\label{eq:physical-abscissa}
       \sigma_c=\sqrt{1/4+3\xi_3}\in(3/2,2).
\end{equation}
\end{proposition}

We prove these propositions in Section~\ref{subsec:proof-estimates}. Below, we first provide the needed technical input from GMC and percolation estimates.

\subsection{GMC estimates}\label{subsec:boundary-chaos}

\paragraph{Integration of the additive constant.}\label{subsec:gaussian-zero-mode}

In the representation of the Brownian annulus \eqref{weld:BAfield}, decompose the field as $h=\mathsf B+h^{\mathrm{lat}}+c$. The circle average $\mathsf B_t$ on $\{t\}\times(\mathbb R/\mathbb Z)$ is a Brownian motion in $t\in[0,\tau]$, with $\mathsf B_0=0$ and $\operatorname{Var}(\mathsf B_t)=2\pi t$. It is independent of $h^{\mathrm{lat}}$, and the additive constant $c$ is integrated against Lebesgue measure. The two boundary lengths are
\[
 L=e^{\gamma c/2}Z_o,
 \qquad m=e^{\gamma c/2+\gamma\mathsf B_\tau/2}Z_i.
\]
For a nonnegative measurable function $A$ of the lateral field and percolation, let $\mathrm{LF}_\tau[A;L,m]$ be the joint boundary-length density obtained by integrating $A$ against $\mathrm{LF}_\tau$ and the independent percolation law. Tonelli's theorem and the substitution $e^{\gamma c/2}=L/Z_o$ give
\begin{equation}\label{eq:zero-mode-endpoint}
 \int_0^\infty m^{-\sigma}\mathrm{LF}_\tau[A;L,m]\,dm
 =\frac2\gamma L^{-\sigma-1}
   e^{2\pi\sigma^2\tau/3}
   \mathbb E[(Z_o/Z_i)^\sigma A].
\end{equation}
The Gaussian expectation is $\mathbb E[e^{-\gamma\sigma\mathsf B_\tau/2}]=e^{2\pi\sigma^2\tau/3}$. Together with $\eta(2i\tau)\asymp e^{-\pi\tau/6}$ as $\tau\to\infty$, this gives the exponent $X(\sigma)$ in \eqref{eq:xi-sigma-conversion-endpoint}.

\paragraph{Mixed GMC boundary-length moments.}

The next lemma estimates the piecewise linear weight $w_\ell(x)=\min\{x,1-x,\ell,1-\ell\}$, summed over scales with the boundary one-arm factor $e^{2n/3}$. For this lemma, absorb $\gamma/2$ into the boundary field, so that its logarithmic covariance coefficient is $4/3$, and write the moment scaling exponent of the boundary GMC measure as
\[
 \zeta(r)=\frac{5r-2r^2}{3}.
\]
The constants below are uniform for a uniformly bounded regular part of the covariance.

\begin{lemma}\label{lem:mixed-chaos}
Let $M$ be the boundary GMC measure on a circle, with the covariance normalization just specified. For an arc $I$, write $I^c$ for its complement in the circle. If $0<p,q<3/2$ and $p+q<5/2$, then
\begin{equation}\label{eq:mixed-chaos}
 \sup_I\mathbb E[M(I)^pM(I^c)^q]<\infty,
\end{equation}
where the supremum is over all arcs. If $3/2<\sigma<2$, $0<x\le1/2$, and $I_n$ is any arc whose length is comparable to $e^{-n}$, then
\begin{equation}\label{eq:chaos-square-root}
 \sum_{n\ge0}e^{2n/3}
 \mathbb E\!\left[(M(I_n)+M(I_n^c))^\sigma
 \min\left\{x,\frac{M(I_n)}{M(I_n)+M(I_n^c)},
                 \frac{M(I_n^c)}{M(I_n)+M(I_n^c)}\right\}\right]
 \le C_\sigma\sqrt{x}.
\end{equation}
The constant is locally uniform in $\sigma\in(3/2,2)$, with the comparability constants for $|I_n|$ fixed.
\end{lemma}

\begin{proof}
Decompose the two complementary arcs into Whitney intervals near their common endpoints. Let $I_n,J_m$ be such intervals of lengths comparable to $e^{-n}$ and $e^{-m}$, with $n\ge m$. Their separation is at least a constant times $e^{-m}$. Comparing their fields with a common Gaussian part at scales larger than $e^{-m}$ and independent parts at smaller scales gives
\begin{equation}\label{eq:whitney-moment}
 \mathbb E[M(I_n)^pM(J_m)^q]
 \le C e^{-m\zeta(p+q)}e^{-(n-m)\zeta(p)}.
\end{equation}
Moments on a fixed interval are finite for powers below $3/2$ by \cite[equation~(1.6) and Theorem~1.1]{Remy} and bounded covariance comparison. For more distant intervals, increasing the cross-covariance to that of the comparison model increases the mixed moment, since the mixed derivative in the two masses is positive. Bounded changes of the covariances within either interval change the estimate by at most a constant factor. Bound powers of sums of masses of the Whitney intervals by sums of their powers with weights $C_\epsilon e^{\epsilon n}$. Since $\zeta(p)$, $\zeta(q)$ and $\zeta(p+q)$ are positive, sufficiently small $\epsilon>0$ makes the double series converge. This proves \eqref{eq:mixed-chaos}. For coordinate intervals, take a comparison circle whose circumference is at least twice their total span.

For the second assertion, fix $\sigma\in(3/2,2)$ and $0<x\le1/2$, and for two masses $A,B>0$ set
\[
 F_x(A,B)=(A+B)^\sigma
              \min\{x,A/(A+B),B/(A+B)\}.
\]
This function is increasing in each mass and homogeneous of degree $\sigma$. For covariance comparisons, we introduce the smooth function
\[
 \widetilde F_x(A,B)
     =\frac{(A+B)^\sigma}{x^{-1}+2+A/B+B/A}.
\]
Note that it satisfies $F_x/3\le\widetilde F_x\le F_x$, and
\[
 |A^iB^j\partial_A^i\partial_B^j\widetilde F_x|
       \le C_\sigma\widetilde F_x,\qquad i,j\in\mathbb Z_{\ge0},\quad i+j\le2.
\]
Gaussian interpolation at fixed GMC regularization and Gronwall's inequality therefore compare $\mathbb E F_x$ up to a constant independent of $x$, whenever the Gaussian covariance kernels differ uniformly by a bounded amount. We will also use this comparison after multiplying by a negative power of a third mass whose underlying Gaussian field has bounded cross-covariances with the first two fields.

Let $\beta$ be a standard Brownian motion, with $\operatorname{Var}(\beta_t)=t$. The contribution from the additional scales of a small arc is represented by
\[
 Y_t=-\frac53t+\sqrt{\frac43}\beta_t,\qquad t\ge0.
\]
The factor $e^{2n/3}$ in the arc count gives the weighted occupation measure
\begin{equation}\label{eq:rewarded-potential}
 \int_0^\infty e^{2t/3}\mathbb P(Y_t\in dy)\,dt
 =\mathbf1_{\{y<0\}}e^{-y/2}\,dy
  +\mathbf1_{\{y>0\}}e^{-2y}\,dy
\end{equation}
by integrating the Gaussian density. For $r>0$, splitting the integral at $y=\log(x/r)$ gives
\begin{equation}\label{eq:truncated-potential}
 \int_0^\infty e^{2t/3}\mathbb E[\min\{x,re^{Y_t}\}]\,dt
       \le4\sqrt{xr}.
\end{equation}
Let $h_x(r)=(1+r)^\sigma\min\{x,r/(1+r),1/(1+r)\}$. On $r\le1$ this is bounded by a constant times $\min\{x,r\}$. On $r>1$ it is bounded by $C_\sigma\sqrt{x}\,r^{\sigma-1/2}$. Since $\zeta(\sigma-1/2)>2/3$, we obtain
\begin{equation}\label{eq:homogeneous-potential}
 \int_0^\infty e^{2t/3}\mathbb E[h_x(re^{Y_t})]\,dt
 \le C_\sigma\sqrt{x}
                  (r^{1/2}+r^{\sigma-1/2}).
\end{equation}
The same bound holds for the discrete-time sum: for $0\le s\le1$, the independent increment $Y_{n+s}-Y_n$ lies in $[-1,1]$ with uniformly positive probability, and $h_x(e^vr)\asymp h_x(r)$ uniformly for $|v|\le1$. Integration over $s$ therefore controls each summand.

Partition $I_n^c$ into groups at distance comparable to $e^{-m}$ from $I_n$ for $0\le m<n$, and a final group at distance $O(e^{-n})$, indexed by $m=n$. Each group is a union of a bounded number of arcs; write its mass as $B_m$. For $\epsilon>0$, put $S_\epsilon=\sum_{m\ge0}e^{-\epsilon m}$. The inequalities $\sum_m B_m\le S_\epsilon\sup_m e^{\epsilon m}B_m$ and the monotonicity and homogeneity of $F_x$ give
\[
 F_x\!\left(A,\sum_m B_m\right)
 \le S_\epsilon^\sigma\sum_m e^{\sigma\epsilon m}F_x(A,B_m).
\]
For separated scales, compare the rescaled fields with a common Gaussian contribution down to scale $e^{-m}$, an additional independent contribution $Y_{n-m}$ in the smaller arc, and chaos masses $U,V$ on intervals of fixed size. The covariance error is uniformly bounded, so the preceding interpolation estimate applies. Integrating the common Gaussian part gives the factor $e^{-m\zeta(\sigma)}$. For a fixed integer $c$ large enough to separate these scales, \eqref{eq:homogeneous-potential} gives
\[
 \begin{split}
 &\sum_{n\ge m+c}e^{2n/3}\,
       \mathbb E F_x\big(M(I_n),B_m\big)\\
 &\qquad\le C_\sigma\sqrt{x}\,
    e^{-m[\zeta(\sigma)-2/3]}
    \mathbb E\!\left[U^{1/2}V^{\sigma-1/2}
                    +U^{\sigma-1/2}V^{1/2}\right],
 \end{split}
\]
The powers on the right are below $3/2$, so \eqref{eq:mixed-chaos} bounds these moments uniformly. The bounded number of neighboring scales is controlled by the same mixed-moment bound and the inequality
\begin{equation}\label{eq:tent-square-root}
 F_x(A,B)\le C_\sigma\sqrt{x}\,
       \big(A^{1/2}B^{\sigma-1/2}
                    +A^{\sigma-1/2}B^{1/2}\big).
\end{equation}
Including the factor $e^{\sigma\epsilon m}$ leaves the convergent series
\[
 \sum_{m\ge0}
   e^{-[\zeta(\sigma)-2/3-\sigma\epsilon]m}.
\]
For $\sigma$ in a fixed compact interval $J\subset(3/2,2)$, choose $\epsilon>0$ with $\epsilon\max J<\inf_{\sigma\in J}[\zeta(\sigma)-2/3]$. This proves \eqref{eq:chaos-square-root} with the stated local uniformity. Apply the interpolation first to finitely many arcs and regularized chaos, and the moment bounds justify the limits.
\end{proof}

\paragraph{The two boundary fields at small modulus.}\label{subsec:thin-cylinders}

For small modulus, decompose the boundary fields into a common Gaussian field and independent remainders. Rescale the circumference to $2\pi$ and put $T=2\pi\tau$. Let $h_o,h_i$ be the outer and inner boundary traces of $h^{\mathrm{lat}}$, and let $K_{ab}(u)$ be the covariance between traces $a,b\in\{o,i\}$ at angular separation $u$. Their covariances are
\begin{align*}
 K_{oo}(u)=K_{ii}(u)
      &=2\sum_{n\ge1}\frac{\coth(nT)}n\cos(nu),\\
 K_{oi}(u)&=2\sum_{n\ge1}\frac{\operatorname{csch}(nT)}n\cos(nu).
\end{align*}
Since $\coth v-\operatorname{csch}v=\tanh(v/2)$, we have the independent decomposition
\begin{equation}\label{eq:thin-decomposition}
 h_o=C_T+E_o,\qquad h_i=C_T+E_i,
\end{equation}
where $C_T,E_o,E_i$ are centered independent Gaussian fields. The covariance coefficients of $C_T$ are $2\operatorname{csch}(nT)/n$, and those of $E_o$ and $E_i$ are $2\tanh(nT/2)/n$. For an arc $J$, write $\operatorname{osc}_J C_T=\sup_J C_T-\inf_J C_T$. For every $a,c>0$, the oscillations of $C_T$ on arcs of length $O(T)$ satisfy
\begin{equation}\label{eq:common-oscillation}
 \sup_{|J|\le cT}\mathbb E\exp\{a\operatorname{osc}_J C_T\}
       \le C_{a,c}.
\end{equation}
The increment variance is at most $C|u-v|^2/T^2$, so the Gaussian supremum bound after rescaling proves~\eqref{eq:common-oscillation}. Splitting the Fourier series at $n=T^{-1}$ shows that, after rescaling an arc of length $O(T)$ to unit size, each of $E_o,E_i$ has logarithmic covariance with uniformly bounded regular part. Its covariance with an arc at distance at least a fixed multiple of $T$ is uniformly bounded.

Partition the boundary circle into arcs $J_j$ of length comparable to $T$, choose $u_j\in J_j$, and write
\[
 a_j=\exp\{\gamma C_T(u_j)/2\}.
\]
On $J_j$, the two chaos measures equal $a_j$ times the respective chaos measures associated with $E_o$ and $E_i$, up to factors between $e^{-\gamma D_j/2}$ and $e^{\gamma D_j/2}$, where $D_j=\operatorname{osc}_{J_j}C_T$. Using the boundary-length regularization $\varepsilon^{\gamma^2/4}e^{\gamma h_\varepsilon/2}\,du$, where $h_\varepsilon$ is the field smoothed at angular scale $\varepsilon$, the common deterministic scale factor of these measures is $c_T=T^{1+\gamma^2/4}=T^{5/3}$. Rescaling $u=Tv$ accounts for both powers of $T$. If $C_j=M_i(J_j)$, then for $p,q>0$,
\begin{equation}\label{eq:denominator-cancellation}
             Z_i^{p+q}\ge C_j^pC_k^q.
\end{equation}
For factors of degree $p,q$ in the outer masses on $J_j,J_k$, this lower bound cancels $a_j^pa_k^q$ and $c_T^{p+q}$ before expectation. The fields $E_o$ and $E_i$ are independent. Lemma~\ref{lem:mixed-chaos} controls the outer positive mixed moments, while the inner masses have all negative moments by \cite{RobertVargas} and bounded covariance comparison. The remaining factors depend only on $D_j,D_k$ and have bounded expectations by \eqref{eq:common-oscillation}. Summing over the arcs gives a bound polynomial in $T^{-1}$. We apply these bounds to the expectations in~\eqref{eq:weighted-selector} and~\eqref{eq:j-positive-definition}.

\subsection{A boundary one-arm estimate for the discarded arcs}
\label{subsec:shadow-count}

Recall from Section~\ref{subsec:shadow-allocation} that $S_C$ is the entire arc of the original boundary cut off by a white cluster $C$. On $H_3$, the maximal arcs $S_C$ have disjoint interiors. Their normalized quantum boundary lengths $\ell_C$ give, by Lemma~\ref{geo:allocation},
\[
 \overline G_\tau(x)
 =\mathbf1_{H_3}\sum_{C\in\mathcal M}
          \min\{x,1-x,\ell_C,1-\ell_C\}.
\]

The following lemma bounds the expected number of maximal discarded boundary arcs at each scale, restricted to the event $H_3$.
\begin{lemma}\label{lem:joint-shadow-count}
For an integer $n\ge0$, let $N_n$ count maximal arcs $S_C$ for which the shorter of $S_C$ and its complement has length comparable to $e^{-n}$ in the cylinder boundary coordinate. There are constants $C,c>0$, independent of $n$ and $\tau$, such that for all sufficiently large $\tau$,
\begin{equation}\label{eq:joint-shadow-count}
 \mathbb E[\mathbf1_{H_3}N_n]
       \le C e^{2n/3}p_3(\tau-c).
\end{equation}
For $\tau\le1$, the corresponding bound for arcs of length $\tau e^{-n}$ is $Ce^{2n/3}$ per boundary arc of length $\tau$.
\end{lemma}

\begin{proof}
Put $\delta=e^{-n}$ and cover the outer circle by $O(\delta^{-1})$ arcs $Q$ of length $\delta/100$. Remove the forced black boundary row and consider a free half-annulus around $Q$, from scale comparable to $\delta$ to a fixed radius $R$ in the outer boundary strip. If it has no black radial crossing, duality gives a white crosscut enclosing $Q$. The white cluster in $D\setminus K$ containing this crosscut has a discarded boundary arc containing $Q$ in its interior. Thus $Q$ cannot contain an endpoint of a maximal discarded arc whose two complementary lengths are at least a constant times $\delta$: overlapping discarded arcs are nested, and on $H_3$ two of them cannot cover the circle.

Each such endpoint therefore requires a black crossing of the half-annulus, with probability at most $C_R\delta^{1/3}$ by Cardy's formula~\cite{Smirnov}. Disjointness and the lower bound on arc lengths give a bounded number of initial endpoints, and likewise terminal endpoints, in each $Q$. On $H_3$, three black arms cross the inner annulus left after removing the outer strip. This event has probability $p_3(\tau-c)$ and is independent of the half-annulus crossing. Summing over $Q$ gives \eqref{eq:joint-shadow-count}. For $\tau\le1$, take $R$ comparable to $\tau$ and rescale each boundary arc to unit size.
\end{proof}

\subsection{Proofs of the three estimates}\label{subsec:proof-estimates}

\begin{proof}[Proof of Proposition~\ref{est:endpoint}]
For small modulus, use $T=2\pi\tau$ and the decomposition of the boundary fields \eqref{eq:thin-decomposition}. Consider an arc of length $Te^{-n}$ in $J_j$ and an arc $J_k$ in its complement, separated from $J_j$. Apply \eqref{eq:homogeneous-potential} conditionally on $C_T,E_i$. The sum over $n$, including its factor $e^{2n/3}$, is at most
\[
 C_\sigma c_T^\sigma\sqrt{x}
 e^{C_\sigma(D_j+D_k)}
 \big(a_j^{1/2}a_k^{\sigma-1/2}
       +a_j^{\sigma-1/2}a_k^{1/2}\big).
\]
Apply \eqref{eq:denominator-cancellation} separately to the two terms. The factors $a_j,a_k$ and $c_T^\sigma$ cancel, and the remaining inner negative moments are bounded. For subarcs in the same or neighboring partition arcs, use their disjointness and subdivide a bounded enlargement of the arc as above. They share a common multiplicative factor, and the convergent sum over common scales keeps both positive powers below $3/2$. For maximal arcs with shorter side at least $cT$, use \eqref{eq:tent-square-root} (there are at most $O(T^{-1})$ such arcs). For smaller arcs, Lemma~\ref{lem:joint-shadow-count} gives the factor $Ce^{2n/3}$ per partition arc. Thus, locally uniformly for $3/2<\sigma<2$,
\begin{equation}\label{eq:thin-endpoint}
 W_\sigma(\tau,x)
       \le C_\sigma\tau^{-K_\sigma}\sqrt{x(1-x)},
       \qquad 0<\tau\le\tau_0,
\end{equation}
for some finite $K_\sigma$ and a fixed $\tau_0>0$. Uniform covariance comparison on compact modulus intervals extends this bound to any fixed upper limit $\tau_0$. Factor out $\sqrt{x}$ before applying H\"older's inequality to the oscillations and negative moments.

For long cylinders, the joint covariance of the boundary traces differs from that of two independent Gaussian fields on the circle by a uniformly bounded smooth function. Condition on percolation to fix the arcs $S_C$. The covariance comparison in Lemma~\ref{lem:mixed-chaos} bounds the two outer masses together with the inner negative moment. Taking the percolation expectation and combining \eqref{eq:joint-shadow-count} with the scale sum \eqref{eq:chaos-square-root} gives the bound for $x\le1/2$. Symmetry in $x$ and $1-x$ then gives \eqref{eq:long-endpoint}.

By \eqref{eq:zero-mode-endpoint} and \eqref{weld:BAfield}, $\varphi_{-\sigma}(x)$ is the integral of $W_\sigma(\tau,x)$ against $(\sqrt2/\gamma)\eta(2i\tau)e^{2\pi\sigma^2\tau/3}\,d\tau$. As $\tau\downarrow0$, the modular identity gives
\begin{equation}\label{eq:eta-thin-endpoint}
 \eta(2i\tau)=(2\tau)^{-1/2}\eta(i/(2\tau))
       \asymp\tau^{-1/2}e^{-\pi/(24\tau)}.
\end{equation}
This exponential decay makes the small-modulus integral finite. As $\tau\to\infty$, the integral converges for $X(\sigma)<\xi_3$. Dividing by $\sqrt{x(1-x)}$ therefore gives a bounded function $u_\sigma$.
\end{proof}

\begin{proof}[Proof of Proposition~\ref{est:inhomogeneous}]
Choose $\max J-3/2<\theta<1/2$. For outer boundary masses $A,B>0$, put $\ell=A/(A+B)$. Then
\begin{equation}\label{eq:fractional-interval-weight}
 (A+B)^\sigma[\ell(1-\ell)]^\theta
 \le C_J\big(A^{\sigma-\theta}B^\theta
                         +A^\theta B^{\sigma-\theta}\big).
\end{equation}
The powers in \eqref{eq:fractional-interval-weight} satisfy Lemma~\ref{lem:mixed-chaos}. Including the inner negative moment bounds the field expectation uniformly over boundary arcs on long cylinders. Let $N_{\partial}$ be the maximal number of disjoint white crossings of a fixed-width strip next to the outer boundary; it bounds the number of white clusters meeting both boundary circles. Disjoint occurrence and Russo--Seymour--Welsh (RSW) crossing estimates give $N_{\partial}$ moments of every finite order, uniformly for large modulus. On $H_3$, a white cluster meeting both boundary circles gives four arms of colors $BBBW$. Integrating the field first therefore bounds the weighted count by $C_J\mathbb E[N_{\partial}\mathbf1_{A_4(\tau-c)}]$, where $A_4$ is the polychromatic four-arm event with color sequence $BBBW$. For $q>1$, H\"older's inequality gives
\[
 \mathbb E[N_{\partial}\mathbf1_{A_4}]
 \le \mathbb E[N_{\partial}^q]^{1/q}
                         \mathbb P(A_4)^{1-1/q}.
\]
The polychromatic four-arm exponent is $5/4$~\cite[Theorem~2]{BN}. Given $\epsilon>0$, choosing $q$ sufficiently large gives
\begin{equation}\label{eq:four-arm-bound}
 \mathbb E\!\left[(Z_o/Z_i)^\sigma\mathbf1_{H_3}
       \sum_{C\in\mathcal M_i}[\ell_C(1-\ell_C)]^\theta\right]
 \le C_{J,\epsilon}e^{-2\pi(5/4-\epsilon)\tau}.
\end{equation}
At small modulus, the argument leading to \eqref{eq:denominator-cancellation}, now with powers $\sigma-\theta,\theta$, gives a polynomial bound. The expected number of white clusters meeting both boundary circles is also bounded by a power of $\tau^{-1}$.

Multiplying \eqref{eq:four-arm-bound} by $\eta(2i\tau)e^{2\pi\sigma^2\tau/3}$, as in~\eqref{eq:j-positive-definition}, leaves, at large modulus, the exponential
\[
 \exp\left\{\left[\frac{2\pi}{3}(\sigma^2-4)
                                  +2\pi\epsilon\right]\tau\right\}.
\]
For sufficiently small $\epsilon$, this decays uniformly on $J$, and \eqref{eq:eta-thin-endpoint} controls small modulus. Since $\theta<1/2$, choose $0<a<1/2-\theta$. The bound
\[
 e^{a|\log(\ell/(1-\ell))|}\sqrt{\ell(1-\ell)}
                       \le C_a[\ell(1-\ell)]^\theta
\]
then proves the asserted weighted norm bound. To prove holomorphy, choose a compact interval $J'\subset(3/2,2)$ containing $J$ in its interior and choose $\theta,a$ above using $J'$. Fix $\sigma_-<\min J\le\max J<\sigma_+$ in $J'$. In the construction of $j_\sigma$, replace the Mellin weight by
\[
 \eta(2i\tau)e^{2\pi\sigma_+^2\tau/3}
 \big((Z_o/Z_i)^{\sigma_-}+(Z_o/Z_i)^{\sigma_+}\big).
\]
The resulting positive measure has finite $e^{a|t|}$-weighted mass by \eqref{eq:four-arm-bound} and the small-modulus bound. For a complex parameter $\sigma'$ in a smaller neighborhood of $J$,
\[
 |e^{2\pi(\sigma')^2\tau/3}|
       \le e^{2\pi(\operatorname{Re}\sigma')^2\tau/3}.
\]
Each fixed-order derivative introduces powers of $\tau$ and $\log(Z_o/Z_i)$. The strict bounds $(\operatorname{Re}\sigma')^2<\sigma_+^2$ and $\sigma_-<\operatorname{Re}\sigma'<\sigma_+$ absorb these factors into the displayed majorant, using $|\log u|^k\le C_{k,\delta}(u^\delta+u^{-\delta})$. At small modulus, \eqref{eq:eta-thin-endpoint} dominates the polynomial factors. For each fixed derivative order, a constant multiple of this measure dominates the derivatives and difference quotients. Dominated differentiation in weighted total variation proves the asserted holomorphy.

Finally, RSW gives positive probability to a separated $BBBW$ configuration in a fixed-modulus annulus. Its white cluster meets both boundary circles, and the discarded arc and its complementary arc both have positive quantum length almost surely. Hence $j_\sigma$ is nonzero, and convolution with the positive kernel $g$ gives the last assertion.
\end{proof}

\begin{proof}[Proof of Proposition~\ref{est:comparison}]\label{subsec:marked-comparison}
It remains to prove the lower bound. As $\tau\to\infty$, the boundary length measures $M_o,M_i$ associated with $h^{\mathrm{lat}}$ converge jointly to independent GMC measures on the circle: couple the Fourier modes of the Gaussian boundary fields with those of two independent limiting fields so that the smooth differences tend uniformly to zero in probability. The limiting measures are finite, nonzero, nonatomic and have full support. There are $\epsilon>0$, $M>1$ and $c_0>0$ such that
\[
 \mathcal R_\tau=\left\{M^{-1}<Z_o,Z_i<M,\quad
       \sup_{|I|\le2\epsilon}M_o(I)<Z_o/3\right\}
\]
has probability at least $c_0$ for all sufficiently large $\tau$. Choose $M,\epsilon$ so that these inequalities hold strictly with positive probability for the limiting measures. Nonatomicity makes the maximum mass of an arc of fixed length continuous under weak convergence, giving the same lower bound for large $\tau$.

On $\mathcal R_\tau$, two marks separated by half the quantum perimeter divide the outer boundary into two arcs, each of coordinate length at least $2\epsilon$. Conditional on the field, arm separation and RSW~\cite[Section~2.1]{BN} in a fixed-width outer strip give three black arms whose outer endpoints $u_1,u_2,u_3$ have counterclockwise order $p,u_1,u_2,q,u_3$, with probability at least $c_\epsilon p_3(\tau-c)$. Their landing intervals are separated by a distance bounded below in terms of $\epsilon$, uniformly in the first mark. The black boundary joins $p$ to $u_1,u_3$ for the condition at $p$, and $q$ to $u_2,u_3$ for the condition at $q$. For the joint condition, join $p$ to $u_1$ and $q$ to $u_2,u_3$. In each case the paths are disjoint except at their marked starting points and avoid the other mark, so $G=1$ by~\eqref{eq:geo:fan}. Averaging over the first mark preserves the bound. Since $(Z_o/Z_i)^\sigma\ge M^{-2\max J}$ on $\mathcal R_\tau$, integration over this event proves the lower bound.

By~\cite[Theorems~4 and~5]{BN}, $p_3(\tau)=\exp\{-2\pi\xi_3\tau+o(\tau)\}$ and $2/3<\xi_3<5/4$. For $3/2<\sigma<2$, equations \eqref{eq:zero-mode-endpoint} and \eqref{eq:marked-unmarked-comparison} imply convergence of $F$ for $X(\sigma)<\xi_3$ and divergence for $X(\sigma)>\xi_3$. The small-modulus contribution is finite by \eqref{eq:eta-thin-endpoint}. This proves \eqref{eq:physical-abscissa}. 
\end{proof}

\section{The spectral threshold and the main theorem}\label{sec:spectral}\label{sec:identification}

We identify the convergence endpoint of the Mellin transform $F$ with the spectral threshold $s_*$ of Theorem~\ref{thm:critical-spectral}.

\subsection{The Mellin transform and the spectral threshold}\label{subsec:proof-inputs}

Recall the Mellin transform $F$ in~\eqref{eq:physical-F}, and the function $u_\sigma$ and measure $j_\sigma$ in Proposition~\ref{weld:transform}. They satisfy
\begin{equation}\label{eq:bridge-equation}
 F(\sigma)=\tfrac12u_\sigma(0),\qquad (L_0-T_{-\sigma})u_\sigma=j_\sigma.
\end{equation}
Propositions~\ref{est:endpoint} and~\ref{est:comparison} show that $u_\sigma$ is bounded, nonnegative and even for $3/2<\sigma<\sigma_c$, where the right endpoint of the interval of convergence of $F$ is
\begin{equation}\label{eq:bridge-abscissa}
 \sigma_c=\sqrt{1/4+3\xi_3}\in(3/2,2).
\end{equation}
By Proposition~\ref{est:inhomogeneous}, $j_\sigma$ is a nonzero positive even measure and extends holomorphically near every compact subinterval of $(3/2,2)$ in the weighted total variation norm $\int e^{\eta_0|t|}|j|(dt)$ for some $\eta_0>0$. In particular, $j_\sigma$ is holomorphic near $\sigma_c$.

Recall $L_0=-d^2/dt^2+1/4$ and $T_s$ from \eqref{eq:intro-kernel}--\eqref{eq:intro-T}. For $|\operatorname{Re}s|<2$, put
\[
 g(t)=e^{-|t|/2},\qquad P_sf=g*(T_sf).
\]
Here $g$ is the Green kernel of $L_0$ on $\mathbb R$, and $*$ denotes convolution. For $R>0$, let $g_R$ be the Dirichlet Green kernel of $L_0$ on $I_R=(-R,R)$, and let $K_{R,s}=g_RT_s$ act on $C([-R,R])$, with both integrations restricted to $I_R$. This finite-interval operator is defined for all $s\in\mathbb C$. Write $r_R(s)$ for its spectral radius, and for real $s$, this is its Perron eigenvalue. Write $T_s(t,u)$ and $P_s(t,u)$ for the kernels of $T_s$ and $P_s$, respectively.

We use the spaces
\[
 \mathcal Y_\eta=\{f\in C_{\rm even}(\mathbb R):\|f\|_\eta:=\sup_t e^{\eta|t|}|f(t)|<\infty\},\qquad
 \mathcal X=\{u\in C_{\rm even}(\mathbb R):\lim_{|t|\to\infty}u(t)\in\mathbb C\},
\]
where $\eta>0$, $C_{\rm even}(\mathbb R)$ denotes continuous even complex-valued functions, and $\mathcal X$ has the supremum norm.

\begin{theorem}[Spectral threshold]\label{thm:critical-spectral}
\label{lem:physical-sign-bounded-resolvent}
There is $s_*\in(1,2)$ such that, for each $\epsilon\in\{-1,1\}$, the equation $P_{\epsilon s_*}u=u$ has a strictly positive bounded even solution $\psi_\epsilon$. The solution $\psi_\epsilon$ has a positive common limit at both ends and is unique up to multiplication by a positive constant. The space of bounded solutions of $P_{\epsilon s_*}u=u$ is one-dimensional. For each sign $\epsilon$, no other $s\in(0,2)$ admits a positive bounded solution of $P_{\epsilon s}u=u$.

For $0\leq s<s_*$, the minimal Green operator $\mathcal G_{\epsilon s}f=\sum_{n\geq0}P_{\epsilon s}^nf$ of $I-P_{\epsilon s}$ is bounded from $\mathcal Y_\eta$ to $\mathcal X$ and is analytic in a complex neighborhood of every compact subinterval of $[0,s_*)$, with $\eta>0$ sufficiently small.  The equation $(I-P_{\epsilon s})u=f$ has at most one bounded solution. The spectral radii $r_R(s)$ satisfy $r_R(s)=r_R(-s)$ and
\begin{equation}\label{eq:positive-threshold-definition}
 s_*=\inf\{s>0:\sup_{R>0}r_R(s)>1\},
 \qquad \sup_{R>0}r_R(s)>1\quad(s_*<s<2).
\end{equation}
\end{theorem}
We postpone the proof to Section~\ref{subsec:two-sign-threshold}. Below, we first finish the proof of Theorem~\ref{thm:main}.

\subsection{Proof of Theorem~\ref{thm:main}}
\label{subsec:proof-identification}

We need the following basic fact on the Laplace transform.
\begin{lemma}\label{lem:positive-Laplace-boundary}
Let $\mu$ be a positive measure on $[0,\infty)$ and suppose $\int e^{sv}\mu(dv)$ has a finite right endpoint $\beta$ of its interval of convergence. The resulting holomorphic function has no analytic continuation to a neighborhood of $\beta$.
\end{lemma}

\begin{proof}
If it were analytic on a disk of radius $r$ about $\beta$, choose $0<\delta<r/3$ and expand at $a=\beta-\delta$. The Taylor series converges at $a+3\delta/2>\beta$, whereas its coefficients are $\int v^n e^{av}\mu(dv)/n!\ge0$. Tonelli's theorem would then make the original integral finite at $a+3\delta/2$, a contradiction.
\end{proof}

\begin{proof}[Proof of Theorem~\ref{thm:main}]
For $3/2<\sigma<\sigma_c$, boundedness of $u_\sigma$ allows us to convolve \eqref{eq:bridge-equation} with $g$:
\begin{equation}\label{eq:id-short-integral-equation}
 u_\sigma=P_{-\sigma}u_\sigma+b_\sigma,
 \qquad b_\sigma=g*j_\sigma>0.
\end{equation}
The difference $w$ between the two sides is bounded and satisfies $L_0w=0$, so the coefficients of $e^{t/2}$ and $e^{-t/2}$ both vanish.

\emph{First, $\sigma_c\le s_*$.} Let $\ell_R$ be a positive left eigenfunctional of $K_{R,-\sigma}$ corresponding to $r_R(-\sigma)$.  Since $g_R\le g$ and the omitted integrals are positive, \eqref{eq:id-short-integral-equation} gives
\[
 (1-r_R(-\sigma))\ell_R(u_\sigma)
       \ge\ell_R(b_\sigma)>0.
\]
All these pairings are on $I_R$.  Hence $r_R(-\sigma)<1$ for every $R$, and transposition gives the same inequality for $r_R(\sigma)$. If $\sigma_c>s_*$, choose $\max\{3/2,s_*\}<\sigma<\sigma_c$. Theorem~\ref{thm:critical-spectral} gives an $R$ with $r_R(\sigma)>1$, a contradiction. Thus $\sigma_c\le s_*$, and in particular $s_*>3/2$.

\emph{Next, $\sigma_c\ge s_*$.} Suppose otherwise, and take a compact interval $J\Subset(3/2,s_*)$ with $\sigma_c$ in its interior. The estimate on $j_\sigma$ and
\[
 \sup_t e^{\eta|t|}|(g*j)(t)|
       \le\int e^{\eta|u|}|j|(du),\qquad 0<\eta<1/2,
\]
make the continuation $b_\zeta=g*j_\zeta$ holomorphic for complex $\zeta$ near $J$ in $\mathcal Y_\eta$. Choose $\eta$ also below the decay rate of $j_\sigma$ and the weight allowed by Theorem~\ref{thm:critical-spectral}.  Uniqueness among bounded solutions in that theorem and \eqref{eq:id-short-integral-equation} give $u_\zeta=\mathcal G_{-\zeta}b_\zeta$ for real $\zeta<\sigma_c$ close to $\sigma_c$. Consequently $\tfrac12(\mathcal G_{-\zeta}b_\zeta)(0)$ analytically continues $F$ through $\sigma_c$.

To exclude this continuation, split the defining integral for $F$ at $m=1$.  By \eqref{weld:densitybound}, the part over $m\ge1$ is analytic for $\operatorname{Re}\zeta>-1/2$.  In the other part, set $v=-\log m$:
\[
 F_<(\zeta):=\int_0^1 m^{-\zeta}z(1,1/2;m)\,dm=\int_0^\infty e^{\zeta v}\mu(dv),\qquad
 \mu(dv)=e^{-v}z(1,1/2;e^{-v})\,dv.
\]
The Laplace transform of this positive measure has $\sigma_c$ as the right endpoint of its interval of convergence.  Lemma~\ref{lem:positive-Laplace-boundary} gives the contradiction.

Thus $\sigma_c=s_*$, and \eqref{eq:bridge-abscissa} yields $\xi_3=(s_*^2-1/4)/3$.  The positive bounded solution, its uniqueness and its limit at infinity follow from Theorem~\ref{thm:critical-spectral}; the simple zero of $\mathscr D$ at $s_*$ follows from Proposition~\ref{prop:simple-resonance}.\qedhere
\end{proof}

\subsection{The limiting convolution kernel and bounded solutions}\label{subsec:free-bounded}

Recall $x(t)=e^t/(1+e^t)$ and $a(t)=1-x(t)$ from \eqref{eq:intro-kernel}. For fixed $t,u$, note that $P_s(t+v,u+v)\to q_0(t-u)$ as $v\to\pm\infty$, where
\begin{equation}\label{eq:critical-kernels}
 q_0(t)=\frac{3\sqrt3}{8\pi(\cosh(t/2)+1/2)},\qquad
 L_0q_0=H,\qquad g*H=q_0.
\end{equation}
The weight in $T_s$ is
\begin{equation}\label{eq:W-s}
 W_s(t,u)=x(t)^{1+s}x(u)^{1-s}+a(t)^{1+s}a(u)^{1-s}.
\end{equation}
Thus $T_s(t,u)=H(t-u)W_s(t,u)$.  Both $T_s$ and $P_s$ preserve even functions.

To prove Theorem~\ref{thm:critical-spectral}, we first study the convolution operator with kernel $q_0(t-u)$.
Write $Qf=q_0*f$. With Fourier convention $\widehat f(k)=\int e^{ikt}f(t)\,dt$, direct integration gives
\begin{equation}\label{eq:critical-symbol}
 \widehat q_0(k)=\frac3{1+2\cosh(4\pi k/3)},\qquad
 \int q_0=1,\qquad
 \sigma_0^2:=\int t^2q_0(t)\,dt=\frac{32\pi^2}{27}.
\end{equation}
The walk with increment density $q_0$ is recurrent. The following Lyapunov argument also applies to the Doob transform below. For $V(t)=\log(1+t^2)$, Taylor expansion and exponential jump tails give $QV(t)-V(t)=-\sigma_0^2/t^2+o(t^{-2})$. Stop on entering a sufficiently large compact interval or exiting $[-R,R]$. Optional stopping bounds the exit probability by $V(t)/\log(1+R^2)$, and strict positivity of the kernel excludes indefinite avoidance within a bounded region. Letting $R\to\infty$ proves recurrence. Stopping on shrinking neighborhoods then shows that every continuous finite nonnegative $Q$-superharmonic function ($Qf\le f$), and every bounded $Q$-harmonic function, is constant.

The kernel of the perturbation $R_s=P_s-Q$ and all its coordinate derivatives decay exponentially in $|t|+|u|$, uniformly on compact subsets of $|\operatorname{Re}s|<2$. To see this, before convolution with $g$ its kernel is $H(v-u)(W_s(v,u)-1)$. For $|\operatorname{Re}s|\leq r$ with $1\leq r<2$ and $v\geq u\geq0$, its terms are bounded by constants times $e^{-(v-u)-u}$ and $e^{-(2-r)(v-u)-2u}$. For $v\geq0\geq u$, use $e^{-v-(2-r)|u|}$ and $e^{-(2-r)v-|u|}$.  Reflection and interchange of coordinates cover the remaining regions.  Parameter derivatives add polynomial factors, which can be absorbed by a smaller exponential rate.

The potential kernel of $Q$ is explicit:
\begin{equation}\label{eq:potential-A}
 A(t)=\frac{27}{32\pi^2}t\coth\frac{3t}{4},\qquad
 (I-Q)A=-q_0,\qquad
 A(t)=\frac{|t|}{\sigma_0^2}+O((1+|t|)e^{-3|t|/2}).
\end{equation}
Define $A(0)$ by continuity.  The identities follow from \eqref{eq:critical-symbol}.  Define an operator on even exponentially decaying functions by
\begin{equation}\label{eq:S-potential}
 (Sf)(t)=f(t)-\int_{\mathbb R}
 \left[\frac{A(t-u)+A(t+u)}2-A(t)\right]f(u)\,du.
\end{equation}
The second difference cancels the linear tails, so $S:\mathcal Y_\eta\to\mathcal Y_\eta$ is bounded for sufficiently small $\eta>0$.  Moreover,
\begin{equation}\label{eq:S-identity}
 (I-Q)Sf=f-q_0\int_{\mathbb R}f.
\end{equation}

Every bounded even solution of $(I-P_s)u=f$, with $f\in\mathcal Y_\eta$ and $|\operatorname{Re}s|<2$, has a common limit at both ends. To see this, set $w=f+R_su$, so $(I-Q)u=w$ and $w$ decays exponentially. Choose a smooth compactly supported function $\chi$ equal to one near zero, and write $\|\cdot\|_1$ for the $L^1(\mathbb R)$ norm. Testing against $\chi(t/R)$ gives
\[
 \left|\int w(t)\chi(t/R)\,dt\right|
 \leq\|u\|_\infty\|(I-Q)\chi(\cdot/R)\|_1
 \leq\frac{\sigma_0^2}{2}\|u\|_\infty
                  \|\partial_t^2\chi(\cdot/R)\|_1=O(R^{-1}).
\]
Thus $\int w=0$, and \eqref{eq:S-identity} makes $u-Sw$ bounded and $Q$-harmonic.  Consequently
\begin{equation}\label{eq:bounded-tail-localization}
 u=c+S(f+R_su),
\end{equation}
where $c\in\mathbb C$ and the second term decays exponentially. Thus bounded even solutions belong to $\mathcal X$.  In the original coordinate, $\varphi(x(t))=\sqrt{x(t)a(t)}u(t)$ has the endpoint behavior corresponding to boundedness of $u$ in Proposition~\ref{est:endpoint}, with no logarithmic term.

\subsection{Proof of Theorem~\ref{thm:critical-spectral}}\label{subsec:two-sign-threshold}

\begin{proof}[Proof of Theorem~\ref{thm:critical-spectral}]
\emph{Comparison on finite intervals.} The Dirichlet Green kernel of $L_0$ on $I_R=(-R,R)$ is
\[
 g_R(t,v)=\frac{2\sinh((\min(t,v)+R)/2)
                  \sinh((R-\max(t,v))/2)}{\sinh R}.
\]
The operator $K_{R,s}$ is compact. For real $s$, its kernel is strictly positive in the interior and $r_R(s)$ increases with $R$, and write $r_\infty(s)=\sup_Rr_R(s)$. Kernel transposition, denoted by $\mathsf T$, gives $T_s^{\mathsf T}=T_{-s}$. Combining this with the equality of the nonzero spectra of $AB$ and $BA$ gives $r_R(s)=r_R(-s)$. H\"older's inequality applied to the iterated integrals for positive kernels gives log convexity in real $s$, also for $r_\infty$.

After translating $I_R$ to either end, $K_{R,s}$ converges to the operator obtained by replacing $T_s(t,u)$ with $H(t-u)$. The spectral radii of these limiting operators tend to one as $R\to\infty$.  Conjugating by the square root of the Dirichlet operator $L_0$ gives, for a smooth compactly supported test function $\chi_R(t)=\chi(t/R)$, the Rayleigh quotient
\[
 \frac{\iint H(t-u)\chi_R(t)\chi_R(u)\,dt\,du}
 {\int(|\chi_R'|^2+\tfrac14\chi_R^2)\,dt}\longrightarrow1
\]
by $\int H=1/4$ and the domination by $Q$. It follows that $r_\infty(s)\geq1$ for every real $s$. For $|s|\leq1$, H\"older gives $W_s\leq1$, so $P_s\leq Q$ and $r_\infty(s)=1$. Whenever $r_\infty(s)=1$, normalize the Perron eigenfunctions on finite intervals at zero.  The bounds for $g_R$ give local comparison and equicontinuity, with comparison constants tending to $e^{|t-u|/2}$. A locally uniform subsequential limit $v_s$ is continuous, even, and satisfies
\begin{equation}\label{eq:spectral-supersolution}
 v_s(0)=1,\qquad P_sv_s\leq v_s,\qquad v_s(t)\geq e^{-|t|/2}.
\end{equation}
Fatou's lemma gives the inequality.

\emph{A return estimate.} We need a lower bound at infinity for positive supersolutions, uniform in the sign and in $1\le s<2$.  Write $M_f$ for multiplication by $f$, and juxtaposition with $g$ denotes convolution. For $1\leq s<2$, retain the $a$ term in $W_s$ to obtain
\[
 P_s\geq\mathsf K_+:=gM_{a^3}(H*),\qquad
 P_{-s}\geq\mathsf K_-:=QM_{a^3},\qquad
 \mathsf K_\pm\leq Q.
\]
Both kernels are obtained from the $q_0$ walk by killing rejected transitions.  Put $h_{\rm ret}(t)=1+\exp(\min(t,0)/4)$, with value zero at the cemetery state.  The moment $\int e^{v/4}q_0(v)\,dv=3/2$ and exponential tails give
\[
 Qh_{\rm ret}(t)-h_{\rm ret}(t)
       \ge\tfrac12e^{t/4}-Ce^{t/2},\qquad t\le0.
\]
Since $1-a(u)^3\le3\min(e^u,1)$ and $\int e^{-|t-u|/2}\min(e^u,1)\,du\le4e^{t/2}$, the differences $(Q-\mathsf K_\epsilon)h_{\rm ret}$ are at most $Ce^{t/2}$. Consequently, for one fixed $R_0>0$ and both signs,
\begin{equation}\label{eq:return-barrier}
 \mathsf K_\epsilon h_{\rm ret}(t)-h_{\rm ret}(t)
       \ge\tfrac12e^{t/4}-Ce^{t/2}>0,\qquad t\le-R_0.
\end{equation}
For the chain $(X_n)_{n\ge0}$ with transition kernel $\mathsf K_\epsilon$, let $\mathbf P_t^{\mathsf K_\epsilon}$ denote its law started at $t$. Let $\tau$ be the first entry into $[-R_0,\infty)$ and $\zeta$ the killing time.  Recurrence of the coupled walk with increment density $q_0$ makes $\tau\wedge\zeta$ finite.  Bounded optional stopping gives $\mathbf P_t^{\mathsf K_\epsilon}(\tau<\zeta)\ge1/2$ for $t<-R_0$. For $J>0$, the exponential tail ratio
\[
 \sup_{b\ge0}\frac{\int_{b+J}^\infty q_0(v)\,dv}
                         {\int_b^\infty q_0(v)\,dv}
       \le Ce^{-J/2}
\]
bounds the overshoot of the walk with increment density $q_0$ at its first crossing. Coupling with this walk up to the killing time, choose a fixed $J$ large enough to obtain
\begin{equation}\label{eq:return-window}
 \mathbf P_t^{\mathsf K_\epsilon}
 (X_\tau\in[-R_0,-R_0+J],\ \tau<\zeta)\ge\tfrac14,
                         \qquad t<-R_0.
\end{equation}
Optional stopping now shows that every positive continuous even $P_{\epsilon s}$-superharmonic function is bounded below at both ends by a positive multiple of its minimum on the fixed interval $[-R_0,-R_0+J]$.

We use this estimate to show that $r_R(s)>1$ for some $R>0$ and $s<2$. If $r_\infty(s)=1$ with $1\leq s<2$, apply \eqref{eq:return-window} to the supersolution in~\eqref{eq:spectral-supersolution}, and then $v_s(t)\ge c_0>0$ for $t<-R_0$, uniformly in $s$.  Restricting the intermediate integration variable to $v\in[0,1]$ and retaining the $x$ term in $W_s$ gives $P_s(0,u)\ge c_1e^{(2-s)u}$ for $u<-R_0$, with $c_1>0$ independent of $s$.  Thus $1\ge c_0c_1e^{-R_0(2-s)}/(2-s)$, which fails as $s\uparrow2$. Fix $\bar s\in(1,2)$ with $r_\infty(\bar s)>1$. The identity $r_R(s)=r_R(-s)$ gives the same conclusion at $-\bar s$.

\emph{Bounded Green potentials below the threshold.} Choose $\eta>0$ below a uniform exponential decay rate of $R_s$ for $s\in[-\bar s,\bar s]$, and put $w_\eta(t)=e^{-\eta|t|}$. The potentials $\sum_{n\ge0}P_s^nw_\eta$ are bounded for $s=0,\pm1$. Indeed, all three kernels are dominated by $Q$ and have a strictly positive one-step killing probability on every compact interval. Let $v$ be a smooth function equal to $e^{-\eta|t|}$ outside a compact interval, with $0<v\leq1$.  Then $Qv-v\geq c w_\eta$ outside a compact interval, since the exponential moment of $q_0$ at $\eta$ exceeds one. For each of the three kernels $P$, a sufficiently large $C$ makes
\[
 (I-P)(C-v)\geq(Qv-v)+(C-1)(1-P1)\geq c'w_\eta
\]
everywhere.  This supersolution bounds the minimal potential.  Recurrence of the walk with increment density $q_0$ and the positive killing probability on compact sets imply almost sure killing, hence uniqueness among bounded solutions.  By \eqref{eq:bounded-tail-localization}, the resulting Green operators $\mathcal G_0,\mathcal G_{\pm1}$ map $\mathcal Y_\eta$ boundedly into $\mathcal X$.

For $\epsilon=\pm1$, let
\[
 h_{\epsilon,s}=\sum_{n\geq0}P_{\epsilon s}^nw_\eta,
 \qquad b_\epsilon=\sup\{s\in[0,\bar s]:
                              \|h_{\epsilon,s}\|_\infty<\infty\}.
\]
H\"older's inequality applied to the iterated integrals for the positive kernels shows that the set of parameters for which $h_{\epsilon,s}$ is bounded is an interval.  Whenever $h_{\epsilon,s}$ is bounded, positivity and $|f|\leq\|f\|_\eta w_\eta$ give a bounded operator $\mathcal G_{\epsilon s}:\mathcal Y_\eta\to\mathcal X$. The interval is open: namely, if $h_{\epsilon,s}$ is bounded, then for $s'$ close to $s$ the positive Neumann series
\[
 [I-\mathcal G_{\epsilon s}|P_{\epsilon s'}-P_{\epsilon s}|]^{-1}
 h_{\epsilon,s}
\]
is a bounded supersolution for the potential at $s'$. Here the operator obtained by taking the absolute value of the kernel difference has norm $O(|s'-s|)$ from bounded even functions to $\mathcal Y_\eta$. This proves $b_\epsilon>1$. A bounded positive potential implies $r_R(\epsilon s)<1$ by pairing on $I_R$ with a positive left eigenfunctional of $K_{R,\epsilon s}$. Since $r_R(\epsilon\bar s)>1$ for some $R$ and $r_R$ is continuous, we have $b_\epsilon<\bar s$.

For $s<b_\epsilon$, the equation $P_{\epsilon s}u=u$ has no nonzero bounded solution. When $s\leq1$, the walk with transition kernel bounded above by $Q$ is killed almost surely. When $s\geq1$, \eqref{eq:return-window} bounds $h_{\epsilon,s}$ below away from zero, so it dominates the absolute value of such a solution.  Iterating $P_{\epsilon s}$ excludes it, because $P_{\epsilon s}^nh_{\epsilon,s} =\sum_{k\geq n}P_{\epsilon s}^kw_\eta\to0$ pointwise.

For $z$ in a complex neighborhood of $[-\bar s,\bar s]$, define the compact analytic family
\begin{equation}\label{eq:compact-Green-reduction}
 \mathcal A_z=\mathcal G_0(P_z-P_0):\mathcal X\to\mathcal X,
 \qquad \mathcal G_z=(I-\mathcal A_z)^{-1}\mathcal G_0.
\end{equation}
The second identity holds wherever $I-\mathcal A_z$ is invertible. Compactness follows from equicontinuity on compact sets and the kernel bounds with exponential decay rate strictly larger than $\eta$. Uniqueness among bounded solutions at $s=0$ makes $(I-P_z)u=f$ equivalent to $(I-\mathcal A_z)u=\mathcal G_0f$.  Below $b_\epsilon$ the homogeneous equation has no nonzero bounded solution, so the Fredholm alternative gives the inverse and its local complex analyticity.

\emph{A positive bounded solution at the threshold.} As $s\uparrow b_\epsilon$, $\|h_{\epsilon,s}\|_\infty$ diverges; otherwise Fatou's lemma and the preceding Neumann-series bound extend finiteness past $b_\epsilon$.  Normalize these potentials by their supremum norms.  Compactness in $h_{\epsilon,s}=\mathcal G_0w_\eta+ \mathcal A_{\epsilon s}h_{\epsilon,s}$ produces a nonzero nonnegative solution of $P_{\epsilon b_\epsilon}u=u$ in $\mathcal X$. Strict positivity of the kernel makes the limit positive.  The return estimate \eqref{eq:return-window} and \eqref{eq:bounded-tail-localization} give a positive common limit at both ends and an exponentially decaying remainder. Normalize that limit to one and denote the resulting function by $\psi_\epsilon$.

\emph{Recurrence and equality of the thresholds.} The Doob transform is the Markov kernel
\begin{equation}\label{eq:critical-Doob}
 \widetilde P_\epsilon(t,u)=P_{\epsilon b_\epsilon}(t,u)
                  \frac{\psi_\epsilon(u)}{\psi_\epsilon(t)}.
\end{equation}
Because $\psi_\epsilon-1$ and the kernel of $R_{\epsilon b_\epsilon}$ decay exponentially, the transition kernel $\widetilde P_\epsilon$ has exponentially small drift, variance tending to $\sigma_0^2$, and uniform exponential jump tails. The Lyapunov argument for $Q$ therefore gives recurrence of the chain with transition kernel $\widetilde P_\epsilon$.  Every continuous finite nonnegative superharmonic function, and every bounded harmonic function, is constant.  Dividing a bounded solution of $P_{\epsilon b_\epsilon}u=u$ by $\psi_\epsilon$ proves uniqueness up to a scalar factor.

The positive bounded solution gives $r_\infty(\epsilon b_\epsilon)=1$. If $b_\epsilon<s<2$ and $r_\infty(\epsilon s)=1$, take $v_{\epsilon s}$ from \eqref{eq:spectral-supersolution} and put $\theta=b_\epsilon/s$. H\"older, retaining the intermediate integration and both terms in $W_s$, gives
\[
 P_{\epsilon b_\epsilon}(v_{\epsilon s}^\theta)
 \leq(P_{\epsilon s}v_{\epsilon s})^\theta(P_01)^{1-\theta}<v_{\epsilon s}^\theta.
\]
Here $P_01(t)<1$ at every finite $t$. Division by $\psi_\epsilon$ contradicts recurrence. Thus $r_\infty(\epsilon s)>1$ for $s>b_\epsilon$, and evenness and log convexity of the spectral radii on finite intervals extend this conclusion to $s\geq2$. It follows that $b_\epsilon=\inf\{s>0:\sup_Rr_R(\epsilon s)>1\}$. The identity $r_R(s)=r_R(-s)$ gives $b_+=b_-=:s_*$ and \eqref{eq:positive-threshold-definition}. A positive bounded solution at any larger parameter would contradict $r_R>1$ on some finite interval; at a smaller parameter it is excluded by uniqueness among bounded solutions. Finally, a bounded nonnegative solution with a nonzero nonnegative exponentially decaying inhomogeneous term at $\epsilon s_*$ would give a nonconstant nonnegative superharmonic function for the same Doob transform.  Recurrence excludes such a solution, so $\mathcal G_s$ is singular at $s=\epsilon s_*$.
\end{proof}

\appendix
\section{Proof of Proposition~\ref{prop:intro-determinant}}
\label{sec:determinant}

To prove Proposition~\ref{prop:intro-determinant}, we separate the limit at infinity from the decaying part of a bounded solution of $P_s\psi=\psi$. This gives the Fredholm determinant $\mathscr D(s)$, whose zeros are exactly the parameters admitting a nonzero bounded even solution. Relative entropy shows that its first positive zero is simple.

\subsection{The bordered determinant}\label{subsec:bordered-determinant}

For $|\operatorname{Re}s|<2$, write a bounded even solution of $P_s\psi=\psi$ as $\psi=c+u$, where $c$ is its common limit at infinity and $u$ decays exponentially. Use the operators $R_s=P_s-Q$ and $S$ from Section~\ref{subsec:free-bounded}.  By \eqref{eq:bounded-tail-localization}, the equation $P_s\psi=\psi$ becomes
\[
 u=SR_s(c+u),\qquad \int_{\mathbb R}R_s(c+u)=0.
\]
The second identity is the solvability condition for $(I-Q)\psi=R_s\psi$.  Set
\begin{equation}\label{eq:bordered-entries}
 \begin{gathered}
 \mathcal K_s=SR_s,\qquad b_s=\mathcal K_s1,\qquad
 c_s(u)=\int_{\mathbb R}R_s(t,u)\,dt,\\
 d_s=\int_{\mathbb R}c_s(u)\,du,\qquad
 \ell_s(f)=\int_{\mathbb R}c_s(u)f(u)\,du.
 \end{gathered}
\end{equation}
Here $\ell_s$ is complex linear and $1$ denotes the constant function. On $\mathbb C\oplus L^2_{\mathrm{even}}(\mathbb R)$ define
\begin{equation}\label{eq:critical-determinant}
 \quad
 \mathscr D(s)=\det
 \begin{pmatrix}
 d_s&\ell_s\\-b_s&I-\mathcal K_s
 \end{pmatrix},\qquad |\operatorname{Re}s|<2.
 \quad
\end{equation}
This block operator is the identity plus a trace-class operator.  Indeed, the kernel of $SR_s$ is smooth and decays exponentially.  Factor $\mathcal K_s$ on both sides by $(1-\partial_t^2+t^2)^{-1}$: the outside factors are Hilbert--Schmidt and the middle factor is bounded.  The same estimates, uniform on compact parameter sets, prove holomorphy in trace norm.

All entries are given by explicit kernels and their integrals.  If $(c,u)$ is a null vector of the block operator in~\eqref{eq:critical-determinant}, the smoothing property of $\mathcal K_s$ implies that $u$ decays exponentially, and the two block equations give $(I-P_s)(c+u)=0$.  Conversely, every bounded even solution of $P_s\psi=\psi$ has this decomposition by \eqref{eq:bounded-tail-localization}.  Hence
\begin{equation}\label{eq:det-equivalence}
 \mathscr D(s)=0
 \quad\Longleftrightarrow\quad
 \text{$P_su=u$ has a nonzero bounded even solution.}
\end{equation}
The scalar row and column account for the constant solutions of $(I-Q)u=0$.  Theorem~\ref{thm:critical-spectral} shows that $s_*$ is the first positive zero of $\mathscr D$.

\subsection{Simplicity of the first positive zero}\label{subsec:resonance-simplicity}

Theorem~\ref{thm:critical-spectral} gives a one-dimensional space of bounded solutions of $P_{s_*}u=u$. We prove that $s_*$ is a simple zero of $\mathscr D$ by constructing a positive adjoint solution and showing that its pairing with $P'_{s_*}\psi$ is strictly positive. For real kernels, $B^*$ denotes the formal adjoint with respect to Lebesgue measure, obtained by transposing the kernel. The notation $B_s^\prime$ denotes differentiation with respect to $s$, and all displayed pairings are bilinear.

\begin{proposition}\label{prop:simple-resonance}
The determinant $\mathscr D$ has a simple zero at $s_*$.  There are strictly positive bounded functions $\psi,\phi$ with $P_{s_*}\psi=\psi$ and $P_{s_*}^*\phi=\phi$ such that, for $f\in\mathcal Y_\eta$ with $\eta>0$ sufficiently small,
\begin{equation}\label{eq:simple-Green-pole}
 \mathcal G_s f=
 \frac{\psi\int_{\mathbb R}\phi f}{M(s_*-s)}+O_{\mathcal X}(1),
 \qquad M=\int_{\mathbb R}\phi P'_{s_*}\psi>0.
\end{equation}
Here $O_{\mathcal X}(1)$ denotes a remainder that is analytic and bounded in $\mathcal X$ near $s_*$.
\end{proposition}

\begin{proof}
\emph{The adjoint solution and the block operator.} Normalize the solutions $\psi_\pm$ of Theorem~\ref{thm:critical-spectral} to have limit one at infinity. Set $\psi=\psi_+$ and
\[
 \phi=T_{-s_*}\psi_-.
\]
Since $g*\phi=\psi_-$, the identity $P_{s_*}^*=T_{-s_*}g*$ gives $P_{s_*}^*\phi=\phi$.  Positivity of $T_{-s_*}$ gives $\phi>0$. The kernel bounds for $T_{-s_*}-H*$ and the exponential decay of $\psi_--1$ further give
\[
 \phi(t)=\frac14+O(e^{-\eta|t|})
\]
for sufficiently small $\eta>0$.  Thus $\psi$ and $\phi$ are bounded with positive limits at infinity.

To compute the coefficient of $(s_*-s)^{-1}$, let $\mathcal B_s$ denote the block operator in \eqref{eq:critical-determinant} and construct its left null vector. Write $v=-R_{s_*}^*\phi$.  This function is even and decays exponentially.  The adjoint equation gives $(I-Q)\phi=R_{s_*}^*\phi$, and therefore $\int v=0$; here one integrates $(I-Q)(\phi-1/4)$, which is absolutely integrable.  From the explicit formula for $S^*$,
\[
 S^*v=v-A*v+\int_{\mathbb R}A(t)v(t)\,dt,
 \qquad (I-Q)S^*v=v.
\]
Evenness and $\int v=0$ make $S^*v$ bounded.  Thus $\phi+S^*v$ is a bounded $Q$-harmonic function, hence a constant $\beta$.  Consequently
\[
 v=\mathcal K_{s_*}^*v-\beta c_{s_*},\qquad
 \beta d_{s_*}-\int v b_{s_*}=0.
\]
The second identity follows by integrating the first.  Hence $(\beta,v)$ spans the left nullspace: the block operator $\mathcal B_{s_*}$ has Fredholm index zero and a one-dimensional right nullspace.  For $(I-P_s)u=f$ with $f\in\mathcal Y_\eta$, the right-hand side of the equation for $\mathcal B_s$ is $(-\int f,Sf)$, whose pairing with $(\beta,v)$ is $-\int\phi f$.

\emph{The derivative pairing.} For one application of $P_s$, retain the intermediate point $y$ and the choice $p\in\{x,a\}$. At $s=0$ the resulting subprobability measure on $\mathbb R^2\times\{x,a\}$ is
\[
 \mu_t(dy\,du,p)=g(t-y)H(y-u)p(y)p(u)\,dy\,du,
 \qquad m_0(t)=\mu_t(1)=P_0 1(t)<1.
\]
Here $\mu_t(1)$ denotes its total mass. Tilting $\mu_t$ by $e^{s_*[\log p(y)-\log p(u)]}$ and applying the Doob transform with $\psi$ gives the probability measure
\[
 \widetilde\mu_t(dy\,du,p)
  =e^{s_*[\log p(y)-\log p(u)]}
       \frac{\psi(u)}{\psi(t)}\mu_t(dy\,du,p).
\]
The Markov kernel $t\mapsto\widetilde\mu_t(u\in\cdot)$ has invariant measure $\pi(t)dt=\psi(t)\phi(t)dt$.  For a probability measure $\nu$ and a finite nonzero subprobability measure $\mu$ with $\nu\ll\mu$, write
\[
 D(\nu\Vert\mu)=\int\log\frac{d\nu}{d\mu}\,d\nu
 =D(\nu\Vert\mu/\mu(1))-\log\mu(1).
\]
Relative entropy and invariance give
\begin{equation}\label{eq:entropy-transversality}
 s_*\int\phi P'_{s_*}\psi
 =\int_{\mathbb R}\pi(t)
          D(\widetilde\mu_t\Vert\mu_t)\,dt
 \ge\int_{\mathbb R}\pi(t)[-\log m_0(t)]\,dt>0.
\end{equation}
The terms involving $\log\psi$ cancel by invariance and absolute integrability. For the two terms indexed by $p$, the logarithmic increment satisfies
\[
 \sum_{p=x,a}\iint H(y-u)p(y)^{1+s}p(u)^{1-s}
           |\log p(y)-\log p(u)|\,dy\,du<\infty,
 \qquad 0<s<2.
\]
In $u\ge y\ge0$ the $a$ term is bounded by $C e^{-2y-(2-s)(u-y)}(u-y)$, and the bound for the $x$ term has an additional factor $e^{-y}\min(1,u-y)$.  The other regions are controlled by the same quadrant bounds. Since $\log\psi$ also decays exponentially, the integral of $\log\psi(u)-\log\psi(t)$ is absolutely convergent and vanishes by invariance.  This proves both finiteness and strict positivity in \eqref{eq:entropy-transversality}.

\emph{The simple pole.} The right null vector of $\mathcal B_{s_*}$ is $(1,\psi-1)$. Applying $\mathcal B_{s_*}^\prime$ gives the pairing
\[
 \big\langle(\beta,v),\mathcal B'_{s_*}(1,\psi-1)\big\rangle
 =\beta\int R'_{s_*}\psi-\int v\,SR'_{s_*}\psi
 =\int\phi P'_{s_*}\psi=M>0.
\]
The nonzero derivative pairing and the one-dimensional kernel give a simple zero of the determinant.  The pairing $-\int\phi f$ with the right-hand side then gives the rank-one pole \eqref{eq:simple-Green-pole}, and subtracting it leaves an analytic remainder.
\end{proof}

\subsection{Numerical evaluation}\label{subsec:critical-numerics}

We approximate the positive bounded even solution $\psi$ of $L_0\psi=T_{s_*}\psi$, normalized to have limit one at infinity. For a truncation length $R>0$, set $\psi(t)=1$ for $t\ge R$ and write
\[
 \psi(t)=1+\int_t^R(u-t)h(u)\,du,\qquad h=\psi'',\qquad 0\le t\le R.
\]
We extend $\psi$ evenly to $\mathbb R$ and discretize the integral equation $h=\psi/4-T_s\psi$ using composite Gauss--Legendre quadrature with 20 nodes on each subinterval of length two. We evaluate $T_s1$ separately, truncating the half-line integral at $R+70$. For each $s$, we solve the resulting linear system for $h$. The Neumann boundary condition $\psi'(0)=0$, equivalently $\int_0^R h(u)\,du=0$, determines $s$. The results are:
\begin{center}
\begin{tabular}{rll}
\toprule
$R$ & $s$ & $(s^2-1/4)/3$\\
\midrule
24 & 1.599328376133 & 0.769283751568\\
30 & 1.599328376097 & 0.769283751530\\
36 & 1.599328376097 & 0.769283751530\\
\bottomrule
\end{tabular}
\end{center}
Refining the quadrature and extending the integration range leave the displayed digits for $R=36$ unchanged. Thus, we have $s_*\approx1.5993283761$ and $\xi_3\approx0.76928375$ by~\eqref{eq:main-value}.

\section{The candidate exponent for general $\CLE_\kappa$}
\label{sec:general}

The welding calculation in Section~\ref{sec:lqg} gives an explicit integral operator for every $\CLE_\kappa$ with $4<\kappa<8$. Here we specify its boundary condition and propose the resulting candidate for the monochromatic three-arm exponent of $\CLE_\kappa$.

Set
\begin{equation}\label{eq:intro-parameters}
 \alpha=\frac\kappa4\in(1,2),\qquad \gamma=\frac4{\sqrt\kappa},\qquad
 e=\alpha-1,\qquad d=2-\alpha,\qquad \sqrt q=-2\cos\frac\pi\alpha.
\end{equation}
Here $\gamma$ is the LQG parameter and $q\in(0,4)$ is the corresponding Fortuin--Kasteleyn (FK) random-cluster weight. Recall $x(t)=e^t/(1+e^t)$, $a(t)=1-x(t)$ and $L_0=-d^2/dt^2+1/4$. Define
\begin{equation}\label{eq:intro-general-H}
 \begin{gathered}
 C_\alpha=-\frac{\alpha\sin(2\pi/\alpha)}\pi,\\
 H_\alpha(t)=\frac{C_\alpha\alpha}{2}
 \frac{\alpha\cosh(\alpha t/2)\sinh(t/2)-\sinh(\alpha t/2)\cosh(t/2)}{\sinh^3(\alpha t/2)},\qquad
 H_\alpha(0)=\frac{C_\alpha(\alpha^2-1)}{6\alpha},
 \end{gathered}
\end{equation}
with the continuous value at zero, and, whenever the integral converges absolutely,
\begin{equation}\label{eq:intro-general-T}
 (T_{\alpha,s}f)(t)=\int_{\mathbb R}H_\alpha(t-u)
 \bigl[x(t)^{1+s}x(u)^{1-s}+a(t)^{1+s}a(u)^{1-s}\bigr]f(u)\,du.
\end{equation}
The powers use the real logarithms of $x$ and $a$. At $\alpha=3/2$, $H_\alpha$ and $T_{\alpha,s}$ reduce to the kernel and operator of Theorem~\ref{thm:main}.

Let $s_A(\alpha)\in(\alpha-1,2\alpha-1)$ be the parameter for which
\begin{equation}\label{eq:intro-general-boundary-problem}
 L_0\psi=T_{\alpha,s}\psi
\end{equation}
has a strictly positive even solution $\psi\in C^2(\mathbb R)$ satisfying, for some $c>0$,
\begin{equation}\label{eq:intro-general-A-tail}
 \psi(t)=\begin{cases}
 c e^{-(\alpha-3/2)|t|}+o\!\left(e^{-|\alpha-3/2||t|}\right),&\alpha\ne3/2,\\
 c+o(1),&\alpha=3/2,
 \end{cases}\qquad |t|\longrightarrow\infty.
\end{equation}
One can prove existence and uniqueness of $s_A(\alpha)$. The candidate exponent is
\begin{equation}\label{eq:intro-candidate}
 \displaystyle X_A(\kappa):=X_{\kappa/4}\bigl(s_A(\kappa/4)\bigr),\qquad
 X_\alpha(s)=\frac{s^2-(\alpha-1)^2}{2\alpha}.
\end{equation}
At $\kappa=6$, Theorem~\ref{thm:main} gives $s_A(3/2)=s_*$ and $X_A(6)=\xi_3$. One can further show that $X_A(\kappa)\to2$ as $\kappa\to8$, and $X_A(\kappa)\to1/8$ as $\kappa\to4$. Numerical evaluation gives
\begin{equation}\label{eq:intro-numerics}
\begin{array}{c|c|c|c}
 q&\kappa&s_A(\kappa/4)&X_A(\kappa)\\ \hline
 0.5&6.504&1.92948&1.0243\\
 1&6&1.59933&0.76928\\
 2&16/3&1.18706&0.48675\\
 3&24/5&0.88176&0.30729
\end{array}
\end{equation}

\bigskip

\noindent
\textsc{Beijing International Center for Mathematical Research, Peking University,\\
No.5 Yiheyuan Road, Haidian District, Beijing 100871, P.R.China}\\
\textit{Email address:} \texttt{caigefei1917@pku.edu.cn}

\footnotesize
\begin{thebibliography}{ACSW26}
\bibitem[ADA99]{ADA} M.~Aizenman, B.~Duplantier, and A.~Aharony, \emph{Path-crossing exponents and the external perimeter in 2D percolation}, Phys. Rev. Lett. \textbf{83} (1999), no.~7, 1359--1362.

\bibitem[AMY25]{AMYMetrics} V.~Ambrosio, J.~Miller, and Y.~Yuan, \emph{Tightness of approximations to metrics on non-simple conformal loop ensemble gaskets}, preprint, \href{https://arxiv.org/abs/2507.15589v1}{arXiv:2507.15589v1} (2025).

\bibitem[ACSW24]{ACSWLoop} M.~Ang, G.~Cai, X.~Sun, and B.~Wu, \emph{SLE loop measure and Liouville quantum gravity}, preprint, \href{https://arxiv.org/abs/2409.16547v2}{arXiv:2409.16547v2} (2024).

\bibitem[ACSW26]{ACSW} M.~Ang, G.~Cai, X.~Sun, and B.~Wu, \emph{Integrability of conformal loop ensemble: imaginary DOZZ formula and beyond}, J. Amer. Math. Soc. (2026), published electronically October~19, 2026, \href{https://doi.org/10.1090/jams/1086}{DOI 10.1090/jams/1086}.

\bibitem[AHS24]{AHSInt} M.~Ang, N.~Holden, and X.~Sun, \emph{Integrability of SLE via conformal welding of random surfaces}, Comm. Pure Appl. Math. \textbf{77} (2024), no.~5, 2651--2707.

\bibitem[AHSY26]{AHSY} M.~Ang, N.~Holden, X.~Sun, and P.~Yu, \emph{Conformal welding of quantum disks and multiple SLE: the non-simple case}, Probab. Theory Related Fields \textbf{194} (2026), no.~1--2, 717--777.

\bibitem[ARS25]{ARSAnnulus} M.~Ang, G.~R\'emy, and X.~Sun, \emph{The moduli of annuli in random conformal geometry}, Ann. Sci. \'Ec. Norm. Sup\'er. (4) \textbf{58} (2025), no.~4, 1037--1087.

\bibitem[ASYZ25]{ASYZ} M.~Ang, X.~Sun, P.~Yu, and Z.~Zhuang, \emph{Boundary touching probability and nested-path exponent for nonsimple CLE}, Ann. Probab. \textbf{53} (2025), no.~3, 797--847.

\bibitem[BN11]{BN} V.~Beffara and P.~Nolin, \emph{On monochromatic arm exponents for 2D critical percolation}, Ann. Probab. \textbf{39} (2011), no.~4, 1286--1304.

\bibitem[BM26]{BettinelliMiermontII} J.~Bettinelli and G.~Miermont, \emph{Compact Brownian surfaces II. Orientable surfaces}, Mem. Eur. Math. Soc., vol.~27, EMS Press, Berlin, 2026.

\bibitem[Cai26]{Cai26} G.~Cai, \emph{Supplementary note on monochromatic $n\ge4$ arm exponents for critical planar percolation}.\newline
\url{https://gefei-cai.github.io/files/Supplementary_note_monochromatic.pdf}

\bibitem[CFSX25]{CFSX} G.~Cai, X.~Fu, X.~Sun, and Z.~Xie, \emph{Disconnection probability of Brownian motion on an annulus}, preprint, \href{https://arxiv.org/abs/2509.14073v1}{arXiv:2509.14073v1} (2025).

\bibitem[CN06]{CNFull} F.~Camia and C.~M.~Newman, \emph{Two-dimensional critical percolation: the full scaling limit}, Comm. Math. Phys. \textbf{268} (2006), no.~1, 1--38.

\bibitem[DMS21]{DMS21} B.~Duplantier, J.~Miller, and S.~Sheffield, \emph{Liouville quantum gravity as a mating of trees}, Ast\'erisque No.~427 (2021), viii+258~pp.

\bibitem[GPS13]{GPS} C.~Garban, G.~Pete, and O.~Schramm, \emph{Pivotal, cluster, and interface measures for critical planar percolation}, J. Amer. Math. Soc. \textbf{26} (2013), no.~4, 939--1024.

\bibitem[GM21]{GMMetric} E.~Gwynne and J.~Miller, \emph{Existence and uniqueness of the Liouville quantum gravity metric for $\gamma\in(0,2)$}, Invent. Math. \textbf{223} (2021), no.~1, 213--333.

\bibitem[JZJ02]{JZJ} J.~L.~Jacobsen and P.~Zinn-Justin, \emph{Monochromatic path crossing exponents and graph connectivity in two-dimensional percolation}, Phys. Rev. E \textbf{66} (2002), no.~5, 055102.

\bibitem[LSW02]{LSWOneArm} G.~F.~Lawler, O.~Schramm, and W.~Werner, \emph{One-arm exponent for critical 2D percolation}, Electron. J. Probab. \textbf{7} (2002), no.~2, 1--13.

\bibitem[LM25]{LeGallMetzDonnadieu} J.-F.~Le Gall and A.~Metz-Donnadieu, \emph{Drilling holes in the Brownian disk: The Brownian annulus}, Electron. J. Probab. \textbf{30} (2025), Paper No.~33, 43~pp.

\bibitem[LSYZ26]{LSYZ} H.~Liu, X.~Sun, P.~Yu, and Z.~Zhuang, \emph{The bulk one-arm exponent for the $\mathrm{CLE}_{\kappa'}$ percolations}, Proc. Lond. Math. Soc. (3) \textbf{133} (2026), no.~1, e70176.

\bibitem[MSW21]{MSW21} J.~Miller, S.~Sheffield, and W.~Werner, \emph{Non-simple conformal loop ensembles on Liouville quantum gravity and the law of CLE percolation interfaces}, Probab. Theory Related Fields \textbf{181} (2021), no.~1--3, 669--710.

\bibitem[MSW22]{MSW20} J.~Miller, S.~Sheffield, and W.~Werner, \emph{Simple conformal loop ensembles on Liouville quantum gravity}, Ann. Probab. \textbf{50} (2022), no.~3, 905--949.

\bibitem[MW26]{MW26} J.~Miller and W.~Werner, \emph{On the direct SLE-martingale approach to the integrability of conformal loop ensembles}, preprint, \href{https://arxiv.org/abs/2610.31648v1}{arXiv:2610.31648v1} (2026).

\bibitem[MY25]{MYBM} J.~Miller and Y.~Yuan, \emph{Existence and uniqueness of the canonical Brownian motion in non-simple conformal loop ensemble gaskets}, preprint, \href{https://arxiv.org/abs/2512.04807v2}{arXiv:2512.04807v2} (2025).

\bibitem[NQSZ23]{NQSZ} P.~Nolin, W.~Qian, X.~Sun, and Z.~Zhuang, \emph{Backbone exponent for two-dimensional percolation}, preprint, \href{https://arxiv.org/abs/2309.05050v2}{arXiv:2309.05050v2} (2023).

\bibitem[Rem20]{Remy} G.~R\'emy, \emph{The Fyodorov--Bouchaud formula and Liouville conformal field theory}, Duke Math. J. \textbf{169} (2020), no.~1, 177--211.

\bibitem[RV10]{RobertVargas} R.~Robert and V.~Vargas, \emph{Gaussian multiplicative chaos revisited}, Ann. Probab. \textbf{38} (2010), no.~2, 605--631.

\bibitem[Sch00]{SchrammSLE} O.~Schramm, \emph{Scaling limits of loop-erased random walks and uniform spanning trees}, Israel J. Math. \textbf{118} (2000), 221--288.

\bibitem[She09]{Sheffield09} S.~Sheffield, \emph{Exploration trees and conformal loop ensembles}, Duke Math. J. \textbf{147} (2009), no.~1, 79--129.

\bibitem[She16]{Sheffield16} S.~Sheffield, \emph{Conformal weldings of random surfaces: SLE and the quantum gravity zipper}, Ann. Probab. \textbf{44} (2016), no.~5, 3474--3545.

\bibitem[Smi01]{Smirnov} S.~Smirnov, \emph{Critical percolation in the plane: conformal invariance, Cardy's formula, scaling limits}, C. R. Acad. Sci. Paris S\'er. I Math. \textbf{333} (2001), no.~3, 239--244.

\bibitem[SW01]{SWArms} S.~Smirnov and W.~Werner, \emph{Critical exponents for two-dimensional percolation}, Math. Res. Lett. \textbf{8} (2001), 729--744.

\bibitem[SXZ24]{SXZ} X.~Sun, S.~Xu, and Z.~Zhuang, \emph{Annulus crossing formulae for critical planar percolation}, preprint, \href{https://arxiv.org/abs/2410.04767v2}{arXiv:2410.04767v2} (2024).

\bibitem[WDL26]{WDL26} Q.~Wang, Y.~Deng, and M.~Li, \emph{Critical geometry of triconnected clusters in percolation}, preprint, \href{https://arxiv.org/abs/2611.13783}{arXiv:2611.13783} (2026).
\end{thebibliography}
\end{document}